%% file: transform.tex
\documentclass[12pt] {article}
\usepackage{amsmath,amsthm}
\usepackage{amssymb,latexsym}
\usepackage{amscd}
\usepackage{epic,eepic}
\title{A Fourier type transform on translation invariant valuations on convex sets.}
\date{}
\author{ Semyon Alesker \footnote{Partially supported by ISF grant 1369/04.}
\\  { \normalsize Department of Mathematics, Tel Aviv University, Ramat Aviv}
 \\  { \normalsize 69978 Tel Aviv,
Israel }
\\ {\normalsize e-mail: semyon@post.tau.ac.il}}

\def\eps{\varepsilon}
\def\alp{\alpha}
\def\ome{\omega}
\def\Ome{\Omega}
\def\lam{\lambda}

\def\to{\rightarrow}
\def\qed { Q.E.D. }
\def\pt{\partial}

\def\grc{{}\!^ {\textbf{C}} Gr}
\def\ffc{{}\!^ {\textbf{C}} \cf}

\def\RR{\mathbb{R}}
\def\CC{\mathbb{C}}

\def\NN{\mathbb{N}}

\def\DD{\mathbb{D}}
\def\PP{\mathbb{P}}
\def\FF{\mathbb{F}}
\def\LL{\mathbb{L}}

\swapnumbers
\newtheorem{theorem}{Theorem}[subsection]
\newtheorem{corollary}[theorem]{Corollary}
\newtheorem{lemma}[theorem]{Lemma}
\newtheorem{proposition}[theorem]{Proposition}
\newtheorem{claim}[theorem]{Claim}
\theoremstyle{definition}
\newtheorem{example}[theorem]{Example}
\newtheorem{definition}[theorem]{Definition}
\newtheorem{remark}[theorem]{Remark}

\theoremstyle{proposition-definition}
\newtheorem{proposition-definition}[theorem]{Proposition-Definition}

\numberwithin{equation}{subsection}

\def\cf{{\cal F}}

\def\ca{{\cal A}} \def\cb{{\cal B}} \def\cc{{\cal C}}
\def\cd{{\cal D}} \def\ce{{\cal E}} \def\cf{{\cal F}}
\def\cg{{\cal G}} \def\ch{{\cal H}} 
 \def\ck{{\cal K}} \def\cl{{\cal L}}
\def\cm{{\cal M}} \def\cn{{\cal N}} \def\co{{\cal O}}
 
\def\car{{\cal R}}
\def\cs{{\cal S}} \def\ct{{\cal T}} 
\def\cv{{\cal V}}  \def\cx{{\cal X}}
\def\cy{{\cal Y}} \def\cz{{\cal Z}}

\def\pvd{PVal_d(V)}
\def\ond{\Omega^{n}_{d}(V)}

\def\pl{\PP_+(V^*)}

\def\inj{\hookrightarrow }
\def\surj{\twoheadrightarrow}

\newcommand \supp{\operatorname{supp} \,}

\def\grc{{}\!^ {\textbf{C}} Gr}
\def\ffc{{}\!^ {\textbf{C}} \cf}
\def\vc{{}\!^ {\textbf{C}} V}
\input diagram
\begin{document}
\maketitle
\begin{abstract}
Let $V$ be a finite dimensional real vector space. Let
$Val^{sm}(V)$ be the space of translation invariant smooth
valuations on convex compact subsets of $V$. Let $Dens(V)$ be the
space of Lebesgue measures on $V$. The goal of the article is to
construct and study an isomorphism $\FF_V\colon
Val^{sm}(V)\tilde\to Val^{sm}(V^*)\otimes Dens(V)$ such that
$\FF_V$ commutes with the natural action of the full linear group
on both spaces, sends the product on the source (introduced in
\cite{alesker-gafa-04}) to the convolution on the target
(introduced in \cite{bernig-fu}), and satisfies a Plancherel type
formula. As an application, a version of the hard Lefschetz
theorem for valuations is proved.
\end{abstract}
\tableofcontents \setcounter{section}{-1}
\section{Introduction.}\label{S:introduction}
\subsection{An overview of the main results.}\label{Ss:overview}
Let $V$ be a finite dimensional real vector space, $\dim V=n$. The
goal of the article is to construct an isomorphism between the
space of translation invariant valuations on convex compact
subsets of $V$ and the space of translation invariant valuations
(twisted by the line of Lebesgue measures) on the dual space
$V^*$. This isomorphism is analogous to the classical Fourier
transform. It has various nice properties studied in detail in
this article. As an application we prove a version of the hard
Lefschetz theorem for translation invariant valuations. To state
the main results more precisely let us fix some notation and
remind definitions.

Let $\ck(V)$ denote the class of all convex compact subsets of
$V$. Equipped with the Hausdorff metric, the space $\ck(V)$ is a
locally compact space.

\begin{definition}
a) A function $\phi :{\cal K}(V) \to \CC$ is called a valuation if
for any $K_1, \, K_2 \in {\cal K}(V)$ such that their union is
also convex one has
$$\phi(K_1 \cup K_2)= \phi(K_1) +\phi(K_2) -\phi(K_1 \cap K_2).$$

b) A valuation $\phi$ is called continuous if it is continuous
with respect to the Hausdorff metric on ${\cal K}(V)$.
\end{definition}
The notion of valuation is very classical in convexity. For the
classical theory of valuations we refer to the surveys by
McMullen-Schneider \cite{mcmullen-schneider} and McMullen
\cite{mcmullen-survey}. For the general background from convexity
we refer to the book by Schneider \cite{schneider-book}.
Approximately during the last decade there was a considerable
progress in the valuation theory. New classification results of
special classes of valuations have been obtained \cite{klain},
\cite{schneider-simple},
\cite{alesker-annals-99}-\cite{alesker-jdg-03},
\cite{alesker-su2}, \cite{ludwig-reitzner}. Also new structures on
valuations have been discovered \cite{alesker-jdg-03},
\cite{alesker-gafa-04}, \cite{bernig-fu}. Moreover some parts of
the classical theory of valuations on affine spaces have been
generalized to more general context of arbitrary manifolds
\cite{part1}-\cite{part4}, \cite{part3} (see also a survey
\cite{alesker-survey} of these results).

Let us denote by $Val(V)$ the space of translation invariant
continuous valuations. Equipped with the topology of uniform
convergence on compact subsets of $\ck(V)$, $Val(V)$ is known to
be a Banach space. In \cite{alesker-gafa-04} there was introduced
a dense subspace of smooth valuations $Val^{sm}(V)\subset Val(V)$.
The definition is recalled in Section \ref{Ss:product-convolution}
below. Note that $Val^{sm}(V)$ is equipped with the natural
Fr\'echet topology which is stronger than the topology induced
from $Val(V)$.
\begin{example}
1) A Lebesgue measure $vol$ on $V$ belongs to $Val^{sm}(V)$.

2) The Euler characteristic $\chi$ belongs to $Val^{sm}(V)$.
(Recall that $\chi(K)=1$ for any $K\in \ck(V)$.)

3) Let us fix a compact strictly convex set $A\subset V$ with
infinitely smooth boundary. The functional $K\mapsto vol(K+A)$ is a
smooth translation invariant valuation. Here $K+A$ is the Minkowski
sum $\{k+a|\,k\in K,a\in A\}$.

4) Let $0\leq i<n=\dim V$. Fix $A_1,\dots,A_i$ compact strictly
convex subsets of $V$ with infinitely smooth boundary. Then the
mixed volume $K\mapsto V(K[n-i],A_1,\dots,A_i)$ belongs to
$Val^{sm}(V)$ (here $K[n-i]$ means that the set $K$ is taken $n-i$
times). For the notion of mixed volume see e.g.
\cite{schneider-book}, especially Ch. 5,6.
\end{example}
The space $Val^{sm}(V)$ carries a canonical structure of
commutative associative topological algebra with unit (the unit is
the Euler characteristic). This structure was constructed by the
author in \cite{alesker-gafa-04}; the main properties of it are
recalled in Section \ref{Ss:product-convolution}.

Let us denote by $Dens(V)$ the complex one dimensional space of
complex valued Lebesgue measures on $V$. The space
$Val^{sm}(V^*)\otimes Dens(V)$ carries a canonical structure of
commutative associative topological algebra with unit. This
structure was recently constructed by Bernig and Fu
\cite{bernig-fu};  the main properties of it are recalled in
Section \ref{Ss:product-convolution}.

Next observe that group $GL(V)$ of all invertible linear
transformations of $V$ acts naturally on $Val(V)$ (and
$Val^{sm}(V)$) as follows: $(g\phi)(K)=\phi(g^{-1}K)$ for any
$g\in GL(V),K\in \ck(V),$ and $\phi\in Val(V)$ or $Val^{sm}(V)$.

Our first main result says that these two topological algebras
with actions of $GL(V)$ are isomorphic. More precisely we prove
the following result.
\begin{theorem}\label{T:main1}
There exists an isomorphism of linear topological spaces
$$\FF_V\colon Val^{sm}(V)\to Val^{sm}(V^*)\otimes Dens(V)$$ which
satisfies the following properties:

1) $\FF_V$ commutes with the natural action of the group $GL(V)$
on both spaces;

2) $\FF_V$ is an isomorphism of algebras when the source is
equipped with the product and the target with the convolution.

3)(Plancherel type formula) Consider the composition $\ce_V$
$$Val^{sm}(V)\overset{\FF_V}{\to} Val^{sm}(V^*)\otimes
Dens(V)\overset{\FF_{V^*}\otimes
Id_{Dens(V)}}{\to}Val^{sm}(V)\otimes Dens(V^*)\otimes
Dens(V)=Val^{sm}(V).$$

This composition $\ce_V$ satisfies
$$(\ce_V\phi)(K)=\phi(-K).$$
\end{theorem}
\begin{remark}\label{R:intro1}
1) On even valuations, the operator $\FF_V$ was first introduced
by the author in \cite{alesker-jdg-03} under a different name and
notation (in \cite{alesker-jdg-03} it was denoted by $\DD$).

2) Part (2) of Theorem \ref{T:main1} was first proved for
{\itshape even} valuations by Bernig and Fu \cite{bernig-fu}.

3) The isomorphism $\FF_V$ from Theorem \ref{T:main1} is not quite
canonical. One can show that in certain precise sense there exist
exactly four different isomorphisms satisfying the theorem
provided $n>1$; for $n=1$ there exist exactly two such
isomorphisms (see Remark \ref{R:non-unique} below for a precise
statement).

\end{remark}

\begin{example}\label{E:dim1-2}
Let us describe the isomorphism $\FF_V$ in dimensions 1 and 2.
First assume that $\dim V=1$. In this case the space of valuations
is two dimensional: $Val(V)=Val^{sm}(V)=\CC\cdot \chi\oplus
\CC\cdot vol_V$ where $vol_V$ is a non-zero Lebesgue measure on
$V$. Let $vol_V^{-1}$ be the corresponding Lebesgue measure on
$V^*$ (see (\ref{D:inverse-measure}) below). Then
$Val^{sm}(V^*)\otimes Dens(V)=\CC \cdot (vol^{-1}_V\otimes
vol_V)\oplus \CC\cdot (\chi\otimes vol(V))$. Then
\begin{eqnarray}\label{E:euler-measure}
\FF_V(\chi)=vol^{-1}_V\otimes vol_V\\
\FF_V(vol_V)=\chi\otimes vol_V.
\end{eqnarray}

Let us assume now that $\dim V=2$. Let us fix a Euclidean metric
on $V$. It induces identifications $V^*\simeq V,\, Dens(V)\simeq
\CC$. Under these identifications $\FF_V\colon
Val^{sm}(V)\tilde\to Val^{sm}(V)$. One has
$Val^{sm}(V)=\CC\chi\oplus Val_1^{sm}(V)\oplus \CC vol_V$ where
$Val_1^{sm}(V)$ denotes the subspace of 1-homogeneous valuations.
Let us fix also an orientation on $V$. Then
\begin{eqnarray*}
\FF_V(\chi)=vol_V,\\
\FF_V(vol_V)=\chi.
\end{eqnarray*}


In order to describe the action of $\FF_V$ on $Val^{sm}_1(V)$ recall
that by Hadwiger's theorem \cite{hadwiger-51} any valuation $\phi\in
Val^{sm}_1(V)$ can be written uniquely in the form
$$\phi(K)=\int_{S^1}h(\ome) dS_1(K,\ome)$$
where $h\colon S^1\to \CC$ is a smooth function which is orthogonal
on $S^1$ to the two dimensional space of linear functionals. Let us
decompose $h$ to the even and odd parts:
$$h=h_++h_-.$$
Let us decompose further the odd part $h_-$ to "holomorphic" and
"anti-holomorphic" parts
$$h_-=h_-^{hol}+h_-^{anti}$$
as follows. Let us decompose $h_-$ to the usual Fourier series on
the circle $S^1$:
$$h_-(\ome)=\sum_{k} \hat h_-(k)e^{ik\ome}.$$
Then by definition
\begin{eqnarray*}
h_-^{hol}(\ome):=\sum_{k>0} \hat h_-(k)e^{ik\ome},\\
h_-^{anti}(\ome):=\sum_{k<0} \hat h_-(k)e^{ik\ome}.
\end{eqnarray*}
Then the Fourier transform of the valuation $\phi$ is equal to
$$(\FF\phi)(K)=\int_{S^1}(h_+(J\ome)+h_-^{hol}(J\ome))dS_1(K,\ome)-\int_{S^1}h_-^{anti}(J\ome)dS_1(K,\ome)$$
where $J$ is the rotation of $\RR^2$ by $\frac{\pi}{2}$
counterclockwise. (Notice the minus sign before the second
integral.) Observe that $\FF$ preserves the class of real valued
even valuations, but for odd real valued valuations this is not
true. This phenomenon also holds in higher dimensions.
\end{example}

\begin{remark}\label{R:intro10}
In dimension higher than 2 the author does not know such a simple
construction of $\FF_V$, especially in the odd case studied in
detail in the article. Note however that the Fourier transform of
a Lebesgue measure and the Euler characteristic can be computed by
the equalities (\ref{E:euler-measure}) in any dimension.

The construction in the odd case is quite involved and uses
various characterization theorems on valuations, in particular
Klain-Schneider characterization of simple valuations
\cite{klain}, \cite{schneider-simple} (see Theorem
\ref{T:klain-schneider} below), Irreducibility Theorem
\cite{alesker-gafa-01} (see Theorem \ref{T:IrrThm} below), and
also some additional representation theoretical computations based
on the Beilinson-Bernstein localization \cite{beilinson-bernstein}
(see Section \ref{Ss:beilinson-bernstein} below) .
\end{remark}

As an application of the Fourier transform (combined with
Bernig-Br\"ocker theorem \cite{bernig-brocker}) we prove a version
of hard Lefschetz theorem for valuations. In order to state it let
us denote by $Val_i^{sm}(V)\subset Val^{sm}(V)$ the subspace of
$i$-homogeneous valuations ($\phi\in Val_i^{sm}(V)$ if and only if
$\phi(\lam K)=\lam^i\phi(K)$ for any $\lam\geq 0, K\in \ck(V)$).
By McMullen's theorem \cite{mcmullen-euler}
$$Val^{sm}(V)=\oplus_{i=0}^nVal_i^{sm}(V).$$
Let us fix a Euclidean metric on $V$. Let us denote by $V_1$ the
first intrinsic volume (see \cite{schneider-book}, p. 210). Here
we just recall that $V_1$ is the only (up to a constant)
continuous isometry invariant 1-homogeneous valuation (this
characterization is due to Hadwiger \cite{hadwiger-book}).
\begin{theorem}[hard Lefschetz theorem]\label{T:hardLef-intro}
Let $0\leq i<n/2$. Then the map
$$Val^{sm}_i(V)\to Val_{n-i}^{sm}(V)$$
given by $\phi\mapsto (V_1)^{n-2i}\cdot \phi$ is an isomorphism.
\end{theorem}
\begin{remark}\label{R:intro2}
1) Theorem \ref{T:hardLef-intro} was proved by the author for
{\itshape even} valuations in \cite{alesker-jdg-03}.

2) There is another version of the hard Lefschetz theorem for
valuations (see Theorem \ref{T:hlold} below). In the even case it
was proved by the author \cite{alesker-jdg-03}, and in the general
one by Bernig and Br\"ocker \cite{bernig-brocker}. Our proof of
Theorem \ref{T:hardLef-intro} (see Section \ref{lefschetz}) uses
in an essential way this result of Bernig-Br\"ocker. Also we use
the fact that the Fourier transform $\FF_V$ establishes an
equivalence of two versions of the hard Lefschetz theorem (Lemma
\ref{L:hlopers} below). This fact was recently observed in the
even case by Bernig and Fu \cite{bernig-fu}; our Lemma
\ref{L:hlopers} is a straightforward generalization of their
observation.
\end{remark}

Another new construction presented in this article is a
construction of pushforward under linear maps of translation
invariant continuous valuations (twisted by the line of Lebesgue
measures). Namely if $f\colon V\to W$ is a linear space of vector
spaces, we define a linear map
$$f_*\colon Val(V)\otimes Dens(V^*)\to Val(W)\otimes Dens(W^*).$$
We refer to Section \ref{pushforward} for the details. Here we
notice that this pushforward map allows to compute the convolution
of valuations in sense of Bernig and Fu \cite{bernig-fu} in two
steps: first one takes the exterior product of valuations in sense
of \cite{alesker-gafa-04}, and then the pushforward under the
addition map $a\colon V\oplus V\to V$.

Another interesting property of the Fourier transform is that it
intertwines the pullback of valuations (obviously defined, see
Section \ref{pullback}) and the pushforward. With an
oversimplification, one can say that the Fourier transform of the
pullback of a valuation is equal to the pushforward of the Fourier
transform. There are however some technical difficulties of making
this statement rigorous due to the fact that the operations of
pullback and pushforward do  not preserve in general the class of
smooth valuations. Nevertheless a rigorous result is possible
though it sounds more technical: see Theorems \ref{T:z3},
\ref{T:z1} below.

\subsection{Organization of the article.}\label{Ss:organization}
In Section \ref{S:RepTheory} we summarize a necessary background
from representation theory. In Section \ref{S:other_background} we
describe necessary facts mostly from valuation theory. These two
sections do not contain new results.

In Section \ref{S:functorial} we introduce operations of pullback
and pushforward on translation invariant valuations. We relate
them to operations of product and convolution (Sections
\ref{relation}, \ref{Ss:homomor-of-push}), and prove a version of
the base change theorem (Section \ref{Ss:base-change}).

In Section \ref{S:isomorphism-val} we prove an isomorphism of
$GL(V)$-modules $Val^{-,sm}_{n-p}(V)$ and
$Val_p^{-,sm}(V^*)\otimes Dens(V)$ (here $Val_i^{-,sm}(V)$ denotes
the space of smooth odd translation invariant $i$-homogeneous
valuations on $V$). The existence of such isomorphism and some
related representation theoretical calculations will be used in
the construction of the Fourier transform in Section
\ref{S:fourier-high-dim}.

In Section \ref{2-dim} we study separately the Fourier transform
on valuations on a two dimensional plane. The two dimensional case
will be used for higher dimensions in Section
\ref{S:fourier-high-dim}.

Section \ref{S:fourier-high-dim} is the main one. Here we
construct the Fourier transform in full generality and prove the
main properties of it.

In Section \ref{lefschetz} a hard Lefschetz type theorem for
valuations is proved.

The appendix contains a slight generalization of the construction of
the exterior product of smooth valuations given in
\cite{alesker-gafa-04}: here we explain how to multiply a smooth
valuation by a continuous one. This generality is necessary in this
article for technical reasons.

\subsection{Notation.}\label{Ss:notation}
$\bullet$ $\ck(V)$ - the family of convex compact subsets of a
vector space $V$.

$\bullet$ $\ck^{sm}(V)$ - the family of strictly convex compact
subsets of a vector space $V$ with infinitely smooth boundary.

$\bullet$ $f^\vee$ - a dual map to a linear map $f$.

$\bullet$ $f\boxtimes g$, $f\times g$, $f\oplus g$ - operations
with linear maps $f$ and $g$, see Section \ref{Ss:linear algebra}.

$\bullet$ $Gr_i(V)$ - the Grassmannian of $i$-dimensional linear
subspaces of a vector space $V$.

$\bullet$ $\cf_{p,p+1}(V)$ - the (real) variety of partial flags
in $V$ $\{(E,F)|\, E\subset F, \dim E=p,\dim F=p+1\}$.

$\bullet$ $\cf$ - the (real) variety of complete flags in a real
vector space.

$\bullet$ $\ffc$ - the (complex) variety of complete flags in a
complex vector space.

$\bullet$ $\ct_{k,V}, \ct_{k,V;i},\ct^0_{k,V},\ct^0_{k,V;i}$ -
certain vector bundles, see Section \ref{Ss:construction-fourier}.

$\bullet$ $|\ome_X|$ - the complex line bundle of densities over a
manifolds $X$.

\section{Background from representation
theory.}\label{S:RepTheory}
\subsection{Some structure theory of reductive
groups.}\label{Ss:Stucture} In this subsection we remind few basic
definitions from the structure theory of real reductive groups.
For simplicity we will do it only in the case of the group
$GL_n(\RR)$. This will suffice for the purposes of this article.

Let $G_0=GL_n(\RR)$. It acts naturally on $\RR^n$. Let us denote
by $H_0\subset G_0$ the subgroup of diagonal invertible matrices
($H_0$ is a Cartan subgroup).
\begin{definition}\label{D:parabolic}
(i) A subgroup $P_0\subset G_0$ is called {\itshape parabolic} if
there exists a partial flag of linear subspaces $0\ne
F_1\subsetneqq F_2\subsetneqq\dots\subsetneqq F_k=\RR^n$ such that
$P_0$ is the stabilizer in $G_0$ of this flag.

(ii) A parabolic subgroup $P_0$ is called {\itshape standard} if
$P_0\supset H_0$.
\end{definition}

\begin{example}\label{E:blocks}
Let us fix a positive integer $k$ and a decomposition $n=n_1+\dots
n_k,\, n_i\in \NN$. Let
$$P_0:=\left\{\left[\begin{array}{c|c|c|c}
                  A_1&*&\dots&*\\\hline
                  0&A_2&\dots&*\\\hline
                  \dots&\dots&\ddots&\dots\\\hline
                  0&0&\dots&A_k
                  \end{array}\right]\big| A_i\in
                  GL_{n_i}(\RR)\right\}$$
be the subgroup of block upper triangular invertible matrices. It
is a standard parabolic.

Note that if one takes $k=1,n_1=n$ one gets $P_0=G_0$. If one
takes $k=n, n_1=\dots =n_n=1$ then $P_0$ is equal to the subgroup
of upper triangular invertible matrices (it is a minimal
parabolic).
\end{example}
\def\pz{P_0}
\def\upz{U_{P_0}}

For a parabolic subgroup $\pz$ let us denote by $U_{\pz}$ its
unipotent radical. Thus $\upz$ is a normal subgroup of $\pz$. More
explicitly, if $\pz$ is the stabilizer of a partial flag $0\ne
F_1\subsetneqq F_2\subsetneqq\dots\subsetneqq F_k=\RR^n$ then
$\upz$ consists of transformations from $\pz$ inducing the
identity map on all consecutive quotients $F_i/F_{i-1},\,
i=0,\dots,k$.
\begin{definition}\label{D:Levi}
Let $\pz$ be a parabolic. An algebraic subgroup $M\subset \pz$ is
called a {\itshape Levi subgroup} if the canonical homomorphism
$\pz\to \pz/\upz$ induces an isomorphism $M\tilde \to \pz/\upz$.
\end{definition}

A Levi subgroup always exists but not unique. In Example
\ref{E:blocks} the unipotent radical is
$$\upz=\left\{\left[\begin{array}{c|c|c|c}
                     Id_{n_1}&*&\dots&*\\\hline
                     0&Id_{n_2}&\dots&*\\\hline
                     \dots&\dots&\ddots&\dots\\\hline
                      0&0&\dots&Id_{n_k}
                     \end{array}\right]\right\}.$$
A Levi subgroup can be chosen to be equal to
\begin{eqnarray}\label{ex01}
   M=\left\{\left[\begin{array}{c|c|c|c}
                  A_1&0&\dots&0\\\hline
                  0&A_2&\dots&0\\\hline
                  \dots&\dots&\ddots&\dots\\\hline
                  0&0&\dots&A_k
                  \end{array}\right]\big| A_i\in
                  GL_{n_i}(\RR)\right\}\simeq GL_{n_1}(\RR)\times\dots \times GL_{n_k}(\RR).
\end{eqnarray}
\def\zem{{}\!^ {0} M}
Let us return back to a general parabolic $\pz$. If $M$ is its
Levi subgroup then $\pz=M\cdot \upz$.



\begin{definition}\label{D:AssParab}
Let $P_0=M_{\pz}\cdot U_{\pz},\, Q_0=M_{Q_0}\cdot U_{Q_0}$ be two
parabolics, where $M_{P_0},\, M_{Q_0}$ are their Levi subgroups.
Then $P_0$ and $Q_0$ are called {\itshape associated} if there
exists $x\in G_0$ such that
$$M_{Q_0}=x^{-1}M_{P_0}x.$$
\end{definition}

\subsection{Admissible and tempered growth representations and a theorem of Casselman-Wallach.}
\label{Ss:casselman-wallach}
\begin{definition}\label{part1-rep-1} Let $\pi$ be a continuous representation of a Lie group $G_0$ in a
Fr\'echet space $F$. A vector $\xi \in F$ is called $G_0$-smooth
if the map $g\mapsto \pi(g)\xi$ is an infinitely differentiable
map from $G_0$ to $F$.
\end{definition}
It is well known (see e.g. \cite{wallach}, Section 1.6) that the
subset $F^{sm}$ of smooth vectors is a $G_0$-invariant linear
subspace dense in $F$. Moreover it has a natural topology of a
Fr\'echet space (which is stronger than the topology induced from
$F$), and the representation of $G_0$ in $F^{sm}$ is continuous.
Moreover all vectors in $F^{sm}$ are $G_0$-smooth.

Let $G_0$ be a real reductive group. Assume that $G_0$ can be
imbedded into the group $GL_N(\RR)$ for some $N$ as a closed
subgroup invariant under the transposition. Let us fix such an
imbedding $p:G_0\hookrightarrow GL_N(\RR)$. (In our applications
$G_0$ will be either $GL_n(\RR)$ or a direct product of several
copies of $GL_n(\RR)$.) Let us introduce a norm $|\cdot |$ on
$G_0$ as follows:

$$|g|:=\max\{||p(g)||,||p(g^{-1})||\}$$
where $||\cdot||$ denotes the usual operator norm in $\RR^N$.
\begin{definition}
Let $\pi$ be a smooth representation of $G_0$ in a Fr\'echet space
$F$ (namely $F^{sm}=F$). One says that this representation has
{\itshape moderate growth} if for each continuous semi-norm
$\lambda$ on $F$ there exists a continuous semi-norm $\nu_\lambda$
on $F$ and $d_{\lambda}\in \RR$ such that
$$\lambda(\pi(g)v)\leq |g|^{d_\lambda}\nu_{\lambda}(v)$$
for all $g\in G,\, v\in F$.
\end{definition}

The proof of the next lemma can be found in \cite{wallach}, Lemmas
11.5.1 and 11.5.2.
\begin{lemma}\label{part1-wallach}
(i) If $(\pi,G_0,H)$ is a continuous representation of $G_0$ in a
Banach space $H$, then $(\pi,G,H^{sm})$ has moderate growth.

(ii) Let $(\pi, G_0,V)$ be a representation of moderate growth.
Let $W$ be a closed $G_0$-invariant subspace of $V$. Then $W$ and
$V/W$ have moderate growth.
\end{lemma}

Remind that a continuous Fr\'echet representation $(\pi,G_0,\cf)$
is said to have {\itshape finite length} if there exists a finite
filtration
$$0=F_0\subset F_1\subset \dots\subset F_m=F$$
by $G_0$-invariant closed subspaces such that $F_i/F_{i-1}$ is
irreducible, i.e. does not have proper closed $G_0$-invariant
subspaces. The sequence $F_1,F_2/F_1,\dots,F_m/F_{m-1}$ is called
the Jordan-H\"older series of the representation $\pi$. It is well
known (and easy to see) that the Jordan-H\"older series of a
finite length representation is unique up to a permutation.

\begin{definition}
A Fr\'echet representation $(\rho,G_0,F)$ of a real reductive
group $G_0$ is called {\itshape admissible} if its restriction to
a maximal compact subgroup $K$ of $G_0$ contains an isomorphism
class of any irreducible representation of $K$ with at most finite
multiplicity. (Remind that a maximal compact subgroup of
$GL_n(\RR)$ is the orthogonal group $O(n)$.)
\end{definition}

\begin{theorem}[Casselman-Wallach, \cite{casselman}]\label{casselman-wallach}
Let $G_0$ be a real reductive group. Let $(\rho,G_0,F_1)$ and
$(\pi,G_0,F_2)$ be smooth representations of moderate growth in
Fr\'echet spaces $F_1, F_2$. Assume in addition that $F_2$ is
admissible of finite length. Then any continuous morphism of
$G_0$-modules $f:F_1\to F_2$ has closed image.
\end{theorem}

The following proposition is essentially a common knowledge; the
proof can be found in \cite{part1}, Proposition 1.1.8.
\begin{proposition}\label{part1-epi}
Let $G_0$ be a real reductive Lie group. Let $F_1,F_2$ be
continuous Fr\'echet $G_0$-modules. Let $\xi:F_1\to F_2$ be a
continuous morphism of $G_0$-modules. Assume that the assumptions
of the Casselman-Wallach theorem are satisfied, namely $F_1$ and
$F_2$ are smooth and have moderate growth, and $F_2$ is admissible
of finite length. Assume moreover that $\xi$ is surjective.

Let $X$ be a smooth manifold. Consider the map
$$\hat \xi:C^\infty(X,F_1)\to C^\infty(X,F_2)$$
defined by $(\hat\xi (f))(x)=\xi(f(x))$ for any $x\in X$.

Then $\hat\xi$ is surjective.
\end{proposition}

\subsection{Induced representations.}\label{Ss:InducedRep}
Let $H\subset G_0$ be a closed subgroup. Let $\pi$ be a
representation of $H$ in a Fr\'echet space $F$. Let us consider
the space of continuous functions
$$\Phi:=\{f\colon G_0\to V|\,\, f(x\cdot h)=\pi(h)^{-1}(f(x))\mbox{
for any } x\in G_0,h\in H\}.$$ The group $G_0$ acts on $\Phi$ by
left translation. The representation of $G_0$ in $\Phi$ is called
the induced representation and denoted by $Ind_H^{G_0}\pi$. Note
that the space $\Phi$ is a space of continuous sections of a
$G_0$-equivariant vector bundle over $G_0/H$ with fiber $F$.

Let $P_0\subset G_0$ be a parabolic subgroup. Let us consider the
natural representation of $G_0$ in the space of (complex valued)
half-densities (see e.g \cite{guillemin-sternberg}, Ch. II \S 6)
on $G_0/P_0$. It is easy to see that there exists a character
$\rho_{\pz}\colon P_0\to \CC^*$ such that this representation is
isomorphic to $Ind_{P_0}^{G_0}\rho_{\pz}$.

\begin{remark}\label{R:character}
It is easy to see that $\upz\subset [P_0,P_0]$. Hence any
character of $P_0$ is trivial on $\upz$, and hence factorizes via
$\pz/\upz$.
\end{remark}

 Let $\pi$ be a
representation of $P_0/\upz$ considered as a representation of
$P_0$. Let us denote by
$$\ch(P_0,\pi):=Ind_{\pz}^{G_0}(\pi\otimes \rho_{\pz}).$$
We will need the next result which is standard in representation
theory.
\begin{theorem}\label{T:finiteness}
Let $P_0\subset G_0$ be a parabolic subgroup. Let $\pi$ be a
character of $P_0/\upz$. Then $Ind_{\pz}^{G_0}(\pi)$ is an
admissible representation of finite length.
\end{theorem}

Let now $P_0=M_{\pz}\cdot U_{\pz},\, Q_0=M_{Q_0}\cdot U_{Q_0}$ be
two associated parabolics (see Definition \ref{D:AssParab}). Thus
$M_{Q_0}=\alp^{-1}M_{P_0}\alp$ for some $\alp\in G_0$. Let $\phi$
be a character of $M_{P_0}$ considered as a character of $P_0$.
Let $\pi'$ be the character of $M_{Q_0}$ defined by
$$\pi'(x)=\pi(\alp\cdot x\cdot \alp^{-1}).$$
The following result is a special case of a theorem by
Harish-Chandra (see e.g. \cite{green-monster}, Proposition
4.1.20).
\begin{theorem}[Harish-Chandra]\label{T:AssHCh}
Let $P_0,Q_0\subset G_0$ be two associated parabolics as above.
Let $\pi$ be a character of $P_0/\upz$. Let $\pi'$ be the
character of $Q_0/U_{Q_0}$ as above. The representations
$\ch(P_0,\pi)$ and $\ch(Q_0,\pi')$ have the same Jordan-H\"older
series.
\end{theorem}
We will use below Theorem \ref{T:AssHCh} in two particular
situations. To describe the first one, let us fix an integer
$k=1,2,\dots n-1$. Let
$$P_0:=\left\{\left[\begin{array}{c|c|c}
                      A&*&*\\\hline
                      0&B&*\\\hline
                      0&0&c
                      \end{array}\right]\big| A\in GL_{n-k-1}(\RR), B\in
                      GL_k(\RR), c\in \RR^*\right\}.$$

Then $G_0/P_0$ is the partial flag space $\{(F,E)|\, F\in
Gr_{n-k-1}(\RR^n), E\in Gr_{n-1}(\RR^n)\}$. Let $p\colon
G_0/P_0\to Gr_{n-k-1}(\RR^n)$ be the natural map given by
$p(F,E)=F$. Let $\cv\to G_0/P_0$ be the line bundle whose fiber
over a partial flag $(F,E)$ is equal to $Dens(E/F)\otimes
or(\RR^n/E)$. Let $|\ome_{Gr_{n-k-1}}|\to Gr_{n-k-1}(\RR^n)$
denote the line bundle of densities.  Let $\ca$ denote the space
of sections $\ca:=C(G_0/P_0,\cv\otimes p^*(|\ome_{Gr_{n-k-1}}|))$.
Clearly $\ca$ is a $G_0$-module.

\begin{lemma}\label{L:normaliz1}
$$\ca =\ch(P_0,\pi)$$ where $\pi\colon P_0\to \CC^*$ is the
character defined by
\begin{eqnarray}\label{def-pi}\pi\left(\left[\begin{array}{c|c|c}
                  A&*&*\\\hline
                      0&B&*\\\hline
                      0&0&c
                      \end{array}\right]\right)=
                      |\det A|^{\frac{k+1}{2}}|\det
                      B|^{-\frac{n-k}{2}-1} |c|^{-\frac{n-1}{2}+k}
                      sgn (c).\end{eqnarray}
\end{lemma}
{\bf Proof} is straightforward computation. \qed

Let us consider another parabolic subgroup
$$Q_0:=\left\{\left[\begin{array}{c|c|c}
                      B&*&*\\\hline
                      0&c&*\\\hline
                      0&0&A
                      \end{array}\right]\big| A\in GL_{n-k-1}(\RR), B\in
                      GL_k(\RR), c\in \RR^*\right\}.$$
It is easy to see that $P_0$ and $Q_0$ are associated. Let
$\xi\colon Q_0\to \CC^*$ be the character defined by
$$\xi\left(\left[\begin{array}{c|c|c}
                      B&*&*\\\hline
                      0&c&*\\\hline
                      0&0&A
                      \end{array}\right]\right)=|\det B|^{-1}sgn (c).$$

\begin{lemma}\label{L:normaliz2}
$$Ind_{Q_0}^{G_0}\xi=\ch(Q_0,\rho)$$
where $\rho\left(\left[\begin{array}{c|c|c}
                      B&*&*\\\hline
                      0&c&*\\\hline
                      0&0&A
                      \end{array}\right]\right)=
|\det A|^{\frac{k+1}{2}}|\det
                      B|^{-\frac{n-k}{2}-1} |c|^{-\frac{n-1}{2}+k}
                      sgn (c).$
\end{lemma}
{\bf Proof} is straightforward computation. \qed

\begin{corollary}\label{C:arepr}
The $G_0$-modules $\ca^{sm}$ and $(Ind_{Q_0}^{G_0}\xi)^{sm}$ have
the same Jordan-H\"older series.
\end{corollary}
{\bf Proof.} Obviously $P_0$ and $Q_0$ are associated, indeed
$M_{P_0}=\alp\cdot M_{Q_0}\cdot \alp^{-1}$ where
$\alp=\left[\begin{array}{cc}
          0&Id_{n-k}\\
          Id_k&0\\
          \end{array}\right]$. Moreover
$\xi(x)=\pi(a^{-1}xa)$ where $\pi$ is defined by (\ref{def-pi}).
Hence the result follows from Theorem \ref{T:AssHCh}. \qed

Let us describe now the second situation where Theorem
\ref{T:AssHCh} will be used.
\begin{corollary}\label{C:arepr2}
Let us consider a parabolic subgroup
\begin{eqnarray*}
P_0'=\left\{\left[\begin{array}{c|c|c}
                    A&*&*\\\hline
                    0&c&*\\\hline
                    0&0&B
                    \end{array}\right]\big| A\in GL_p(\RR),c\in
                    \RR^*,B\in GL_{n-p-1}(\RR)\right\}.
\end{eqnarray*}
Let us denote by $\ome\colon P_0'\to \CC^*$ the character such
that $Ind_{P_0'}^{G_0}(\ome)$ is isomorphic to the representation
of $G_0$ in densities on $G_0/P'_0$. Let $\alp\colon P_0'\to
\CC^*$ be the character defined by
$$\alp\left(\left[\begin{array}{c|c|c}
                    A&*&*\\\hline
                    0&c&*\\\hline
                    0&0&B
                    \end{array}\right]\right)=|\det(A)|\cdot sgn(c)\cdot
                    \ome.$$
Let us consider another parabolic subgroup
\begin{eqnarray*}
Q_0'=\left\{\left[\begin{array}{c|c|c}
                    B&*&*\\\hline
                    0&c&*\\\hline
                    0&0&A
                    \end{array}\right]\big| A\in GL_p(\RR),c\in
                    \RR^*,B\in GL_{n-p-1}(\RR)\right\}.
\end{eqnarray*}
Let $\beta\colon Q_0'\to \CC^*$ be the character defined by
$$\beta\left(\left[\begin{array}{c|c|c}
                    B&*&*\\\hline
                    0&c&*\\\hline
                    0&0&A
                    \end{array}\right]\right)=|\det(A)|\cdot
                    sgn(c).$$
Then $(Ind_{P_0'}^{G_0}(\alp))^{sm}$ and
$(Ind_{Q_0'}^{G_0}(\beta))^{sm}$ have the same Jordan-H\"older
series.
\end{corollary}
{\bf Proof} is a straightforward computation similar to the proof
of Corollary \ref{C:arepr}. \qed

\subsection{The Beilinson-Bernstein
localization.}\label{Ss:beilinson-bernstein}
\def\dla{\cd_\lam}
We recall the Beilinson-Bernstein theorem on localization of
$\textsl{g} $-modules following \cite{beilinson-bernstein} (see
also \cite{bien}). Then we recall the version of this result for
dominant but not regular characters following \cite{kashiwara}. We
denote by capital letters the Lie groups, and by the corresponding
small letters their Lie algebras.

Let $G$ be a complex reductive algebraic group. Let $T$ denote a
Cartan subgroup of $G$. In our examples $G=GL_n(\CC)$.  Let $B$ be
a Borel subgroup of $G$ containing $T$. Let $t$ denote the Lie
algebra of $T$. Let $b$ denote the Lie algebra of $B$. Let $n$
denote the nilpotent radical of $b$.

In the case when $G$ is complexification of a real reductive group
$G_0$ let us denote by $K$ the complexification of a maximal
compact subgroup of $G_0$. Thus if $G=GL_n(\CC)$ then $K=O(n,\CC)$
is the group of complex orthogonal matrices.

Let $R(t)\subset t^*$ be the set of roots of $t$ in $\textsl{g}$.
The set $R(t)$ is naturally divided into the set of roots whose
root spaces are contained in $n$ and its complement. Let $R^+(t)$
be the set of roots of $t$ in $g/b$. If $\alpha$ is a root of $t$
in $\textsl{g}$ then the dimension of the corresponding root
subspace $\textsl{g}_{\alpha}$ is called the multiplicity of
$\alpha$. Let $\rho _b$ be the half sum of the roots contained in
$R^+(t)$ counted with their multiplicities.

\begin{definition}\label{D:regular-dominant}
We say that $\lambda \in t^*$ is {\itshape dominant} if for any
root $\alpha \in R^+(t_b)$ we have $<\lambda, \alpha ^{V}>\ne -1,
-2, \dots$. We shall say that $\lambda \in t_b^*$ is $B$-{\itshape
regular} if for any root $\alpha \in R^+(t_b)$ we have $<\lambda,
\alpha ^{V}>\ne 0$.
\end{definition}

For the definitions and basic properties of the sheaves of twisted
differential operators we refer to \cite{bien}, \cite{kashiwara}.
Here we will present only the explicit description of the sheaf
$\dla$ in order to agree about the normalization.

\def\ug{U(\texttt{g})}
\def\ox{{\cal O}_{\cf}}
\def\tx{{\cal T}_{\cf}}
Let $\ffc$ be the complete flag variety of $G$ (then $\ffc=G/B$).
Let $\ox$ denote the sheaf of regular functions on $\ffc$. Let
$\ug$ denote the universal enveloping algebra of $\textsf{g}$. Let
$U^o$ be the sheaf $\ug \otimes _{\CC} \ox$, and
$\textsf{g}^o:=\textsf{g}\otimes _{\CC}\ox$. Let $\tx$ be the
tangent sheaf of $\ffc$. We have a canonical morphism
$\alpha\colon \textsf{g}^o\to \tx$. Let also $b^o:=Ker \,\alpha
=\{\xi\in\textsf{g}^o|\,\xi_x\in b_x \forall x\in \ffc\}. $ Let
$\lambda\colon b\to \CC$ be a linear functional which is trivial
on $n$ (thus $\lambda \in t^*$). Then $\lambda$ defines a morphism
$\lambda ^o\colon b^o\to \ox$. We will denote by ${\cal
D}_\lambda$ the sheaf of twisted differential operators
corresponding to $\lambda -\rho _b$, i.e. $\dla$ is isomorphic to
$U^o/{\cal I}_{\lambda}$, where ${\cal I}_{\lambda}$ is the two
sided ideal generated by the elements of the form $\xi-(\lambda
-\rho_b)^o(\xi)$ where $\xi$ is a local section of $b^o$.

Let $D_\lambda:=\Gamma (\ffc,\dla)$ denote the ring of global
sections of $\dla$. We have a canonical morphism $U(\textsl{g})\to
D_\lambda$. For the complete flag variety $\ffc$ this map is onto
(\cite{beilinson-bernstein}). The kernel of this homomorphism was
also described in  \cite{beilinson-bernstein}. To describe it,
remind that we have the Harish-Chandra isomorphism
$Z(U(\textsl{g}))\tilde \to (Sym^\bullet(t))^W$ where
$Z(U(\textsl{g}))$ denotes the center of $U(\textsl{g})$,
$(Sym^\bullet(t))^W$ is the algebra of elements of the full
symmetric algebra of $t$ invariant under the Weyl group $W$. Hence
the element $\lam\in t^*$ defines a homomorphism
$Z(U(\textsl{g}))\to \CC$ called the infinitesimal character. Let
$I_\lam$ be the two-sided ideal in $U(\textsl{g})$ generated by
the kernel of this homomorphism. Then by
\cite{beilinson-bernstein}, $I_\lam$ is equal to the kernel of the
homomorphism $U(\textsl{g})\to D_\lam$.
\begin{remark}\label{R:infchar}
The category of $D_\lam$-modules coincides with the category of
$\textsl{g}$-modules with the given infinitesimal character.
\end{remark}

In this notation one has the following result proved in
\cite{beilinson-bernstein}.

\begin{theorem}[Beilinson-Bernstein]\label{T:beilinson-bernstein}
(1) If $\lambda \in t^*$ is dominant then the functor $\Gamma:
{\cal D}_{\lambda}-mod \to U(\textsl{g})-mod$ is exact.

(2) If $\lambda \in t^*$ is dominant and regular then the functor
$\Gamma$ is also faithful.
\end{theorem}

Note also that always the functor $\Gamma$ has a left adjoint
functor (called the localization functor) $\Delta: D_\lambda
-mod\to \dla -mod$. It is defined as $\Delta (M)=\dla
\otimes_{D_\lambda} M$. $\Delta(M)$ is called the {\itshape
localization} of $M$.

The proof of the next lemma can be found in \cite{bien},
Proposition I.6.6.

\begin{lemma}
Suppose $\Gamma: {\cal D}_{\lambda}-mod \to D_{\lambda}-mod$ is
exact. Then the localization functor
$\Delta:D_{\lambda}-mod\to\dla -mod$ is the right inverse of
$\Gamma$:
$$\Gamma\circ \Delta =Id.$$


\end{lemma}
We have the following immediate corollary (see \cite{bien}, p.
24).
\begin{corollary}\label{C:BBequivalence}
Let $\lambda \in t^*$ is dominant and regular. Then the functor
$$\Gamma\colon \cd_\lam-mod\to D_\lam-mod$$
is an equivalence of categories. Moreover this equivalence holds
for $K$-equivariant versions of these categories.
\end{corollary}

The following sufficient condition for $\lam$ being regular and
dominant will be useful.
\begin{proposition}\label{P:reg-dom}
Let $G_0=GL_n(\RR),\, G=GL_n(\CC)$. Let $B_0\subset G_0$ be the
subgroup of real invertible upper triangular matrices, $B$ be its
complexification . Let $\lam \in t^*$. Let $\chi\colon B_0\to
\CC^*$ be a character such that its (complexified) differential
$b\to \CC$ is equal to $\lam -\rho_b$. Assume that the
representation $Ind_{P_0}^{G_0}\chi$ has a non-zero finite
dimensional $G_0$-submodule. Then $\lam$ is dominant and regular.
\end{proposition}
\subsection{Some analysis on
manifolds.}\label{Ss:analysis-on-manifolds} Let $X$ be a smooth
manifold countable at infinity. Let $\ce\to X$ be a finite
dimensional vector bundle. We denote by $C(X,\ce),\,
C^\infty(X,\ce)$ the spaces of continuous, $C^\infty$-smooth
sections respectively. The space $C(X,\ce)$ being equipped with
the topology of uniform convergence on compact subsets of $X$ is a
Fr\'echet space; if $X$ is compact then it is Banach space. The
space $C^\infty(X,\ce)$ being equipped with the topology of
uniform convergence on compact subsets of $X$ of all partial
derivatives is a Fr\'echet space.

The following result is well known (see e.g.
\cite{gelfand-vilenkin}).

\begin{theorem}\label{T:bilinear-forms}
Let $X_1$and $X_2$ be compact smooth manifolds. Let $\ce_1$ and
$\ce_2$ be smooth finite dimensional vector bundles over $X_1$ and
$X_2$ respectively. Let $\cg$ be a Fr\'echet space. Let
$$B:C^\infty(X_1,\ce_1)\times C^\infty(X_2,\ce_2)\to \cg$$
be a continuous bilinear map. Then there exists unique continuous
linear operator
$$b:C^\infty(X_1\times X_2,\ce_1\boxtimes \ce_2)\to \cg$$
such that $b(f_1\otimes f_2)=B(f_1,f_2)$ for any $f_i\in
C^\infty(X_i,\ce_i),\, i=1,2$.
\end{theorem}

\section{Background on valuation theory.}\label{S:other_background} We collect
in this section some necessary notation and results from the
valuation theory and linear algebra.
\subsection{Linear algebra.}\label{Ss:linear algebra} Let $V$ be a
finite dimensional real vector space of dimension $n$. We denote
\begin{eqnarray}\label{D:det}
\det V:=\wedge^nV.
\end{eqnarray}
Also we denote by $Dens(V)$ the space of complex valued Lebesgue
measures on $V$; thus $Dens(V)$ is a complex line. ($Dens$ stays
for {\itshape densities}.)

Let us define $or(V)$ the {\itshape orientation line} of $V$ as
follows. Let us denote by $Bas(V)$ the set of all basis in $V$.
The group $GL(V)$ of linear invertible transformations acts
naturally on $Bas(V)$ by
$g((x_1,\dots,x_n))=(g(x_1),\dots,g(x_n))$. Then set
\begin{eqnarray*}\label{D:or}
or(V)=\{f\colon Bas(V)\to \CC|\, \, f(g(x))=\\sgn(\det(g))\cdot
f(x)\mbox{ for any } g\in GL(V),\, x\in Bas(V)\}.
\end{eqnarray*}
Clearly $or(V)$ is a one dimensional complex vector space. The
operation of passing to the biorthogonal basis gives an
identification $Bas(V)\simeq Bas(V^*)$. It induces an isomorphism
of vector spaces
\begin{eqnarray}\label{or-dual}
or(V)\simeq or(V^*)
\end{eqnarray}
which will be used throughout the article. Also one has another
natural isomorphism
\begin{eqnarray}\label{det-dens}
Dens(V)\simeq \det(V^*)\otimes_\CC or(V).
\end{eqnarray}

Next let
\begin{eqnarray}\label{Exact:Seq}
0\to U\to V\to W\to 0
\end{eqnarray}
be a short exact sequence of finite dimensional vector spaces.
Then one has canonical isomorphisms
\begin{eqnarray}\label{Can:1}
\det U\otimes \det W\tilde\to\det V,\\\label{Can:2} or(U)\otimes
or(W)\tilde \to or(V),\\\label{Can:3} Dens(U)\otimes
Dens(W)\tilde\to Dens(V).
\end{eqnarray}
The isomorphism (\ref{Can:1}) is given by $x\otimes y\mapsto
x\wedge \tilde y$ where $\tilde y$ is an arbitrary lift of $y\in
\det W$ to $\wedge^{\dim W}V$. To describe the isomorphisms
(\ref{Can:2}) and (\ref{Can:3}) let us fix an arbitrary linear
splitting of the exact sequence (\ref{Exact:Seq}), $s:W\to V$.
Thus $V\simeq U\oplus W$. Then $Bas(U)\times Bas(W)\subset
Bas(V)$. Thus the restriction from functions on $Bas(V)$ to
functions on $Bas(U)\times Bas(W)$ defines the isomorphism
(\ref{Can:2}) which is in fact independent of the choice of a
splitting $s$. Next the usual product measure construction defines
a map $Dens(U)\otimes Dens(W)\to Dens(V)$ which is an isomorphism
independent of a choice of a splitting $s$.

Also we will use an isomorphism
\begin{eqnarray}\label{Can:10}
Dens(V^*)\tilde\to Dens(V)^*
\end{eqnarray}
which can be described as follows. Let us fix a basis of $V$, and
let $(x_1,\dots,x_n)$ denote coordinates of a vector in $V$ in
this basis. Let $(y_1,\dots,y_n)$ denote coordinates of a vector
in $V^*$ in the biorthogonal basis. Let us choose the isomorphism
(\ref{Can:10}) such that the Lebesgue measure on $|dy_1\wedge\dots
\wedge dy_n|$ on $V^*$ goes to the linear functional on $Dens(V)$
whose value on the Lebesgue measure $|dx_1\wedge\dots \wedge
dx_n|$ on $V$ is equal to 1. It is easy to see that this
isomorphism does not depend on a choice of a basis. In this
situation we also write
\begin{eqnarray}\label{D:inverse-measure}
|dy_1\wedge\dots \wedge dy_n|=|dx_1\wedge\dots \wedge dx_n|^{-1}.
\end{eqnarray}

For a linear map of vector spaces $f\colon V\to W$ we denote by
$f^\vee\colon W^*\to V^*$ the dual space. (This notation is
different from the probably more standard notation $f^*$, since we
keep the symbol $f^*$ to denote the pullback of valuations, see
Section \ref{pullback}.)

Let $f_i\colon V_i\to W_i,\, i=1,2,$ be two linear maps of vector
spaces. We denote by
\begin{eqnarray}\label{D:boxtimes-lin}
f_1\boxtimes f_2\colon V_1\oplus V_2\to W_1\oplus W_2
\end{eqnarray}
defined, as usual, by $(f_1\boxtimes
f_2)(x_1,x_2)=(f_1(x_1),f_2(x_2))$.

Let $g_i\colon V\to W_i,\, i=1,2,$ be two linear maps. We denote
by
\begin{eqnarray}\label{D:times-lin}
g_1\times g_2\colon V\to W_1\oplus W_2
\end{eqnarray}
the map defined by $(g_1\times g_2)(x)=(g_1(x),g_2(x))$.

Let $h_i\colon V_i\to W,\, i=1,2,$ be two linear maps. We denote
by
\begin{eqnarray}\label{D:oplus-lin}
h_1\oplus h_2\colon V_1\oplus V_2\to W
\end{eqnarray}
the map defined by $(h_1\oplus h_2)(x_1,x_2)=h_1(x_1)+h_2(x_2)$.

\subsection{McMullen's decomposition of
valuations.}\label{Ss:mcmullen-decomp} Let $V$ be an
$n$-dimensional real vector space. Let $Val(V)$ denote the space
of translation invariant continuous valuations on $V$.  Let $\alp$
be a complex number. We say that $\phi$ is $\alp$-homogeneous if
$$\phi(\lam K)=\lam^\alp\phi(K) \mbox{ for any } \lam>0,\, K\in
\ck(V).$$ Let us denote by $Val_\alp(V)$ the subspace of $Val(V)$
of $\alp$-homogeneous convex valuations. The following result is
due to P. McMullen \cite{mcmullen-euler}.
\begin{theorem}[\cite{mcmullen-euler}]\label{T:mcmullen-decomp}
Let $n=\dim V$. Then
$$Val(V)=\bigoplus_{i=0}^n Val_i(V).$$
\end{theorem}

\subsection{Characterization theorems on
valuations.}\label{Ss:charact-homogeneous} In this section we
describe several theorems on translation invariant continuous
valuations which will be used in the article.
\begin{proposition}\label{P:zero-n}
(1) $Val_0(V)$ is spanned by the Euler characteristic $\chi$.

(2)(Hadwiger, \cite{hadwiger-book}) $Val_n(V)$ is spanned by a
Lebesgue measure $vol$.
\end{proposition}
Note that part (1) of the proposition is obvious. Let us remind
now a description of $(n-1)$-homogeneous valuations due to P.
McMullen \cite{mcmullen-80}.

\def\pph{\PP^\vee_+(V)}
Let $\pph$ denote the manifold of cooriented linear hyperplanes in
$V$ (recall that coorientation of $E\subset V$ is just an
orientation of $V/E$). Let $\cl\to \pph$ denote the complex line
bundle whose fiber over $E\in \pph$ is equal to $Dens(E)$. Let us
construct a continuous linear map
\begin{eqnarray}\label{map-mcmullen}
\Psi\colon C(\pph, \cl)\to Val_{n-1}(V).
\end{eqnarray}
For let us fix $\xi\in C(\pph,\cl)$. Let us define first a
valuation $\phi_\xi$ on convex compact {\itshape polytopes} in
$V$. Let $P\subset V$ be a convex compact polytope. Any
$(n-1)$-dimensional face of $P$ carries a coorientation such that
the exterior normal of it has a positive direction. Set
\begin{eqnarray}\label{phixi}
\Psi(\xi)(P):=\sum_{F\in\{(n-1)-\mbox{faces of }P\}}\xi(F).
\end{eqnarray}
\begin{lemma}\label{L:phixi}
$\Psi(\xi)$ extends (uniquely) to a continuous valuation on
$\ck(V)$. This valuation (also denoted by $\Psi(\xi)$) is
continuous, translation invariant, and $(n-1)$-homogeneous.
\end{lemma}
{\bf Proof.} In fact this claim is due to Schneider
\cite{schneider-75} in a slightly different language. We are going
to explain this. Let us fix a Euclidean metric on $V$. It induces
an identification of $\pph$ with the unit sphere $S^{n-1}$, and an
isomorphism of $\cl$ with the trivial line bundle. Then if $\xi\in
C(\pph,\cl)\simeq C(S^{n-1})$ then for a polytope $P$
$$\Psi(\xi)(P)=\int_{S^{n-1}}\xi(\ome)dS_{n-1}(P,\ome)$$
where $dS_{n-1}(P,\bullet)$ dentes the $(n-1)$-th area measure
(see e.g. \cite{schneider-book}, p.203). But the functional on
$\ck(V)$
$$K\mapsto \int_{S^{n-1}}\xi(\ome)dS_{n-1}(K,\ome)$$
is a continuous translation invariant $(n-1)$-homogeneous
valuation by \cite{schneider-75}. \qed

The next theorem is due to McMullen.
\begin{theorem}[\cite{mcmullen-80}]\label{T:mcmullen-80}
The map $\Psi$ is onto. The kernel is $n$-dimensional.
\end{theorem}
\begin{remark}\label{R:invariance}
 The line bundle $\cl$ is obviously $GL(V)$-equivariant. The map
$\Psi$ is $GL(V)$-equivariant.
\end{remark}
Let us describe explicitly the kernel of $\Psi$. First let us
construct a linear map $$\wedge^{n-1}V^*\otimes or(V)\to
C(\pph,\cl).$$ Fix an arbitrary cooriented hyperplane $E\in \pph$.
By (\ref{det-dens}), (\ref{or-dual}), (\ref{Can:2}) one has
$$Dens(E)=\det E^*\otimes or(E)=\det E^*\otimes or(V)\otimes
or(V/E).$$ Since $E$ is cooriented $or(V/E)=\CC$. Hence
$$Dens(E)=\det E^*\otimes or(V).$$
But one has a canonical map $\wedge ^{n-1}V^*\otimes or(V)\to \det
E^*\otimes or(V)$ induced by the map $V^*\to E^*$ dual to the
identity imbedding $E\to V$. When $E$ varies this defines the
desired map $\wedge^{n-1}V^*\otimes or(V)\to C(\pph,\cl)$. Clearly
it is $GL()$-equivariant and injective.

The claim is that the image of the above map is equal precisely to
the kernel of $\Psi$. This is a $GL(V)$-equivariant interpretation
of the well known fact from from convexity that the closed (in the
weak topology on measures) linear span of all area measures
$dS_{n-1}(K,\bullet)$ on $S^{n-1}$, when $K$ runs through $\ck(V)$,
is equal to all measures $\mu$ on $S^{n-1}$ satisfying
$$\int_{S^{n-1}}\ome\cdot d\mu(\ome)=0.$$
(This fact is an easy consequence of the Minkowski existence
theorem, see e.g. \cite{schneider-book}, Theorem 7.1.2.)

Let $Val^+(V)$ denote the subspace of even valuations, i.e. such
that $\phi(-K)=\phi(K)$ for any $K\in \ck(V)$, and let $Val^-(V)$
denote the subspace of odd valuations, i.e. such that
$\phi(-K)=-\phi(K)$ for any $K\in \ck(V)$. Similarly
$Val^\pm_i(V)$ denote the analogous subspaces in $i$-homogeneous
valuations.

Let $C^+(\pph,\cl)$ denote the subspace of $C(\pph,\cl)$ of even
section, i.e. sections which do not change when one reverses a
coorientation of a hyperplane. Let $C^-(\pph,\cl)$ denote the
subspace of $C(\pph,\cl)$ of odd section, i.e. sections which
change the sign when one reverses a coorientation of a hyperplane.
Let us denote $\Psi^+$ (resp. $\Psi^-$) the restriction of $\Psi$
to $C^+(\pph,\cl)$ (resp. $C^-(\pph,\cl)$).

Let $\PP^\vee(V)$ be the manifold of linear hyperplanes in $V$
(without coorientation). Let $\cl^+\to \PP^\vee(V)$ be the line
bundle whose fiber over $E$ is equal to $Dens(E)$. Let $\cl^-\to
\PP^\vee(V)$ be the line bundle whose fiber over $E$ is equal to
$Dens(E)\otimes or(V/E)$. It is easy to see that one has canonical
isomorphisms
\begin{eqnarray*}
C^+(\pph,\cl)\simeq C(\PP^\vee(V),\cl^+),\\
C^-(\pph,\cl)\simeq C(\PP^\vee(V),\cl^-).
\end{eqnarray*}

Then obviously one has
\begin{eqnarray}\label{Psi-plus}
\Psi^+\colon C(\PP^\vee(V),\cl^+)\to
Val^+_{n-1}(V),\\\label{Psi-minus} \Psi^-\colon
C(\PP^\vee(V),\cl^-)\to Val^-_{n-1}(V).
\end{eqnarray}
The following claim is obvious from the previous discussion.
\begin{claim}\label{Cl:psi-plus-minus}
$\Psi^+$ is an isomorphism of vector spaces. The kernel of
$\Psi^-$ is equal to $\wedge^{n-1}V^*\otimes or(V)$.
\end{claim}

The next theorem is very useful. In the even case it was proved by
Klain \cite{klain}, and in the odd case by Schneider
\cite{schneider-simple}. First recall that a valuation is called
{\itshape simple} if it vanishes on all convex compact sets of
dimension less than $n$.
\begin{theorem}[\cite{klain},\cite{schneider-simple}]\label{T:klain-schneider}
A translation invariant continuous valuation is simple if and only
if it is representable as a sum of a Lebesgue measure and an odd
$(n-1)$-homogeneous translation invariant continuous valuation.
\end{theorem}

The group $GL(V)$ acts continuously and linearly in the space
$Val(V)$ as follows: $(g\phi)(K)=\phi(g^{-1}K)$ for any $g\in
GL(V), \phi\in Val(V),K\in \ck(V)$. Obviously this action
preserves degree of homogeneity and parity of valuations. The next
result was proved by the author \cite{alesker-gafa-01}, it will be
used in the article many times.
\begin{theorem}[Irreducibility Theorem,
\cite{alesker-gafa-01}]\label{T:IrrThm} The natural representation
of $GL(V)$ in $Val^\pm_i(V),\, i=0,1,\dots,n,$ is irreducible, i.e.
there is no $GL(V)$-invariant proper closed subspace.
\end{theorem}

\subsection{Klain-Schneider realizations of
valuations.}\label{Ss:klain-schneider} In this section we describe
$GL(V)$-equivariant realizations of $Val_i^+(V)$ as a subspace of
the space of sections of certain line bundle over the Grassmannian
$Gr_i(V)$, and of $Val_i^-(V)$ as a subspace of a quotient of the
space of sections of certain line bundle over the partial flag
space $\cf_{i,i+1}(V)$. The exposition follows
\cite{alesker-adv-00}, \cite{alesker-gafa-01}. The even case was
also considered in \cite{klain-00} using a slightly different
language. These realizations of even and smooth valuations we call
respectively Klain and Schneider realizations. The reason for such
terminology is that behind these constructions stays a deep
Klain-Schneider theorem \ref{T:klain-schneider}.

Let us start with the even case. Let us denote by $\cm\to Gr_i(V)$
the complex line bundle whose fiber over $E\in Gr_i(V)$ is equal to
$Dens(E)$. Let $\phi\in Val_i^+(V)$. For any $E\in Gr_i(V)$ let us
consider the restriction $\phi|_E$. Clearly $\phi|_E\in Val_i(E)$.
But by Proposition \ref{P:zero-n}(2) (due to Hadwiger)
$Val_i(E)=Dens(E)$. Then $\phi$ defines a section in
$C(Gr_i(V),\cm)$. We get a continuous $GL(V)$-equivariant map
\begin{eqnarray}\label{D:klain-map}
Val_i^+(V)\to C(Gr_i(V),\cm).
\end{eqnarray}
The key fact is that this map is injective. This can be easily
deduced by induction on the dimension from the even case (due to
Klain) of the Klain-Schneider theorem \ref{T:klain-schneider} (see
\cite{alesker-adv-00}, Proposition 3.1, or \cite{klain-00}). We
call this imbedding the Klain imbedding.

Let us consider the odd case. Let $\cf_{i,i+1}(V)$ denote the
partial flag variety
$$\cf_{i,i+1}(V):=\{(E,F)|\, E\subset F, \dim E=i,\dim F=i+1\}.$$

Let us denote by $\cx\to Gr_{i+1}(V)$ the (infinite dimensional)
vector bundle whose fiber over $F\in Gr_{i+1}(V)$ is equal to
$Val^-_i(F)$. Let $\phi\in Val^-_i(V)$. For any $F\in Gr_{i+1}(V)$
let us consider the restriction $\phi|_F\in Val_i^-(F)$. Thus we get
a $GL(V)$-equivariant continuous map
$$Val_i^-(V)\to C(Gr_{i+1}(V),\cx).$$
The key point is that this map is injective. This easily follows
by the induction on dimension from the odd case (due to Schneider)
of the Klain-Schneider theorem \ref{T:klain-schneider} (see
\cite{alesker-gafa-01}, Proposition 2.6 for the details).

Let $\cn\to \cf_{i,i+1}(V)$ denote the line bundle whose fiber
over $(E,F)\in \cf_{i,i+1}(V)$ is equal to $Dens(E)\otimes
or(F/E)$. Applying the map $\Psi^-$ (see (\ref{Psi-minus})) to
every subspace $F\in Gr_{i+1}(V)$ (instead of $V$) we get a
continuous map
\begin{eqnarray*}
C(\cf_{i,i+1}(V),\cn)\to C(Gr_{i+1}(V),\cx).
\end{eqnarray*}
This map is onto. This follows from the fact that $\Psi^-$ is onto
and has finite dimensional kernel.

Thus $Val_i^-(V)$ is realized as a subspace of a quotient of
$C(\cf_{i,i+1}(V),\cn)$. We call this realization the Schneider
realization.

\subsection{Product and convolution of smooth
valuations.}\label{Ss:product-convolution} Let us denote by
$Val^{sm}(V)$ the space of smooth valuations (in sense of
Definition \ref{part1-rep-1}) under the natural action of the
group $GL(V)$ on the space $Val(V)$.

Let $\ck^{sm}(V)$ denote the space of strictly convex compact
subsets of $V$ with $C^\infty$-smooth boundary. A typical example
of smooth valuation is a functional $K\mapsto vol(K+A)$ where
$A\in \ck^{sm}(V)$ is fixed.

Product on smooth translation invariant valuations was defined by
the author in \cite{alesker-gafa-04}. Let us summarize the main
properties of the product in the following theorem. Let us fix on
$V$ a positive Lebesgue measure $vol_V$. Below we denote by
$vol_{V\times V}$ the product measure on $V\times V$
\begin{theorem}[\cite{alesker-gafa-04}]\label{T:product-val}
There exists a bilinear map
$$Val^{sm}(V)\times Val^{sm}(V)\to Val^{sm}(V)$$
which is uniquely characterized by the following two properties:

1) continuity;

2) if $\phi(\bullet)=vol_V(\bullet+A),\, \psi=vol_V(\bullet+B)$
then $$(\phi,\psi)\mapsto vol_{V\times V}(\Delta(\bullet)+(A\times
B))$$ where $\Delta\colon V\to V\times V$ is the diagonal
imbedding.

This bilinear map defines a product making $Val^{sm}(V)$ a
commutative associative algebra with unit (which is the Euler
characteristic).
\end{theorem}
\begin{example}[\cite{alesker-gafa-04}, Proposition 2.2]\label{E:product-2dim}
Assume that $\dim V=2$. Let $\phi(K)$=V(K,A), $\psi(K)=V(K,B)$.
Then $$(\phi\cdot\psi)(K)=\frac{1}{2} V(A,-B)vol(K).$$
\end{example}

Convolution on $Val^{sm}(V)\otimes Dens(V^*)$ was defined by
Bernig and Fu in \cite{bernig-fu}. Let us summarize their result
in the following theorem.
\begin{theorem}[\cite{bernig-fu}]\label{T2}
There exists a bilinear map
$$Val^{sm}(V)\otimes
Dens(V)^*\times Val^{sm}(V)\otimes Dens(V)^*\to Val^{sm}(V)\otimes
Dens(V)^*$$ which is uniquely characterized by the following two
properties:

1) continuity;

2) if $\phi(\bullet)=vol_V(\bullet+A)\otimes vol_V^{-1},\,
\psi=vol_V(\bullet+B)\otimes vol_V^{-1}$ then $$(\phi,\psi)\mapsto
vol_{V}(\bullet+A+ B)\otimes vol_V^{-1}.$$

This bilinear map defines a product making $Val^{sm}(V)\otimes
Dens(V)^*$ a commutative associative algebra with unit (which is
equal to $vol_V\otimes vol_V^{-1}$).
\end{theorem}

\subsection{A technical lemma.}\label{Ss:aLemma}
For a future reference we state a simple and well known lemma (see
\cite{schneider-book}, p. 294, particularly equality (5.3.23)).
\begin{lemma}\label{L:aLemma}
Let $f\colon V\surj W$ be a linear epimorphism of finite
dimensional real vector spaces. Let $k:=\dim Ker (f)$. Let
$vol_{Ker}, vol_W$ be Lebesgue measures on $Ker(f), W$
respectively. Let $vol_V:=vol_{Ker}\otimes vol_W$ be the
corresponding Lebesgue measure on $V$. Let $A\in \ck(V), B\in
\ck(Ker(f))$. Then
$$\frac{1}{k!}\frac{d^k}{d\eps^k}|\big|_{\eps=0} vol_V(A+\eps
B)=vol_{Ker}(B)\cdot vol_W(f(A)).$$
\end{lemma}
{\bf Proof.} It is an application of the Fubini theorem. \qed

\section{Functorial properties of translation invariant valuations.}\label{S:functorial}
\subsection{Pullback of valuations.}\label{pullback}
\begin{definition}\label{D:pullback}
Let $f\colon V\to W$ be a linear map of vector spaces. Let us
define a map, called pullback,
$$f^*\colon Val(W)\to Val(V)$$
by $(f^*\phi)(K)=\phi(f(K))$ for any $K\in \ck(V)$.
\end{definition}
\begin{proposition}\label{P:pullback1}
(i) $f^*$ is a continuous map of Banach spaces.

(ii) $f^*$ preserves degree of homogeneity and parity of
valuations.

(iii) $(f_1\circ f_2)^*=f_2^*\circ f_1^*$.
\end{proposition}
{\bf Proof} is obvious. \qed

\subsection{Pushforward of valuations twisted by
densities.}\label{pushforward} For a linear map $f\colon V\to W$
we are going to define in this section a canonical map, called
pushforward,
$$f_*\colon Val(V)\otimes Dens(V^*)\to Val(W)\otimes Dens(W^*).$$
The main result of this section is the following proposition which
is also a definition.
\begin{proposition}\label{P:pushforward}
(i) Let $f\colon V\to W$ be a linear map of vector spaces. Then
there exists a continuous linear map, called pushforward,
$$f_*\colon Val(V)\otimes Dens(V^*)\to Val(W)\otimes Dens(W^*)$$
which is uniquely characterized by the following property. Let us
fix Lebesgue measures $vol_V$ on $V$, $vol_{Ker}$ on $Ker f$,
$vol_{CoKer}$ on $CoKer f=W/Im (f)$, and let $vol_{Im}$ be the
induced Lebesgue measure on $Im (f)$ (it is obtained as the image
of the measure $\frac{vol_V}{vol_{Ker}}$ on $V/ Ker(f)$ under the
isomorphism $f\colon V/Ker(f)\tilde\to Im(f)$). Then for any $A\in
\ck(V)$
\begin{eqnarray*}
f_*(vol_V(\bullet +A) \otimes vol_V^{-1})=\\\left(\int_{z\in
CoKer(f)}vol_{Im}((\bullet +f(A))\cap z) dvol_{CoKer}(z)\right)
\otimes (vol_{Im}\otimes vol_{CoKer})^{-1}
\end{eqnarray*}
where $vol_{Im}\otimes vol_{CoKer}$ is considered as a Lebesgue
measure on $W$ under the isomorphism $Dens(W)\simeq
Dens(Im(f))\otimes Dens(CoKer(f))$.

(ii)$f_*$ preserves the parity, and $f_*(Val_\bullet(V))\subset
Val_{\bullet+\dim W-\dim V}(W)$.

(iii) $(f_1\circ f_2)_*=f_{1*}\circ f_{2*}$.
\end{proposition}
\begin{remark}\label{R:pushforward}
Before we prove this proposition let us discuss two special cases
of the pushforward. Let us fix Lebesgue measures $vol_V$ on $V$
and $vol_W$ on $W$.  This choice induces isomorphisms
$Dens(V)\tilde\to \CC$, $Dens(W)\tilde\to \CC$. Under these
identifications, $f_*\colon Val(V)\to Val(W)$.

(1) Let us assume that $V$ is a subspace of $W$, and $f\colon
V\inj W$ is the imbedding map. Consider the Lebesgue measure
$vol_{W/V}:=\frac{vol_{W}}{vol_V}$ on $W/V$. For any $\phi\in
Val(V)$
$$(f_*\phi)(K)=\int_{z\in W/V}\phi(K\cap z)d vol_{W/V}(z).$$

(2) Let us assume that $W$ is a quotient space of $V$, and
$f\colon V \twoheadrightarrow W$ is the quotient map. Let $\phi\in
Val(V)$ has the form $\phi(K)=vol_V(K+A)$ where $A\in \ck(V)$ is
fixed. Then for any $K\in \ck(W)$
$$(f_*\phi)(K)=vol_W(K+f(A)).$$
\end{remark}

{\bf Proof} of Proposition \ref{P:pushforward}. The uniqueness
follows immediately from the McMullen's conjecture. Let us prove
the existence. Let us decompose $f$ as a composition of a linear
surjection $p\colon V\surj X$ followed by a linear injection
$j\colon X\inj W$; thus $f=j\circ p$. Such a decomposition is
unique up to an isomorphism (in the obvious sense).

Let us define $j_*\colon Val(X)\otimes Dens(X^*)\to Val(W)\otimes
Dens(W^*)$ as in Remark \ref{R:pushforward}(1). Let us define now
$p_*\colon Val(V)\otimes Dens(V^*)\to Val(X)\otimes Dens(X^*)$.
Let us fix Lebesgue measures $vol_{Ker}$ on $Ker(p)$ and $vol_X$
on $X$. This induces a Lebesgue measure $vol_V:=vol_{Ker}\otimes
vol_X$ on $V$. Let $\phi\in Val(V)$. Fix a set $L\in \ck(V)$.
Consider a valuation $\tau$ on $Ker(p)$ defined by
$$\tau(S)=\frac{1}{k!}\frac{d^k}{d\eps^k}\big|_{\eps=0}\phi(L+\eps
S)$$ where $S\in \ck(Ker(p))$, $k=\dim (Ker(p))$. Recall that by a
result of McMullen \cite{mcmullen-euler}, $\phi(L+\eps S)$ is a
polynomial in $\eps\geq 0$ of degree at most $k$. It is easy to
see that $\tau$ is a $k$-homogeneous translation invariant
continuous valuation on $Ker(p)$. By a result of Hadwiger
\cite{hadwiger-book} $\tau$ must be proportional to $vol_{Ker}$
with a constant depending on $L$ and $\phi$:
$$\tau=C(L,\phi)vol_{Ker}.$$
It is easy to see that $C(L,\phi)$ depends continuously on $L\in
\ck(V)$ and $\phi\in Val(V)$, and linearly on $\phi$. Let $K\in
\ck(X)$. Let $\tilde K\in \ck(V)$ be an arbitrary convex compact
set such that $p(\tilde K)=K$. Define
$$\tilde \phi(K):=C(\tilde K,\phi).$$
\begin{claim}\label{claim}
$\tilde\phi (K)$ does not depend on a choice of $\tilde K$ such
that $p(\tilde K)=K$.
\end{claim}
{\bf Proof.} Since $C(\tilde K,\phi)$ is continuous in $\phi\in
Val(V)$ when $\tilde K$ in fixed, by the McMullen's conjecture, it
is enough to prove the claim for $\phi(\bullet)=vol_V(\bullet +A)$
for $A\in \ck(V)$. Let us fix $S\in \ck(Ker(p))$ with
$vol_{Ker}(S)=1$. Then $C(\tilde
K,\phi)=\frac{1}{k!}\frac{d^k}{d\eps^k}\big|_{\eps=0}vol_V(\tilde
K+A+\eps S)$. But the last expression is equal to $vol_X(p(\tilde
K+A))=vol_X(K+p(A))$ by Lemma \ref{L:aLemma}. \qed

Next $\tilde\phi$ is a continuous translation invariant valuation
on $X$. Translation invariance is obvious. In order to prove
continuity and valuation property let us fix a linear right
inverse of $p$, $s\colon X\to V$. For any $K\in\ck(X)$ let us
choose $\tilde K:=s(K)$. Thus $\tilde \phi(K)=C(s(K),\phi)$.
Clearly the last expression is a continuous valuation in $K\in
\ck(X)$. Let us define
\begin{eqnarray}\label{p-star}
p_*(\phi\otimes vol_V^{-1}):=\tilde\phi\otimes vol_X^{-1}.
\end{eqnarray}
It is easy to see that the definition of $p_*$ does not depend on
a choice of Lebesgue measures $vol_V, vol_X$. Finally define
$$f_*:=j_*\circ p_*.$$
It follows from the construction that $f_*$ satisfies the
assumptions of the proposition.

(ii) is obvious.

(iii) Let $U\overset{f_2}{\to}V\overset{f_1}{\to} W$. We have to
show
\begin{eqnarray}\label{func}
(f_1\circ f_2)_*=f_{1*}\circ f_{2*}.\end{eqnarray} First we will
show this in a number of special cases. If $f_1$ is an injection
and $f_2$ is a surjection the equality (\ref{func}) is clear by
the construction of pushforward from the proof of part (i).

\begin{lemma}\label{inj_inj}
The equality (\ref{func}) holds if $f_1,f_2$ are injections.
\end{lemma}
{\bf Proof.} We may and will assume for simplicity that $f_1,f_2$
are imbeddings of linear subspaces, thus $U\subset V\subset W$.
Let us fix Lebesgue measures $vol_U,vol_{V/U},vol_{W/V}$ on
$U,V/U,W/V$ respectively. Let $vol_V:=vol_U\otimes vol_{V/U},\,
vol_W:=vol_V\otimes vol_{W/V}$ be Lebesgue measures on $V$, $W$
respectively. Let us fix $\phi\in Val(U)$, $K\in \ck(W)$. We have
\begin{eqnarray*}
(f_{2*}(\phi\otimes vol_U^{-1}))(K)=\left(\int_{z\in
V/U}\phi(K\cap z)d vol_{V/U}(z)\right)\otimes
vol_V^{-1},\end{eqnarray*}
\begin{eqnarray*}
(f_{1*}(f_{2*}(\phi\otimes vol_U^{-1})))(K)=\left(\int_{w\in
W/V}\left(f_{2*}(\phi\otimes vol_U^{-1})(K\cap w)\otimes
vol_V\right)\cdot
d vol_{W/V}(w)\right)\otimes vol_W^{-1}=\\
\left(\int_{x\in W/U}\phi(K\cap x)d\left(
\frac{vol_W}{vol_U}\right)(x)\right)\otimes vol_W^{-1}=(f_1\circ
f_2)_*(\phi\otimes vol_U^{-1}).
\end{eqnarray*}
Lemma is proved. \qed
\begin{lemma}\label{surj_surj}
The equality (\ref{func}) holds if $f_1,f_2$ are surjections.
\end{lemma}
{\bf Proof.} Let us fix Lebesgue measures $vol_1,vol_2,vol_W$ on
$Ker(f_1),Ker(f_2),W$ respectively. Set $vol_U:=vol_1\otimes
vol_2\otimes vol_W\in Dens(U)$, $vol_V=vol_1\otimes vol_W\in
Dens(V)$. Assume that $\phi\in Val(U)$ has the form
$$\phi(\bullet)=vol_U(\bullet +A)$$ where $A\in \ck(U)$ is fixed.
Then
\begin{eqnarray*}
(f_{2*}(\phi\otimes vol_U^{-1}))(K)=vol_V(K+f_2(A))\otimes
vol_V^{-1},\end{eqnarray*}
\begin{eqnarray}\label{s1}
\left(f_{1*}(f_{2*}(\phi\otimes
vol_U^{-1}))\right)(K)=vol_W(K+f_1(f_2(A)))\otimes vol_W^{-1}.
\end{eqnarray}
On the other hand
\begin{eqnarray}\label{s2}
(f_1\circ f_2)_*(\phi\otimes vol_U^{-1})(K)=vol_W(K+(f_1\circ
f_2)(A))\otimes vol_W^{-1}
\end{eqnarray}
Comparing (\ref{s1}) and (\ref{s2}) one concludes the lemma using
the McMullen's conjecture. \qed

\begin{lemma}\label{surj_inj}
The equality (\ref{func}) holds if $f_1$ is a surjection, $f_2$ is
an injection.
\end{lemma}
{\bf Proof.} We will assume for simplicity and without loss of
generality that $U\subset V$ and $f_2\colon U\inj V$ is the
identity imbedding. Also we may assume that $W$ is a quotient
space of $V$, and $f_1\colon V\surj W$ is the canonical quotient
map.

\underline{Step 1.} We consider the case
$$Ker(f_1)\subset U.$$ Denote $k:=\dim Ker (f_1)$. Let us fix
Lebesgue measures $vol_{Ker}$ on $Ker(f_1)$, $vol_U$ on $U$,
$vol_{V/U}$ on $V/U$. Let $vol_V=vol_U\otimes vol_{V/U},$
$vol_W:=\frac{vol_V}{vol_{Ker}}$ be the Lebesgue measure on $W$.
Let $\phi\in Val(U)$ has the form
$$\phi(\bullet)=vol_U(\bullet +A)$$ for some $A\in \ck(U)$. Then
for any $L\in \ck(V)$ we have
\begin{eqnarray*}
\left(f_{2*}(\phi\otimes vol_U^{-1})\right)(L)=\left(\int_{z\in
V/U}vol_U((L\cap z)+A) d vol_{V/U}(z)\right) \otimes vol_V^{-1}.
\end{eqnarray*}
Now let us fix a linear right inverse of $f_1$
$$s\colon  W\to V.$$
Let $M\in \ck(W)$. Let $\tilde M:= s(M)$. Let us fix a set $S\in
\ck(Ker(f_1))$ with $vol_{Ker}(S)=1$. By the construction in the
proof of part (i) we get
\begin{eqnarray}\label{s3.5}
\left(f_{1*}(f_{2*}(\phi\otimes vol_U^{-1}))\right)
(M)=\\\label{s4}\frac{1}{k!}\frac{d^k}{d\eps^k}\big|_{\eps=0}\left(\int_{z\in
V/U} vol_U\left(((\tilde M+\eps S)\cap z)+A\right)d
vol_{V/U}(z)\right)\otimes vol_W^{-1}=\\\label{s4.5}
\frac{1}{k!}\frac{d^k}{d\eps^k}\big|_{\eps=0}\left(\int_{z\in V/U}
vol_U\left((\tilde M\cap z)+\eps S+A\right)d
vol_{V/U}(z)\right)\otimes vol_W^{-1}=\\\label{s4.6}
\left(\int_{z\in (V/Ker f_1)/(U/Ker f_1)}vol_W\left((M\cap
z)+f_1(A)\right) d vol_{V/U}(z)\right)\otimes vol_W^{-1}.
\end{eqnarray}
where the equality (\ref{s4}) follows from the fact that $(\tilde
M+\eps S)\cap z=(\tilde M\cap z)+\eps S$ since $S\subset
Ker(f_1)\subset U$, and in (\ref{s4.5}) we used the identification
$V/U=(V/Ker(f_1))/(U/Ker(f_1))$.

Let us compute now $(f_1\circ f_2)_*(\phi\otimes vol_U^{-1})$. The
map $f_1\circ f_2$ factorizes as
$$U\overset{p}{\surj} U/Ker(f_1) \overset{j}{\inj} V/Ker(f_1)=W.$$
By the construction from part (i)
$$(f_1\circ f_2)_*=j_*\circ p_*.$$ We have for any $L\in \ck(V)$
$$\left(p_*(\phi\otimes
vol_U^{-1})\right)(L)=vol_{U/Ker}(L+p(A))\otimes
vol_{U/Ker}^{-1}.$$ Next for any $M\in \ck(W)$
\begin{eqnarray}\label{s5}
\left(j_*(p_*(\phi\otimes vol_U^{-1}))\right) (M)=\\\label{s6}
\left( \int_{z\in V/U}vol_{W} \left((M\cap z)+p(A)\right) d
vol_{V/U}(z)\right)\otimes vol_W^{-1}.
\end{eqnarray}
Comparing (\ref{s4.5}) and (\ref{s6}) and observing that
$p(A)=f_1(A)$, we conclude Step 1.

\underline{Step 2.} At this step we will assume that $f_1\circ
f_2$ is injective.

Let us choose $s\colon W\to V$ a right inverse of $f_1$ so that
$s(W)\supset f_2(U)$. Let $H:=Ker (f_1)$. To simplify the notation
and without loss of generality we will identify $U$ with $f_2(U)$,
and $W$ with $s(W)$. Let us fix Lebesgue measures $vol_U$ on $U$,
$vol_{W/U}$ on $W/U$, $vol_H$ on $H$. Let $vol_W:=vol_U\otimes
vol_{W/U}$, $vol_V:=vol_W\otimes vol_H$,
$vol_{V/U}:=vol_{W/U}\otimes vol_H$ be the corresponding Lebesgue
measures on $W,V,V/U$ respectively.

Let $\phi\in Val(U)$. Let us fix $K\in \ck(W)$. Then we have
\begin{eqnarray}\label{s7}
\left((f_1\circ f_2)_*(\phi\otimes
vol_U^{-1})\right)(K)=\int_{z\in W/U}\phi(K\cap z)d
vol_{W/U}(z)\otimes vol_W^{-1}.
\end{eqnarray}
Let us fix a subset $S\subset H$ with $vol_H(S)=1$. Denote
$k:=\dim H$. Then for any $L\in \ck(V)$ we have
\begin{eqnarray*}
\left(f_{2*}(\phi\otimes vol_U^{-1})\right)(L)=\left(\int_{x\in
V/U}\phi(L\cap x)d vol_{V/U}(x)\right)\otimes vol_V^{-1}.
\end{eqnarray*}
Next for any $K\in \ck(W)$
\begin{eqnarray}
\left(f_{1*}(f_{2*}(\phi\otimes
vol_U^{-1}))\right)(K)=\\\label{s8}
\frac{1}{k!}\frac{d^k}{d\eps^k}\big|_{\eps=0} \left(\int_{x\in
V/U}\phi\left((K+\eps S)\cap x\right) d vol_{V/U}(x)\right)
\otimes vol_W^{-1}.
\end{eqnarray}
But since $K\subset W$ and $S\subset H$, $x\in V/U$,
$$(K+\eps S)\cap x=(K\cap x)+(\eps S\cap x)$$
and $\eps S\cap x$ is either a point or the empty set. Hence
(\ref{s8}) can be continued as
\begin{eqnarray*}
\frac{1}{k!}\frac{d^k}{d\eps^k}\big|_{\eps=0}
\left(\eps^k\int_{x\in W/U}\phi(K\cap x)d vol_{W/U}
(x)\right)\otimes vol_W^{-1}=\\
\left(\int_{x\in W/U}\phi(K\cap x)d vol_{W/U} (x)\right)\otimes
vol_W^{-1}\overset{(\ref{s7})}{=}\left((f_1\circ
f_2)_*(\phi\otimes vol_U^{-1})\right)(K).
\end{eqnarray*}
This completes Step 2.

\underline{Step 3.} Let us consider finally the case of general
surjection $f_1$ and general injection $f_2$.

We will assume again that $U$ is a subspace of $V$, and $W$ is a
quotient space of $V$. Set $A:=U\cap Ker(f_1)$, $X:=V/A$. We can
decompose uniquely $f_2\colon V\surj W$ as a composition of two
surjections
$$V\overset{q}{\surj} X\overset{p}{\surj} W$$ where $p\colon
V\surj V/A=X$ is the canonical surjection. By Lemma
\ref{surj_surj}
\begin{eqnarray}\label{s8.01}
f_{1*}=p_*\circ q_*.
\end{eqnarray}
Let us denote $Y:=U/A$. Let $t\colon
U\surj Y$ be the canonical surjection, and $j\colon Y\inj X$ be
the natural imbedding. Note that $Ker(q)=A\subset U$. Then by Step
1 and the construction of the pushforward from part (i) we get
\begin{eqnarray}\label{s9}
q_*\circ f_{2*}=(q\circ f_2)_*=(j\circ t)_*=j_*\circ t_*.
\end{eqnarray}
Clearly $p\circ j\colon Y\to W$ is injective. Hence by Step 2
\begin{eqnarray}\label{s10}
p_*\circ j_*=(p\circ j)_*.
\end{eqnarray}
Using (\ref{s8.01}),(\ref{s9}), and (\ref{s10}) we obtain
\begin{eqnarray*}
f_{1*}\circ f_{2*}=p_*\circ q_*\circ f_{2*}=p_*\circ j_*\circ
t_*=(p\circ j)_*\circ t_*=(f_1\circ f_2)_*
\end{eqnarray*}
where the last equality follows from the facts that $f_1\circ
f_2=(p\circ j)\circ t$, $t$ is surjective, $p\circ j$ is
injective, and from the construction of the pushforward from the
proof of part (i). Lemma \ref{surj_inj} is proved. \qed


Now let us finish the proof of Proposition \ref{P:pushforward}(ii)
for general $f_1,f_2$. Let us decompose
$$f_1=j_1\circ p_1,\, f_2=j_2\circ p_2$$
where $j_1,j_2$ are injections, $p_1,p_2$ are surjections. By
definition
\begin{eqnarray}\label{s3}
f_{1*}\circ f_{2*}=j_{1*}\circ p_{1*}\circ j_{2*}\circ p_{2*}.
\end{eqnarray}
Let us decompose $p_1\circ j_2=j_3\circ p_3$ where $j_3$ is an
injection, $p_3$ is a surjection. By Lemma \ref{surj_inj},
$p_{1*}\circ j_{2*}=(p_1\circ j_2)_*=j_{3*}\circ p_{3*}$. Hence
using this, (\ref{s3}), and Lemmas \ref{inj_inj}, \ref{surj_surj}
we get
\begin{eqnarray*}
f_{1*}\circ f_{2*}=j_{1*}\circ j_{3*}\circ p_{3*} \circ p_{2*}=
(j_1\circ j_3)_*\circ (p_3\circ p_2)_*=(f_1\circ f_2)_*
\end{eqnarray*}
where the last equality follows from the construction of the
pushforward given in the proof of part (i). Proposition
\ref{P:pushforward} is proved. \qed

\subsection{Relations to product and convolution.}\label{relation} We will explain the relation of
pullback and pushforward to product and convolution. But first we
will have to remind the notion of exterior product of smooth
translation invariant valuations. Let $V,W$ be finite dimensional
real vector spaces. In \cite{alesker-gafa-04} the author has
defined a continuous linear map
\begin{eqnarray}\label{r1}
Val^{sm}(V)\times Val^{sm}(W)\to Val(V\times W)
\end{eqnarray}
called the {\itshape exterior product}. For $\phi\in Val^{sm}(V),
\psi\in Val^{sm}(W)$ their exterior product is denoted by
$\phi\boxtimes \psi$. The map (\ref{r1}) is uniquely characterized
by the following property. Let $\phi(\bullet)=vol_V(\bullet +A),\,
\psi(\bullet)=vol_W(\bullet +B)$ where $vol_V,\, vol_W$ are
Lebesgue measures on $V,W$ respectively, $A\in \ck^{sm}(V),\, B\in
\ck^{sm}(W)$. Then
$$(\phi\boxtimes \psi)(K)=(vol_V\boxtimes vol_W)(K+(A\times B))$$
for any $K\in \ck(V\times W)$, and where $vol_V\boxtimes vol_W$
denotes the usual product measure. Note that the exterior product
of two smooth valuations may not be smooth.

In Appendix of this article we show a slightly more precise
statement which will be useful later for some technical reasons.
Namely it is shown that the exterior product extends (uniquely) to
a continuous bilinear map
$$Val(V)\times Val^{sm}(W)\to Val(V\times W),$$
i.e. the first variable may be replaced by continuous valuations
instead of smooth.

The following proposition is essentially the definition of the
product of smooth valuations from \cite{alesker-gafa-04}.
\begin{proposition}\label{P:r1}
For any $\phi,\psi\in Val^{sm}(V)$
$$\phi\cdot \psi=\Delta^*(\phi\boxtimes \psi)$$
where $\Delta\colon V\inj V\times V$ denotes the diagonal
imbedding.
\end{proposition}

Let us explain now the relation between the convolution and the
pushforward. Clearly $Dens(V\times W)^*=Dens(V^*)\otimes
Dens(W^*)$. Hence the exterior product (\ref{r1}) tensored with
$Id_{Dens(V\times W)^*}$ gives a continuous bilinear map
\begin{eqnarray*}\label{r2}
\left(Val^{sm}(V)\otimes
Dens(V^*)\right)\times\left(Val^{sm}(W)\otimes
Dens(W^*)\right)\to\\ Val(V\times W)\otimes Dens(V\times W)^*
\end{eqnarray*}
which will also be called exterior product and denoted by
$\boxtimes$. Let us denote by $a\colon V\times V\to V$ the
addition map, i.e. $a(x,y)=x+y$.
\begin{proposition}\label{P:r2}
For any $\phi,\psi\in Val^{sm}(V)\otimes Dens(V^*)$ one has
$$\phi\ast\psi=a_*(\phi\boxtimes \psi).$$
\end{proposition}
{\bf Proof.} Let us fix a Lebesgue measure $vol_V$ on $V$. By
continuity and the McMullen's conjecture it is enough to prove the
proposition for $\phi(\bullet)=vol_V(\bullet+A),\,
\psi(\bullet)=vol_V(\bullet+B)$ with $A,B\in \ck^{sm}(V)$. Then
\begin{eqnarray*}
(\phi\boxtimes \psi)(\bullet)=(vol_V\boxtimes
vol_V)(\bullet+(A\times B))\otimes (vol_V\otimes vol_V)^{-1}.
\end{eqnarray*}
Next
\begin{eqnarray*}
(a_*(\phi\boxtimes \psi))(\bullet)=vol_V(\bullet+a(A\times
B))\otimes vol_V^{-1}=vol_V(\bullet +A+B)\otimes
vol_V^{-1}=(\phi\ast\psi)(\bullet).
\end{eqnarray*}
\qed
\subsection{Homomorphism property of
pushforward.}\label{Ss:homomor-of-push} The main result of this
section is the following proposition.
\begin{proposition}\label{P:homom-of-push}
Let $p\colon X\surj Y$ be a linear epimorphism of vector spaces.
Then for any $\phi\in Val^{sm}(X)\otimes Dens(X^*)$ the
pushforward $p_*\phi$ is smooth, i.e. $p_*\phi\in
Val^{sm}(Y)\otimes Dens(Y^*)$, and
$$p_*\colon Val^{sm}(X)\otimes Dens(X^*)\to Val^{sm}(Y)\otimes
Dens(Y^*)$$ is a homomorphism of algebras (when both spaces are
equipped with convolution).
\end{proposition}

It is easy to see that $p_*\phi$ is smooth if $p$ is surjective
and $\phi$ is smooth. In order to prove the second statement of
the proposition we will need another proposition.
\begin{proposition}\label{L:ext-push}
Let
\begin{eqnarray*}
f_1\colon V_1\to W_1,\\
f_2\colon V_2\to W_2
\end{eqnarray*}
be linear maps. Let $\phi_i\in Val(W_i)\otimes Dens(W_i^*)$,
$i=1,2$. Assume that $f_1$ is surjective and $\phi_1$ is smooth.
Then
$$(f_1\boxtimes f_2)_*(\phi_1\boxtimes \phi_2)=f_{1*}\phi_1\boxtimes
f_{2*}\phi_2.$$
\end{proposition}
{\bf Proof.} Let us fix Lebesgue measures $vol_{Ker}$ on
$Ker(f_1)$, $vol_{W_1}$ on $W_1$, $vol_{W_2}$ on $W_2$,
$vol_{V_2}$ on $V_2$. Let $vol_{V_1}=vol_{Ker}\otimes vol_{W_1}$
be the induced Lebesgue measure on $V_1$. Observe that by the
Appendix to this article, both sides of the last equality are
continuous with respect to $\phi_1\in Val^{sm}(V_1)\otimes
Dens(V_1^*),\phi_2\in Val(V_2)\otimes Dens(V_2^*)$.

 Hence, by the McMullen's conjecture, we may assume that
$$\phi_i(\bullet)=vol_{V_i}(\bullet +A_i)\otimes vol_{V_i}^{-1},\,
i=1,2.$$ Then
\begin{eqnarray}\label{hop1}
(\phi_1\boxtimes \phi_2)(\bullet)=(vol_{V_1}\boxtimes
vol_{V_2})(\bullet +(A_1\times A_2))\otimes (vol_{V_1}\boxtimes
vol_{V_2})^{-1}.
\end{eqnarray}
Then
$$
(f_{i*}\phi_i)(\bullet)=vol_{W_i}(\bullet +f_i(A_i))\otimes
vol_{W_i}^{-1},i=1,2,$$
\begin{eqnarray*}
(f_1\boxtimes f_2)_*(\phi_1\boxtimes\phi_2)(\bullet)=
(vol_{W_1}\boxtimes vol_{W_2})(\bullet + (f_1(A_1)\times
f_2(A_2)))\otimes (vol_{W_1}\boxtimes
vol_{W_2})^{-1}=\\(f_{1*}\phi_1\boxtimes f_{2*}\phi_2)(\bullet).
\end{eqnarray*}
Proposition \ref{L:ext-push} is proved. \qed

{\bf Proof} of Proposition \ref{P:homom-of-push}. Let
\begin{eqnarray*}
a_X\colon X\times X\to X,\\
a_Y\colon Y\times Y\to Y
\end{eqnarray*}
be the addition maps. Then
\begin{eqnarray*}
p_*(\phi\ast\psi)=p_*(a_{X*}(\phi\boxtimes \psi))=(p\circ
a_X)_*(\phi\boxtimes \psi)=\\
(a_Y\circ (p\times p))_*(\phi\boxtimes \psi)=a_{Y*}((p\times
p)_*(\phi\boxtimes \psi)\overset{\mbox{Lemma }
\ref{L:ext-push}}{=}\\ a_{Y*}(p_*\phi\boxtimes p_*\psi)=
p_*\phi\ast p_*\psi.
\end{eqnarray*}
\qed

\subsection{Base change theorem.}\label{Ss:base-change} Recall
that a commutative diagram of linear maps of vector spaces
$$\square[A`B`C`D; f` g`u`v]$$
is called a {\itshape Cartesian square} if it is isomorphic to a
diagram
$$\square[Y\times_Z X`X`Y`Z; pr_X` pr_Y`u`v]$$
where $Y\times_Z X:=\{(y,x)\in Y\times X|\, v(y)=u(x)\}$ and
$pr_X\colon Y\times_Z X\to X$ and $pr_Y\colon Y\times_Z X\to Y$
are the natural maps.
\begin{lemma}\label{L:bch}
(i) Let
$$\square[\tilde X`\tilde Y`X`Y;\tilde f`\tilde g`g`f]$$
be a Cartesian square of vector spaces such that $f\oplus g\colon
X\oplus \tilde Y\to Y$ is onto. Then there exists a canonical
isomorphism
$$\frac{Dens(\tilde
X^*)}{Dens(X^*)}\tilde \to \frac{Dens(\tilde Y^*)}{Dens(Y^*)}.$$

(ii) The following transitivity property of the isomorphism from
part (i) holds. Assume that we have the following commutative
diagram:
$$\bfig \putsquare<1`1`1`1;500`500>(0,500)[\tilde X`\tilde Y`X`Y;\tilde f_1`\tilde
g`g`f_1]\putsquare<1`0`1`1;500`500>(500,500)[\phantom{\tilde
Y}`\tilde Z`\phantom{Y}`Z;\tilde f_2``\hat g`f_2]\efig .$$ If the
two small squares are Cartesian then the exterior contour is
Cartesian. If $f_1\oplus g$ and $f_2\oplus \hat g$ are onto, then
$(f_2\circ f_1)\oplus \hat g$ is onto.

Moreover in the last case the isomorphism $\frac{Dens(\tilde
X)}{Dens(X)}\tilde\to \frac{Dens(\tilde Z)}{Dens(Z)}$
corresponding to the exterior contour of the diagram by the part
(i) of the lemma, is equal to the composition of the isomorphisms
$$\frac{Dens(\tilde
X)}{Dens(X)}\tilde\to \frac{Dens(\tilde Y)}{Dens(Y)}\tilde \to
\frac{Dens(\tilde Z)}{Dens(Z)}$$ corresponding to the small
Cartesian squares.
\end{lemma}
{\bf Proof.} (i) We have the short exact sequence of vector spaces
$$0\to \tilde X\overset{\tilde f\times -\tilde g}{\to}X\oplus
\tilde Y\overset{f\oplus g}{\to}Y\to 0.$$ Hence
$$Dens(X\oplus \tilde Y)\simeq Dens(\tilde X)\otimes Dens(Y).$$
But on the other hand $Dens(X\oplus \tilde Y)=Dens(X)\otimes
Dens(\tilde Y)$. Hence
$$Dens(\tilde X)\otimes Dens(Y)\simeq Dens(X)\otimes
Dens(\tilde Y).$$ Dualization of this isomorphism implies part
(i).

The proof of part (ii) we leave to the reader. \qed

We have the following result we call the base change property.
Roughly put, it says that for a Cartesian square as above one has
$$g^*\circ f_*=\tilde f_*\circ \tilde g^*.$$
\begin{theorem}[Base change theorem]\label{T:base-change}
Let
\begin{eqnarray}\label{ob0}
\square[\tilde X`\tilde Y`X`Y;\tilde f`\tilde g`g`f]
\end{eqnarray} be a Cartesian square of vector spaces such that
$f\oplus g\colon X\oplus \tilde Y\to Y$ is onto. Consider the
following two maps
$$Val(X)\otimes Dens(\tilde X^*)\to Val(\tilde Y)\otimes
Dens(\tilde Y^*),$$ the first map given is by the composition
\begin{eqnarray*}
Val(X)\otimes Dens(\tilde X^*)=\left(Val(X)\otimes
Dens(X^*)\right)\otimes \frac{Dens(\tilde
X^*)}{Dens(X^*)}\overset{f_*\otimes Id}{\to}\\\left(Val(Y)\otimes
Dens(Y^*)\right)\otimes\frac{Dens(\tilde
X^*)}{Dens(X^*)}=Val(Y)\otimes Dens(\tilde Y^*)\overset{g^*\otimes
Id_{Dens(\tilde Y^*)}}{\to}\\Val(\tilde Y)\otimes Dens(\tilde Y^*)
\end{eqnarray*} where we have used the identification
$Dens(Y^*)\otimes \frac{Dens(\tilde X^*)}{Dens(X^*)}\simeq
Dens(\tilde Y^*)$ from Lemma \ref{L:bch}(i); and the second map is
given by the composition
\begin{eqnarray*}
Val(X)\otimes Dens(\tilde X^*)\overset{\tilde g^*\otimes
Id_{Dens(\tilde X^*)}}{\to}Val(\tilde X)\otimes Dens(\tilde
X^*)\overset{\tilde f_*}{\to}Val(\tilde Y)\otimes Dens(\tilde
Y^*).
\end{eqnarray*}
Then these two maps coincide.
\end{theorem}
{\bf Proof.} We will prove this result in several steps. In the
first two steps we will show that it is enough to prove the
theorem under the assumption that each of $f$ and $g$ is either
injection or surjection. Then we will prove the theorem in each of
these cases.

\underline{Step 1.} Transitivity with respect to $f$.

Assume that we have a commutative diagram
\begin{eqnarray*}
\bfig \putsquare<1`1`1`1;500`500>(0,500)[\tilde X`\tilde
Y`X`Y;\tilde f_1`\tilde
g`g`f_1]\putsquare<1`0`1`1;500`500>(500,500)[\phantom{\tilde
Y}`\tilde Z`\phantom{Y}`Z;\tilde f_2``\hat g`f_2]\efig
\end{eqnarray*}
such that $f_1\oplus g$ and $f_2\oplus \hat g$ are onto and two
small squares satisfy conclusions of the theorem. Then the diagram
of the exterior contour also satisfies these conclusions by
Propositions \ref{P:pullback1}(iii), \ref{P:pushforward}(iii), and
Lemma \ref{L:bch}(ii).

\underline{Step 2.} Transitivity with respect to $g$.

Assume that we have a commutative diagram
\begin{eqnarray*}
\bfig \putsquare<1`1`1`1;600`400>(0,750)[\tilde X` \tilde Y `X`
Y;\tilde f`\tilde g_1`g_1`f]
\putsquare<0`1`1`1;600`400>(0,350)[\phantom{X}`\phantom{Y}`\tilde
Z `Z;`\tilde g_2`g_2`\hat f] \efig
\end{eqnarray*}
such that $f\oplus g_1$ and $\hat f\oplus g_2$ are onto, and the
small squares satisfy the conclusions of the theorem. Then, as in
Step 1, the diagram of the exterior contour also satisfies these
conclusions by Propositions \ref{P:pullback1}(iii),
\ref{P:pushforward}(iii), and Lemma \ref{L:bch}(ii).

\underline{Step 3.} Let us assume that $g$ is surjective.

Then we may and will assume that $$\tilde Y=Y\oplus L$$ and $g$ is
the projection $pr_Y\colon Y\oplus L\to Y$. Hence $\tilde
X=X\oplus L$, $\tilde g$ is the projection $pr_X\colon X\oplus
L\to X$, and $\tilde f=f\boxtimes Id_L\colon X\oplus L\to Y\oplus
L$. Thus the diagram (\ref{ob0}) becomes equal to the diagram
\begin{eqnarray}\label{ob1}\square<1`1`1`1;800`400>[ X\oplus L`Y\oplus
L`X`Y;f\boxtimes Id_L`pr_X`pr_Y`f].\end{eqnarray}

Observe that we have canonical isomorphisms
$$\frac{D(\tilde X^*)}{D(X^*)}\simeq \frac{D(\tilde
Y^*)}{D(Y^*)}\simeq D(L^*).$$

Let us denote $l:=\dim L$. Let us fix Lebesgue measures $vol_L$ on
$L$ and $vol_X$ on $X$. Let us fix also $S\in \ck(L)$ such that
$vol_L(S)=1$. Abusing the notation we will denote the first map in
the statement of the theorem by $"pr_Y^*\circ f_*"$, and the
second map by $"(f\times Id_L)_*\circ pr_X^*"$. We have to show
that they coincide. We may decompose $f$ into a composition of
injection and surjection and, using Step 1, prove the result
separately in each case.

\underline{Case a.} Assume that $f$ is surjective.

Then $X$ may and will be assumed to be equal to $Y\oplus M$, and
$f\colon Y\oplus M\to Y$ is the natural projection. Thus $\tilde
X=Y\oplus M\oplus L$. Then the diagram (\ref{ob1}) becomes equal
to
\begin{eqnarray}\label{ob2}\square<1`1`1`1;800`400>[ Y\oplus M\oplus L`Y\oplus
L`X=Y\oplus M`Y;f\boxtimes Id_L`pr_X`pr_Y`f].\end{eqnarray} Let us
fix a Lebesgue measure $vol_M$ on $M$. Let
$vol_Y=\frac{vol_X}{vol_M}\in Dens(Y)$. Let us denote $m:=\dim M$.
Let us fix $T\in \ck(M)$ such that $vol_M(T)=1$. Let us fix
$\phi\in Val(Y\oplus M)\otimes Dens((Y\oplus M\oplus L)^*)$.
Finally let us fix an arbitrary subset $K\in \ck(Y\oplus L)$. Then
we have
\begin{eqnarray}
"(pr_Y^*\circ f_*)"(\phi)(K)=("f_*"\phi)(pr_Y(K))=\\
\frac{1}{m!}\frac{d^m}{d\eps^m}\big|_0\phi(pr_Y(K)\times \eps
T)\otimes vol_M
=\\\label{ob3}\frac{1}{m!}\frac{d^m}{d\eps^m}\big|_0\phi\left(pr_X(K\times
\eps T)\right)\otimes vol_M.
\end{eqnarray}
On the other hand
\begin{eqnarray}
"(f\boxtimes Id_L)_*\circ
pr_X^*"(\phi)(K)=\\\frac{1}{m!}\frac{d^m}{d\eps^m}\big|_0"pr_X^*"(\phi)(K\times
\eps T)\otimes vol_M=\\\label{ob4}
\frac{1}{m!}\frac{d^m}{d\eps^m}\big|_0\phi(pr_X(K\times \eps
T))\otimes vol_M.
\end{eqnarray}
Since (\ref{ob3})=(\ref{ob4}), Case a is proved.

\underline{Case b.} Assume that $f\colon X\to Y$ is injective.

It suffices to consider the case when $X$ is a subspace of $Y$,
and $f$ is the identity imbedding. Let us fix Lebesgue measures
$vol_X$ on $X$, $vol_Y$ on $Y$. Set
$vol_{Y/X}:=\frac{vol_Y}{vol_X}\in Dens(Y/X)$. Then the diagram
(\ref{ob1}) becomes equal to
\begin{eqnarray}\label{ob5}\square<1`1`1`1;800`400>[ X\oplus L`Y\oplus
L`X`Y;f\boxtimes Id_L`pr_X`pr_Y`f].\end{eqnarray} We have
\begin{eqnarray}
"pr^*_Y\circ f_*"(\phi)(K)=("f_*"\phi)(pr_Y(K))=\\
\int_{z\in Y/X}\phi(pr_Y(K)\cap z)dvol_{Y/X}(z)\otimes
vol_{Y/X}^{-1}=\\\label{ob6} \int_{z\in Y/X}\phi(pr_Y(K\cap
(z\times L)))dvol_{Y/X}(z)\otimes vol_{Y/X}^{-1}.
\end{eqnarray}
On the other hand we get
\begin{eqnarray}
"(f\boxtimes Id_L)_*\circ pr_X^*"(\phi)(K)=\\
\int_{z\in Y/X}(pr_X^*\phi)(K\cap (z\times L))dvol_{Y/X}\otimes
vol_{Y/X}^{-1}=\\\label{ob7} \int_{z\in Y/X}\phi(pr_X(K\cap
(z\times L)))dvol_{Y/X}\otimes vol_{Y/X}^{-1}
\end{eqnarray}
Comparing (\ref{ob6}) and (\ref{ob7}) and making appropriate
identifications of subsets of $X$ inside $Y$, we conclude Case b.
Thus Step 3 is completed.

\underline{Step 4.} Assume that $g$ is injective.

Then we may and will assume that $\tilde Y\subset Y$ and $g\colon
\tilde Y\to Y$ is the identity imbedding. Let us fix Lebesgue
measures $vol_{\tilde Y}$ on $\tilde Y$ and $vol_Y$ on $Y$. Set
$$vol_{Y/\tilde Y}:=\frac{vol_Y}{vol_{\tilde Y}}\in Dens(Y/\tilde
Y).$$ By Step 1, it suffices to prove the result in two cases:
either $f$ is surjective or injective.

\underline{Case a.} Assume that $f$ is surjective.

Then we may assume that $X=Y\oplus M$ and $f=pr_Y\colon Y\oplus
M\to Y$ is the natural projection. Then the diagram (\ref{ob0})
becomes equal to
\begin{eqnarray}\label{ob8}
\square<1`1`1`1;800`400>[ \tilde Y\oplus M`\tilde Y`Y\oplus M
`Y;pr_{\tilde Y}`\tilde g=g\boxtimes Id_M`g`pr_Y].
\end{eqnarray}
Let us fix a Lebesgue measure $vol_M$ on $M$, and $T\in \ck(M)$
such that $vol_M(T)=1$. Set $vol_X:=vol_Y\otimes vol_M\in
Dens(Y\oplus M)=Dens(X)$. Let $\phi\in Val(Y\oplus M)\otimes
Dens(\tilde Y\oplus M)$. Let us fix also $K\in \ck(\tilde Y)$.

We have
\begin{eqnarray}
("g^*\circ pr_{Y*}")(\phi)(K)=("pr_{Y*}"(\phi))(K)=\\\label{ob9}
\frac{1}{m!}\frac{d^m}{d\eps^m}\big|_0\phi(K\times\eps T).
\end{eqnarray}
On the other hand we have
\begin{eqnarray}
("pr_{\tilde Y*}\circ \tilde
g^*")(\phi)(K)=\frac{1}{m!}\frac{d^m}{d\eps^m}\big|_0("\tilde
g^*"\phi)(K\times \eps T)=\\\label{ob10}
\frac{1}{m!}\frac{d^m}{d\eps^m}\big|_0\phi(K\times \eps T).
\end{eqnarray}
Comparing (\ref{ob9}) and (\ref{ob10}) we conclude Case a.

\underline{Case b.} Assume that $f$ is injective.

We may and will assume that $X\subset Y$ and $f$ is the identity
imbedding. Under these assumptions $\tilde X=X\cap \tilde Y$. Let
us fix decompositions
\begin{eqnarray*}
X=\tilde X\oplus L,\\
\tilde Y=\tilde X\oplus M.
\end{eqnarray*}
Then $Y=\tilde X\oplus L\oplus M$ (since we assume that $f\oplus
g\colon X\oplus \tilde Y\to Y$ is onto). Then the diagram
(\ref{ob0}) becomes the following diagram of imbeddings
\begin{eqnarray}\label{ob11}
\square<1`1`1`1;1100`400>[ \tilde X`\tilde X\oplus M`\tilde
X\oplus L `\tilde X\oplus L\oplus M;Id_{\tilde X}\boxtimes
0_M`Id_{\tilde X}\boxtimes 0_L`Id_{\tilde X\oplus M}\boxtimes
0_L`Id_{\tilde X\oplus L}\boxtimes 0_M].
\end{eqnarray}
Let us fix Lebesgue measures $vol_{\tilde X}, vol_L,vol_M$ on
$\tilde X,L,M$ respectively. Let
\begin{eqnarray*}
vol_X:=vol_{\tilde X}\otimes vol_L\in Dens(X),\\
vol_{\tilde Y}:=vol_{\tilde X}\otimes vol_M\in Dens(\tilde Y),\\
vol_Y:=vol_{\tilde X}\otimes vol_L\otimes vol_M\in Dens(Y).
\end{eqnarray*}
Fix $\phi\in Val(X)\otimes Dens(\tilde X^*)$. Fix $K\in \ck(\tilde
Y)$. We have
\begin{eqnarray}
("g^*\circ f_*")(\phi)(K)=("f_*"\phi)(K)=\\ \int_{x\in
M}\phi(K\cap (x+(\tilde X\oplus L)))dvol_M(x)\otimes vol_M^{-1}
=\\\label{ob12} \int_{x\in M}\phi((K\cap (x+\tilde X))\times
\{0_L\}) \otimes vol_M^{-1}.
\end{eqnarray}
On the other hand we have
\begin{eqnarray}
("\tilde f_*\circ \tilde g^*")(\phi)(K)=\\ \int_{x\in M} ("\tilde
g^*"\phi)(K\cap (x+\tilde X))dvol_M(x)\otimes
vol_M^{-1}=\\\label{ob13} \int_{x\in M}\phi((K\cap (x+\tilde
X))\times \{0_L\})\otimes vol_M^{-1}.
\end{eqnarray}
Comparing (\ref{ob12}) and (\ref{ob13}) we conclude Case b. Thus
Step 4 is proved.

\underline{Step 5.} Let us consider the general case. By Step 2 we
may assume that $g$ is either  injective of surjective. Now the
theorem follows from Steps 3,4. \qed


The next result is in fact an equivalent reformulation of the base
change theorem. It will be needed later.
\begin{theorem}\label{T:base-change2}
Let
$$\square[\tilde X`\tilde Y`X`Y;\tilde f`\tilde g`g`f]$$
be a Cartesian square of vector spaces such that $f\oplus g\colon
X\oplus \tilde Y\to Y$ is onto. Consider the following two maps
$$Val(\tilde Y)\otimes Dens(\tilde Y^*)\to Val(X)\otimes
Dens(Y^*).$$ The first map  is the composition
\begin{eqnarray*}
Val(\tilde Y)\otimes Dens(\tilde
Y^*)\overset{g_*}{\to}Val(Y)\otimes Dens(Y^*)\overset{f^*\otimes
Id_{Dens(Y^*)}}{\to}Val(X)\otimes Dens(Y^*);
\end{eqnarray*}
and the second map is the composition
\begin{eqnarray*}
Val(\tilde Y)\otimes Dens(\tilde Y^*)\overset{\tilde f^*\otimes
Id_{Dens(\tilde Y^*)}}{\to}\\Val(\tilde X)\otimes Dens(\tilde
Y^*)=\left(Val(\tilde X)\otimes Dens(\tilde X^*)\right)\otimes
\frac{Dens(\tilde Y^*)}{Dens(\tilde X^*)}\overset{\tilde
g_*\otimes Id}{\to}\\\left(Val(X)\otimes
Dens(X^*)\right)\otimes\frac{Dens(\tilde Y^*)}{Dens(\tilde
X^*)}=Val(X)\otimes Dens(Y^*)
\end{eqnarray*}
where in the last equality we have used the identification from
Lemma \ref{L:bch}(i).

Then these two maps coincide.
\end{theorem}
{\bf Proof.} This result is obtained from Theorem
\ref{T:base-change} by flipping the diagram in the latter theorem
with respect to the diagonal, twisting all the spaces by
$\frac{Dens(\tilde Y^*)}{Dens(\tilde X^*)}$, and using the
isomorphism $\frac{Dens(\tilde Y^*)}{Dens(\tilde X^*)}\otimes
Dens(X^*)\simeq  Dens(Y^*)$ from Lemma \ref{L:bch}(i). \qed

\section{An isomorphism of $GL(V)$-modules $Val^{-,sm}_{n-p}(V)$
and $Val_p^{-,sm}(V^*)\otimes Dens(V)$}\label{S:isomorphism-val}
\setcounter{subsection}{1}
\setcounter{theorem}{0}\setcounter{equation}{0} The main result of
this section is the following proposition.
\begin{proposition}\label{isomorphism}
Let $1\leq p\leq n-1$. The $GL(V)$-modules $Val^{-,sm}_{n-p}(V)$
and $Val_p^{-,sm}(V^*)\otimes Dens(V)$ are isomorphic.
\end{proposition}
{\bf Proof.} \underline{Step 1.} By \cite{alesker-gafa-04} the
product on valuations
$$Val_{n-p}^{-,sm}\times Val_p^{-,sm}(V)\to Val_n(V)=Dens(V)$$
is a perfect pairing. It follows that the induced map
\begin{eqnarray}\label{1}
Val_{n-p}^{-,sm}\to (Val_p^{-,sm}(V))^{*,sm}\otimes Dens(V)
\end{eqnarray}
is an isomorphism of $GL(V)$-modules. Thus to prove the
proposition we have to show that the $GL(V)$-modules
$(Val_p^{-,sm}(V))^{*,sm}$ and $Val_p^{-,sm}(V^*)$ are isomorphic.

\underline{Step 2.} For a vector space $W$ let us denote by
$\cf_{k,k+1}(W)$ the manifold of partial flags
$$\cf_{k,k+1}(W):=\{(E,F)\,\,\,| \,\,\, E\in Gr_k(W),\, F\in
Gr_{k+1}(W),\, E\subset F\}.$$ Let
$\cm_{k,k+1}(W)\to\cf_{k,k+1}(W)$ be the vector bundle such that
its fiber over a pair $(E\subset F)\in \cf_{k,k+1}(W)$ is equal to
$E$. Similarly let $\cn_{k,k+1}(W)\to\cf_{k,k+1}(W)$ be the vector
bundle such that its fiber over $(E\subset F)\in \cf_{k,k+1}(W)$
is equal to $F$. Thus $\cm_{k,k+1}(W)\subset \cn_{k,k+1}(W)$.

By \cite{alesker-adv-00}, $Val_k^{-,sm}(W)$ is isomorphic to an
irreducible subquotient of $$C^\infty(\cf_{k,k+1}(W),\det
\cm_{k,k+1}^*(W)\otimes or (\cn_{k,k+1}(W))).$$ Thus
$Val_p^{-,sm}(V^*)$ is isomorphic to an irreducible subquotient of
$$C^\infty(\cf_{p,p+1}(V^*), \det \cm_{p,p+1}^*(V^*)\otimes or
(\cn_{p,p+1}(V^*))).$$

\underline{Step 3.} Now we will show that both
$(Val_p^{-,sm}(V))^{*,sm}$ and $Val_p^{-,sm}(V^*)$ appear in the
Jordan-H\"older series of the same degenerate principal series
representation.

Step 2 implies that $(Val_p^{-,sm}(V))^{*,sm}$ is isomorphic to an
irreducible subquotient of
\begin{eqnarray}\label{2}
C^\infty(\cf_{p,p+1}(V),\det \cm_{p,p+1}(V)\otimes or(
\cn_{p,p+1}(V))\otimes |\ome_{\cf_{p,p+1}(V)}|)
\end{eqnarray}
where $|\ome_X|$ denotes the line bundle of densities over a
manifold $X$.

Recall that by Step 2 $Val_p^{-,sm}(V^*)$ is isomorphic to an
irreducible subquotient of $$C^\infty(\cf_{p,p+1}(V^*), \det
\cm_{p,p+1}^*(V^*)\otimes or (\cn_{p,p+1}(V^*))).$$

By taking the orthogonal complement, $\cf_{p,p+1}(V^*)$ is
identified with $\cf_{n-p-1,n-p}(V)$. Then $\cm^*_{p,p+1}(V^*)$ is
identified with $\underline{V}/\cn_{n-p-1,n-p}(V)$, and
$\cn_{p,p+1}(V^*)$ is identified with the bundle
$(\underline{V}/\cm_{n-p-1,n-p}(V))^*$ where
$\underline{V}=\cf_{n-p-1,n-p}(V)\times V$. Of course, all the
identifications are $GL(V)$-equivariant. Hence
\begin{eqnarray}\label{2.1}
C^\infty(\cf_{p,p+1}(V^*), \det \cm_{p,p+1}^*(V^*)\otimes or
(\cn_{p,p+1}(V^*)))=\\ \label{3}
C^\infty\left(\cf_{n-p-1,n-p}(V),\det(\underline{V}/\cn_{n-p-1,n-p}(V))
\otimes or (\underline{V}/\cm_{n-p-1,n-p}(V))\right).
\end{eqnarray}

By Corollary \ref{C:arepr2} the natural representations of $GL(V)$
in the spaces (\ref{2}) and (\ref{3}) have the same
Jordan-H\"older series. Hence both $(Val_p^{-,sm}(V))^{*,sm}$ and
$Val_p^{-,sm}(V^*)$ appear in the Jordan-H\"older series of
(\ref{3}) which is isomorphic to (\ref{2.1}), i.e. to
\begin{eqnarray}\label{4}
\cx:=C^\infty(\cf_{p,p+1}(V^*), \det \cm_{p,p+1}^*(V^*)\otimes or
(\cn_{p,p+1}(V^*))).
\end{eqnarray}

\underline{Step 4.} In this last and the most technical step we
will show that $(Val_p^{-,sm}(V))^{*,sm}$ cannot be isomorphic to
any constituent of the Jordan-H\"older series of (\ref{4})
different from $Val_p^{-,sm}(V^*)$. This will finish the proof of
the proposition.
\def\ttg{\texttt{g}}
Remind that to any finitely generated $U(\ttg)$-module $M$ one can
attach an algebraic subvariety of the variety of nilpotent element
of $\ttg$ which is called {\itshape associated variety} or
{\itshape Bernstein variety}. We refer to \cite{borho-brylinski1}
for the detail on this notion. It turns out that the associated
variety of $(Val_p^{-,sm}(V))^{*,sm}$ is equal to the variety of
complex symmetric nilpotent matrices of rank at most 1. Indeed by
the Poincar\'e duality (\ref{1}) this space is isomorphic to
$Val_{n-p}^{-,sm}(V)\otimes Dens(V)^*$. Clearly the associated
variety of the last space coincides with that of
$Val_{n-p}^{-,sm}(V)$. But by \cite{alesker-gafa-01}, Theorem 3.1,
the associated variety of $Val_{n-p}^{-,sm}(V)$ is equal to the
variety of complex symmetric nilpotent matrices of rank at most 1.

Let us remind few facts about the structure of the space
(\ref{4}). Let
$$q\colon \cf_{p,p+1}(V^*)\to Gr_{p+1}(V^*)$$ be the canonical
projection. Let $\ct_{p+1}(V^*)\to Gr_{p+1}(V^*)$ be the
tautological bundle, i.e. the bundle whose fiber over $F\in
Gr_{p+1}(V^*)$ is equal to $F$. It is clear that
$\cn_{p,p+1}(V^*)=q^*(\ct_{p+1}(V^*))$. Hence
\begin{eqnarray*}
\wedge^p\cn_{p,p+1}^*(V^*)=q^*(\wedge^p\ct_{p+1}^*(V^*)),\\
or(\cn_{p,p+1}(V^*))=q^*(or(\ct_{p+1}(V^*)).
\end{eqnarray*}
Consider the map of vector bundles $\cn_{p,p+1}^*(V^*)\to
\cm_{p,p+1}^*(V^*)$ dual to the natural imbedding
$\cm_{p,p+1}(V^*)\inj \cn_{p,p+1}(V^*)$. The $p$-th exterior power
of this map induces a map
\begin{eqnarray*}
\cy:=C^\infty(Gr_{p+1}(V^*),\wedge^p\ct_{p+1}^*(V^*)\otimes
or(\ct_{p+1}(V^*)))\to\\
C^\infty(\cf_{p,p+1}(V^*),\det\cm_{p,p+1}^*(V^*)\otimes
or(\cn_{p,p+1}(V^*)))=\cx.
\end{eqnarray*}
Clearly this map is $GL(V)$-equivariant (and in fact injective).
The Casselman-Wallach theorem implies that the image of this map
is a closed subspace. We will identify $\cy$ with its image in
$\cx$ under this map. By \cite{alesker-adv-00},
$Val_p^{-,sm}(V^*)$ imbeds $GL(V)$-equivariantly as a subspace in
$\cx/\cy$. Moreover, by \cite{alesker-gafa-01}, Section 5,
$Val_p^{-,sm}(V^*)$ is the only irreducible subquotient of
$\cx/\cy$ whose associated variety consists of complex symmetric
nilpotent matrices of rank at most 1.

Hence it remains to show that $(Val_p^{-,sm}(V))^{*,sm}$ cannot be
isomorphic to any of the irreducible subquotients of
$\cy=C^\infty(Gr_{p+1}(V^*),\wedge^p\ct_{p+1}^*(V^*)\otimes
or(\ct_{p+1}(V^*)))$.

Before we will treat the general case, let us observe now that if
$p+1=n$ then the last statement is trivial since
$(Val_p^{-,sm}(V))^{*,sm}$ is infinite dimensional while
$\cy=\wedge^pV^*\otimes or(V^*)$ is finite dimensional (since
$Gr_{p+1}(V^*)$ is just a point). Hence for $p=n-1$ Proposition
\ref{isomorphism} is proved. By symmetry, replacing $V$ by $V^*$,
Proposition \ref{isomorphism} follows also for $p=1$.

Let us assume now that $2\leq p\leq n-2$, hence $n\geq 4$. In this
case we are going to use the Beilinson-Bernstein localization
theorem. First we will have to introduce more notation and remind
some constructions from \cite{alesker-gafa-01}.

Let us fix a Euclidean metric on $V^*$. Let $G_0=GL(V^*)$. Let
$\texttt{g}_0=Lie(G_0)$ be the Lie algebra of $G_0$. Let
$\ttg:=\ttg_0\otimes _\RR\CC$ be its complexification. Let
$K_0\subset G_0$ be the subgroup of the orthogonal transformations
of $V^*$. Let $K$ be the complexification of $K_0$. Let $G$ be the
complexification of $G_0$. Thus $G\simeq GL_n(\CC)$,
$Lie(G)=\ttg$.

Let $\grc_k$ denote the Grassmannian of complex $k$-dimensional
subspaces of $V^*\otimes_\RR\CC=:\vc$. Denote
$$\ffc_{p,p+1}:=\{(E,F)|\, E\subset F,E\in\grc_p,F\in
\grc_{p+1}\}.$$ Let $\ffc$ be the variety of complete flags in
$\vc^*$. We have the canonical projection
\begin{eqnarray}\label{obe1}
\bar q\colon \ffc\to \ffc_{p,p+1}.
\end{eqnarray}
It is well known that the group $K$ acts on $\ffc$ (and hence on
$\ffc_{p,p+1},\grc_p$) with finitely many orbits.

Let us fix a basis $e_1,e_2,\dots,e_n$ in $V^*$. Let $T_0\subset
G_0$ be the subgroup of diagonal transformations with respect to
this basis. Let $B_0\subset G_0$ be the subgroup of upper
triangular transformations. Let $T$ and $B$ be the
complexifications of $T_0$ and $B_0$ respectively. Thus $T\subset
G$ is a Cartan subgroup, $B\subset G$ is a Borel subgroup.

Let $P_0\subset G_0$ be the subgroup of transformations of
preserving the flag $span_\RR\{e_1,\dots,e_p\}\subset
span_\RR\{e_1,\dots,e_p,e_{p+1}\}$. Let $P\subset G$ be its
complexification which is a parabolic subgroup of $G$. Then
clearly $T_0\subset B_0\subset P_0$, $T\subset B\subset P$.
\def\lit{\texttt{t}}
\def\lib{\texttt{b}}
\def\lip{\texttt{p}}
By $\lit,\lib,\lip$ we will denote the Lie algebras of $T,B,P$
respectively.

Clearly in the basis $e_1,\dots,e_n$ the subgroup $T$ consists of
complex diagonal invertible matrices, $B$ consists of complex
upper triangular invertible matrices, and
\begin{eqnarray*}
P=\left\{ \left[\begin{array}{ccc}
                             A&*&*\\
                             0&b&*\\
                             0&0&C
                \end{array}\right]\big|\,
                A\in GL_p(\CC),b\in \CC^*,C\in
                GL_{n-p-1}(\CC)\right\}.
\end{eqnarray*}
Let us consider the character $\chi\colon \lip/[\lip,\lip]\to \CC$
given by
\begin{eqnarray}\label{obe2}
\chi\left(\left[\begin{array}{ccc}
                             A&*&*\\
                             0&b&*\\
                             0&0&C
                \end{array}\right]\right)=-Tr(A).
\end{eqnarray}
Let $\hat\chi\colon \lib/[\lib,\lib]\to \CC$ be the composition of
$\chi$ with the canonical map
$\lib/[\lib,\lib]\to\lip/[\lip,\lip]$. It is easy to see that
\begin{eqnarray*}
\hat\chi\left(\left[\begin{array}{cccccc}
                                 x_1&*&&&&*\\
                                 0&\ddots&*&&&\\
                                 0&\dots&x_p&*&&\\
                                 0&\dots&0&x_{p+1}&&\\
                                 0&\dots&\dots&0&\ddots&*\\
                                 0&\dots&\dots&\dots&0&x_n
                                 \end{array}\right]\right)=-(x_1+\dots+x_p).
\end{eqnarray*}
\def\hchi{\hat\chi}

We will be interested in $\cd_\chi$-modules on $\ffc_{p,p+1}$ and
$\cd_{\hchi}$-modules on $\ffc$. Clearly we have the pullback
functor
\begin{eqnarray*}
\bar q^*\colon \cd_\chi(\ffc_{p,p+1})-mod\to\cd_{\hchi}(\ffc)-mod.
\end{eqnarray*}
\begin{lemma}\label{obe3}
The character $\hchi\colon \lib/[\lib,\lib]\to\CC$ is dominant and
regular in sense of Definition \ref{D:regular-dominant}.
\end{lemma}
{\bf Proof.} Let $\bar \chi\colon B_0/[B_0,B_0]\to \CC^*$ be the
character of the group defined by
\begin{eqnarray*}
\bar\chi\left(\left[\begin{array}{ccc}
                                  z_1&*&*\\
                                  0&\ddots&*\\
                                  0&\dots&z_n
                                  \end{array}\right]\right)=
                                  (z_1\dots z_p)^{-1}
\end{eqnarray*}
Consider the representation $Ind_{B_0}^{G_0}\bar\chi$. By
Proposition \ref{P:reg-dom}, in order to prove that $\hchi$ is
dominant and regular it is enough to show that
$Ind_{B_0}^{G_0}\bar\chi$ has a non-zero finite dimensional
submodule. But we have a natural non-zero map $\wedge^pV\to
Ind_{B_0}^{G_0}\bar\chi$. Hence lemma is proved. \qed

\begin{corollary}\label{obe4}
The functor of global sections
\begin{eqnarray*}\Gamma\colon
\cd_\chi(\ffc_{p,p+1})-mod \to \texttt{g} -mod\end{eqnarray*} is
exact and faithful.
\end{corollary}
{\bf Proof.} Since the morphism $\bar q$ is projective and smooth
$$\Gamma=\hat\Gamma\circ \bar q^*$$
where $\hat\Gamma$ is the functor of global sections on the
category $\cd_{\hchi}(\ffc)-mod$. $\hat\Gamma$ is exact and
faithful by Lemma \ref{obe3} and the Beilinson-Bernstein theorem.
The functor $\bar q^*$ is exact and faithful too since $\bar q$ is
a smooth morphism. Hence $\Gamma$ is also exact and faithful. \qed

Let us now remind, following \cite{alesker-gafa-01}, the
construction of a $K$-equivariant $\cd_\chi$-module $\cm$ on
$\ffc_{p,p+1}$ such that the space of its global sections is
isomorphic, as a $(\ttg,K)$-module, to the Harish-Chandra module
of $\cx$ (defined in (\ref{4})).

Let $B$ be the complexification of the Euclidean form on $V^*$.
Thus $B\colon \vc\times\vc\to \CC$ is a symmetric non-degenerate
bilinear form. Let $U\subset \ffc_{p,p+1}$ denote the open
$K$-orbit of $\ffc_{p,p+1}$. Let $j\colon U\inj \ffc_{p,p+1}$
denote the identity imbedding. Explicitly one has
\begin{eqnarray*}
U=\{(E,F)\in \ffc_{p,p+1}|\, \mbox{the restrictions of } B \mbox{
to } E \mbox{ and to } F \mbox{ are non-degenerate}\}.
\end{eqnarray*}
Let us fix an element $(E_0,F_0)\in U$. The stabilizer $S\subset
K$ of this element is isomorphic to the group $O(p,\CC)\times
O(1,\CC)\times O(n-p-1,\CC)$. Note that $S$ is a reductive group,
and hence $U=K/S$ is an affine variety.

The category of $K$-equivariant $\cd_\chi$-modules on the orbit
$U$ is equivalent to the category of representations of the group
of connected components of $S$. Let $\cm_0$ be the $K$-equivariant
$\cd_\chi$-module on $U$ corresponding to the representation of
$S\simeq O(p,\CC) \times O(1,\CC)\times O(n-p-1,\CC)$ given by
$$(A,B,C)\mapsto \det A\cdot \det B.$$
Let us define $\cm:=j_*\cm_0$. It was shown in
\cite{alesker-gafa-01} that the $(\ttg,K)$-module
$\Gamma(\ffc_{p,p+1}\cm)$ is isomorphic to the Harish-Chandra
module of $\cx$.

Let us describe now the $K$-equivariant sub-$\cd_\chi$-module
$\cn\subset\cm$ corresponding to the Harish-Chandra module of
$\cy$ (more precisely, $\Gamma(\ffc_{p,p+1},\cn)$ coincides with
the Harish-Chandra module of $\cy\subset \cx$).

Let us denote by $V\subset \ffc_{p,p+1}$ the (open) subvariety
$$V:=\{(E,F)\in \ffc_{p,p+1}|\mbox{ s.t. } B|_F \mbox{ is
non-degenerate}\}.$$ Then $U\subset V$. Let $j'\colon U\inj V$ and
$j''\colon V\inj \ffc_{p,p+1}$ denote the identity imbedding
morphisms. Set
$$\cn:=j''_*(j'_{!*}\cm_0)$$ where $j'_{!*}$ denotes the minimal
(Goresky-Macpherson) extension of $\cm$ under the open imbedding
$j'$. It is easy to see that the morphism $j''\colon V\to
\ffc_{p,p+1}$ is affine. Hence the functor $j''_*$ is exact and we
have
$$\cn\subset j''_*(j'_{*}\cm_0)=j_*\cm_0=\cm.$$
As it was shown in \cite{alesker-gafa-01},
$\Gamma(\ffc_{p,p+1},\cn)$ coincides with the Harish-Chandra
module of $\cy\subset \cx$. Also let us define
$$\ck:=j_{!*}\cm_0.$$
Thus $\ck\subset \cn$.
\begin{remark}\label{ny1}
Since $Val_p^{-,sm}(V^*)\inj \cx/\cy$, it corresponds to a
sub-$\cd_\chi$-module of $\cm/\cn$. Note also that $\supp
(\cm/\cn)\ne \ffc_{p,p+1}$.
\end{remark}
\begin{lemma}\label{ny2}
No irreducible subquotient of
$\Gamma(\cn/\ck)=\Gamma(\cn)/\Gamma(\ck)$ has associated variety
contained in the variety of complex symmetric nilpotent matrices
of rank at most 1.
\end{lemma}
{\bf Proof.} Let $\cl\to \ffc_{p,p+1}$ be the (algebraic) line
bundle whose fiber over $(E,F)\in \ffc_{p,p+1}$ is equal to
$\wedge^pE^*$. Then the sheaf $\cd_\chi$ is the sheaf of
differential operators with values in $\cl$. Tensoring by $\cl^*$
establishes an equivalence of the categories
$\cd_\chi(\ffc_{p,p+1})-mod$ and $\cd(\ffc_{p,p+1})-mod$ where
$\cd(X)$ denotes the sheaf of rings of usual (untwisted)
differential operators on a variety $X$. Since the line bundle
$\cl$ is $G$-equivariant, the analogous equivalence holds for
$K$-equivariant versions of these categories. This equivalence
preserves singular supports of the corresponding modules.

Let $\car$ be a $\cd_\chi$-module. Let $\car':=\cl^*\otimes \car$
be the corresponding $\cd$-module. The associated variety of
$\Gamma(\car,\ffc_{p,p+1})$ (resp. $\Gamma(\car',\ffc_{p,p+1})$)
is equal to the image under the moment map of the singular support
of $\car$ (resp. $\car'$) (see \cite{alesker-gafa-01} where the
discussion follows \cite{borho-brylinski2}). Hence it follows that
$\Gamma(\car,\ffc_{p,p+1})$ and $\Gamma(\car',\ffc_{p,p+1})$ have
the same associated variety.

Set $\cm_0':=\cl^*\otimes \cm_0$. Then $\cm_0'$ is a
$K$-equivariant $\cd$-module on the open orbit $U\subset
\ffc_{p,p+1})$ corresponding to the same as $\cm_0$ representation
of the group of connected components of the group $S\simeq
O(p,\CC)\times O(1,\CC)\times O(n-p-1,\CC)$, i.e. to
$(A,B,C)\mapsto \det A \cdot \det B$.

Then
\begin{eqnarray*}
\cm':=\cl^*\otimes \cm=j_*\cm_0',\\
\cn':=\cl^*\otimes \cn=j''_*(j'_{!*}\cm_0')
\end{eqnarray*}
where now all the functors $j_*, j'_{!*},j''_*$ are in the
category of $\cd$-modules rather than $\cd_\chi$-modules.

Let $f\colon \ffc_{p,p+1}\to \grc_{p+1}$ be the canonical
projection $(E,F)\mapsto F$. Let us denote by $\co$ the open
$K$-orbit in $\grc_{p+1}$.  Explicitly
$$\co=\{F\in \grc_{p+1}|\, B|_F \mbox{ is non-degenerate}\}.$$
It is clear that $V=f^{-1}(\co)$.

The stabilizer $S'\subset K$ of a point from $\co$ is isomorphic
to the group $O(p+1,\CC)\times O(n-p-1,\CC)$. Let us denote by
$\ca$ the $\cd$-module on $\co$ corresponding to the
representation of the group of connected components of $S'\simeq
O(p+1,\CC)\times O(n-p-1,\CC)$ given by $(M,N)\mapsto \det M$.

\begin{claim}\label{ny3}
The $K$-equivariant $\cd$-modules $f^*\ca$ and $j'_{!*}\cm_0'$ on
$V=f^{-1}(\co)$ are isomorphic.
\end{claim}
{\bf Proof.} It is clear from the definitions of $\ca$ and
$\cm_0'$ that the restriction of $f^*\ca$ to $U$ is isomorphic to
$\cm_0'$. Hence we have a morphism of $K$-equivariant
$\cd$-modules $f^*\ca\to j'_*\cm_0'$ which is an isomorphism over
$U$. It is clear that $f^*\ca$ is an irreducible $\cd$-module,
hence its image is equal to $j_{!*}'\cm_0'$. \qed

Let $l\colon \co\inj \grc_{p+1}$ denote the identity imbedding
morphism. Since the morphism $f$ is smooth (in particular flat) by
the flat base change theorem (see e.g \cite{hartshorne},
Proposition 9.3) we have
\begin{eqnarray}\label{ny4}
\cn':=\cl^*\otimes \cn=j''_*(j'_{!*}\cm_0')\simeq
j''_*(f^*\ca)=f^*(l_*\ca).
\end{eqnarray}
Set $\ck':=\cl^*\otimes \ck=j_{!*}\cm_0'$. Since
$j_{!*}=j''_{!*}\circ j'_{!*}$ we have
\begin{eqnarray}\label{ny5}
\ck'\simeq j''_{!*}(f^*\ca)=f^*(l_{!*}\ca)
\end{eqnarray}
where the last equality also follows from the smooth base change
theorem for $\cd$-modules \cite{??}.

Since the morphism $f$ is projective and smooth, one has
$\Gamma(\ffc_{p,p+1},f^*(\bullet))=\Gamma(\grc_{p+1},\bullet)$.
Then
$\Gamma(\ffc_{p,p+1},\cn'/\ck')=\Gamma(\ffc_{p,p+1},f^*(l_*\ca/l_{!*}\ca))=
\Gamma(\grc_{p+1},l_*\ca/l_{!*}\ca)$. Thus to finish the proof of
Lemma \ref{ny2} it remains to show that the associated variety of
any irreducible $(\ttg,K)$-subquotient of $\Gamma(\grc_{p+1},
l_*\ca/l_{!*}\ca)$ is {\itshape not} contained in the variety of
complex symmetric nilpotent matrices of rank at most 1. But this
statement was proved in fact in Theorem 4.3 in
\cite{alesker-gafa-01}. \qed

To finish the proof of Proposition \ref{isomorphism} it remains to
show that the Harish-Chandra module of $(Val_p^{-,sm}(V))^*$
cannot be isomorphic to the $(\ttg, K)$-module
$\Gamma(\ffc,j_{!*}\cm_0)$.

Since $\bar q\colon \ffc\to \ffc_{p,p+1}$ is a smooth projective
morphism we have
\begin{eqnarray}\label{ny6}
\Gamma(\ffc,\bar q^*(\bullet))=\Gamma(\ffc_{p,p+1},\bullet).
\end{eqnarray}
Clearly $\bar q^*(j_{!*}\cm_0)$ is an irreducible
$\cd_{\hat\chi}$-module with support
\begin{eqnarray}\label{ny6.1}
\supp(\bar q^*(j_{!*}\cm_0))=\ffc.
\end{eqnarray}
It is clear from (\ref{ny6}) and from the fact that $\hat\chi$ is
dominant and regular, that $\bar q^*(j_{!*}\cm_0)$ is the
Beilinson-Bernstein localization of
$\Gamma(\ffc_{p,p+1},\bullet)$. Thus it suffices to prove the
following claim.
\begin{claim}\label{C:ny7}
The support on $\ffc$ of the Beilinson-Bernstein localization of
the Harish-Chandra module of $(Val_p^{-,sm}(V))^*$ is not equal to
$\ffc$.
\end{claim}
{\bf Proof.} Remind that $\ffc$ denotes the variety of complete
flags in $V^*\otimes_\RR\CC$. By taking the orthogonal complement
we can identify $\ffc$ with the variety $\ffc'$ of complete flags
in $V\otimes_\RR\CC$. Also the group $G_0=GL(V^*)$ can be
identified with the group $G_0'=GL(V)$ via the isomorphism
$g\mapsto (g^t)^{-1}$. Let $G'$ denote the complexification of
$G_0'$. Now we will work with $G_0',G'$ instead of $G_0,G$, and we
will consider the Beilinson-Bernstein localization on $\ffc'$. Let
us denote similarly $\ttg':=Lie(G')$, $K'\subset G'$ is the
subgroup preserving the non-degenerate form on $V\otimes_\RR\CC$.
Let us denote
$$\ffc'_{k,k+1}:=\{(E,F)|\, E\subset F, E\in
\grc_k(V\otimes_\RR\CC),F\in\grc_{k+1}(V\otimes_\RR\CC)\}.$$
Remind that we have the Poincar\'e duality isomorphism
$(Val_p^{-,sm}(V))^{*,sm}\simeq Val_{n-p}^{-,sm}(V)\otimes
Dens(V^*)$. Thus it suffices to prove that the support in $\ffc'$
of the Beilinson-Bernstein localization of $Val_{n-p}^{-,sm}(V)$
is not equal to $\ffc'$.

By Remark \ref{ny1} applied to $V$ instead of $V^*$ and $n-p$
instead of $p$, the Harish-Chandra module of $Val_{n-p}^{-,sm}(V)$
is isomorphic to the $(\ttg',K')$-module of global section of
certain $K'$-equivariant $\cd_{\bar \chi}$-module $\tau$ on
$\cf_{n-p,n-p+1}'$ whose support is not equal to
$\ffc'_{n-p,n-p+1}$. Here we denote by $\bar \chi$ an appropriate
character of a stabilizer of a Lie algebra $Q$ of a point from
$\ffc'_{n-p,n-p+1}$. Note that the lifting $\bar{\bar\chi}$ of
$\bar \chi$ to the Cartan algebra $t'\subset \ttg'$ is dominant
and regular.

Let $g\colon \ffc'\to \ffc'_{n-p,n-p+1}$ denote the canonical
projection. Then we have the functor  $g^*\colon
\cd_{\bar\chi}(\ffc'_{n-p,n-p+1})-mod\to
\cd_{\bar{\bar\chi}}(\ffc')-mod$. Since $g$ is projective and
smooth we have
$$\Gamma(\ffc',g^*(\bullet))=\Gamma(\ffc'_{n-p,n-p+1},\bullet).$$
Hence the Beilinson-Bernstein localization of the Harish-Chandra
module of $Val^{-,sm}_{n-p}(V)$ is isomorphic to $g^*\tau$. Since
$\supp \tau\ne \ffc'_{n-p,n-p+1}$ then $\supp g^*\tau\ne \ffc'$.
Thus Claim \ref{C:ny7} is proved. \qed

Hence Proposition \ref{isomorphism} is proved as well. \qed

\begin{remark}\label{R:mult1thm}
Note that in Step 4 of the proof of Proposition \ref{isomorphism}
we have actually proven the following result: for any
$n$-dimensional vector space $W$ and any $p=1,\dots,n-1$ the
Jordan-H\"older series of the $GL(W)$-module
$C^\infty\left(\cf_{p,p+1}(W),\det\cm^*_{p,p+1}(W)\otimes
or(\cn_{p,p+1}(W))\right)$ contains $Val^{-,sm}_{p}(W)$ with
multiplicity one. This fact will be used below in the proof of
Proposition \ref{P:w10}.
\end{remark}

\def\unv{\underline{V}}
\def\tv{\ct(V)}
\def\tvs{\ct(V^*)}


\section{The two dimensional case.}\label{2-dim} 
The goal of this section is to construct the Fourier transform on
two-dimensional spaces and to establish its homomorphism property.
Let $V$ be a {\itshape two dimensional} real vector space.
\subsection{A canonical isomorphism
in two dimensions.}  First we are going to construct a {\itshape
canonical} isomorphism
\begin{eqnarray}\label{can-isomor}
Val_1^{-,sm}(V)\tilde\to Val_1^{-,sm}(V^*)\otimes \det V^*.
\end{eqnarray}
Let $\ct(V)\to \PP(V)$ be the tautological vector bundle, i.e. the
fiber of $\ct(V)$ over $E\in \PP(V)$ is equal to $E$. Let
$\underline{V}\to \PP(V)$ be the trivial bundle, i.e.
$\underline{V}=\PP(V)\times V$. Then we have a canonical imbedding
of vector bundles $\ct(V)\inj \underline{V}$.

We have the canonical epimorphism
\begin{eqnarray}\label{2-1}
C^\infty(\PP(V),Dens\tv\otimes or(\unv/\tv))\surj Val_1^{-,sm}(V).
\end{eqnarray}
Recall that the kernel of this map is two dimensional irreducible
$GL(V)$-module. Observe that
\begin{eqnarray*}
Dens \tv=\tv^*\otimes or(\tv),\\
or(\unv/\tv)=or(\unv)\otimes or(\tv).
\end{eqnarray*}
Using these identifications we get a canonical epimorphism
\begin{eqnarray}\label{2-2}
C^\infty(\PP(V),\det \tv^*)\otimes or(V)\surj Val_1^{-,sm}(V).
\end{eqnarray}
By taking the orthogonal complement, we have a natural
identification $\PP(V)=\PP(V^*)$. Let $E\subset V$ be a line. Then
$E=(V^*/E^\perp)^*=(\det V^*\otimes ( E^\perp)^*)^*.$ Hence
\begin{eqnarray}\label{2-3}
E= E^\perp\otimes \det V.
\end{eqnarray}
Similarly
\begin{eqnarray}\label{2-4}
or E=or E^\perp\otimes or V.
\end{eqnarray}
Hence
\begin{eqnarray*}
C^\infty(\PP(V),\tv^*)=C^\infty(\PP(V^*),\det
\ct(V^*)^*\otimes\unv^*)=\\
C^\infty(\PP(V^*),\ct(V^*)^*)\otimes\det V^*.
\end{eqnarray*}
Hence (\ref{2-2}) can be rewritten as
\begin{eqnarray}\label{2-5}
C^\infty(\PP(V^*),\ct(V^*)^*)\otimes Dens(V)\surj Val_1^{-,sm}(V).
\end{eqnarray}
Replacing $V$ by $V^*$ in (\ref{2-2}) we get an epimorphism
\begin{eqnarray}\label{2-6}
C^\infty(\PP(V^*),\ct(V^*)^*)\otimes orV^*\surj Val_1^{-,sm}(V^*).
\end{eqnarray}
Tensoring (\ref{2-6}) by $\det V^*$ and observing that $or V=or
V^*$ we get
\begin{eqnarray}\label{2-7}
C^\infty(\PP(V^*),\ct(V^*)^*)\otimes Dens (V)\surj
Val_1^{-,sm}(V^*)\otimes \det V^*.
\end{eqnarray}
Comparing (\ref{2-5}) and (\ref{2-7}) we conclude that there
exists a unique isomorphism of $GL(V)$-modules
\begin{eqnarray}\label{2-8}
\tilde \FF_V\colon Val_1^{-,sm}(V)\tilde\to
Val_1^{-,sm}(V^*)\otimes \det V^*=Val_1^{-,sm}(V^*)\otimes
Dens(V)\otimes or V
\end{eqnarray}
which makes the following diagram commutative:
\def\aaa{C^\infty(\PP(V^*),\ct(V^*)^*)\otimes Dens
(V)}
\def\bbb{Val_1^{-,sm}(V)}
\def\ccc{Val_1^{-,sm}(V^*)\otimes \det V^*}
\begin{eqnarray*}
\Atriangle[\aaa`\bbb`\ccc;``\tilde\FF_V]
\end{eqnarray*}
Remind also that in the even case we have an isomorphism
\begin{eqnarray}\label{2-9}
\FF_V\colon Val^{+,sm}(V)\tilde\to Val^{+,sm}(V^*)\otimes Dens(V).
\end{eqnarray}
Combining (\ref{2-8}) and (\ref{2-9}) together we obtain a
canonical isomorphism
\begin{eqnarray}\label{2-10}
\bar\FF_V\colon Val^{sm}(V)\tilde \to (Val^{+,sm}(V^*)\oplus
(Val_1^{-,sm}(V^*)\otimes or V))\otimes Dens(V).
\end{eqnarray}

\subsection{Construction of the convolution product on $Val(V^*)\otimes Dens(V)$, $\dim
V=2$.}\label{construction} Let us denote
\begin{eqnarray*}
\ca_0:=Val_2(V^*)\otimes Dens(V),\\
\ca_2:=Dens(V),\\
\ca_1^+:=Val_1^{+,sm}(V^*)\otimes Dens(V),\\
\ca_1^-:=Val_1^{-,sm}(V^*)\otimes Dens(V)\otimes or V,\\
\ca:=(Val^{+,sm}(V^*)\oplus (Val_1^{-,sm}(V^*)\otimes or V))\otimes
Dens(V)=\ca_0\oplus \ca_2\oplus \ca_1^+\oplus \ca_1^-.
\end{eqnarray*}

Using the isomorphism $\bar \FF_V$ we can define the product $\star$
on $\ca$ by
$$\phi\star\psi:=\bar\FF_V(\bar\FF_V^{-1}\phi\cdot
\bar\FF_V^{-1}\psi).$$ Then $\ca$ becomes a commutative graded
algebra such that the graded components are
$\ca_0,\ca_1^+\oplus\ca_1^-,\ca_2$. Note that
\begin{eqnarray}\label{2-11}
\ca_1^-\star\ca_1^+=\ca_1^-\star\ca_2=0,\\\label{2-11.5}
\ca^+\star\ca_1^-\subset \ca_1^-.
\end{eqnarray}
Let us denote
$$\ca^+:=\ca_0\oplus\ca_1^+\oplus \ca_2.$$

Let us define now an algebra structure on
$\cb:=Val^{sm}(V^*)\otimes Dens(V)$. First write
\begin{eqnarray}\label{2-12}
\cb=\ca_0\oplus(\ca_1^+\oplus (\ca_1^-\otimes or (V)))\oplus
\ca_2.
\end{eqnarray}
Thus $\cb=\ca^+\oplus (\ca_1^-\otimes or (V))$. Let us define
$$\ast\colon \cb\otimes \cb\to \cb$$ as follows. First we have
(using the equality $or V\otimes or V=\CC$)
$$\cb\otimes \cb=(\ca^+\otimes \ca^+)\oplus(\ca^+\otimes
\ca_1^-\otimes or V)\oplus(\ca_1^-\otimes \ca^+\otimes or
V)\oplus(\ca_1^-\otimes \ca_1^-).$$ We define $\ast$ to be equal
to $\star$ on the first summand, to $\star\otimes id_{or V}$ on
the second and the third summands, and to $-\star$ on the fourth
summand (with the minus sign!).

It is easy to see that $(\cb,\ast)$ is a commutative associative
graded algebra with unit (since $(\ca,\star)$ is). The
decomposition of $\cb$ into the graded components is given by
(\ref{2-12}).

\subsection{Explicit computation of the convolution.}
In this subsection we show that in the two dimensional case the
convolution $\ast$ defined in Section \ref{construction} coincides
with the convolution introduced by Bernig and Fu \cite{bernig-fu}.

Let us fix a Lebesgue measure $vol$ on $V^*$. It gives an
isomorphism $Dens(V)\tilde \to \CC$. For $\ca\in \ck^{sm}(V^*)$
let us denote $\mu_A(K)=vol(K+A)\otimes vol^{-1}$ for any $K\in
\ck(V^*)$. It is easy to see that $\mu_A\in Val^{sm}(V^*)\otimes
Dens(V)$.
\begin{proposition}\label{3-1}
Let $V$ be a two dimensional real vector space. For the
convolution product $\ast$ on $Val^{sm}(V^*)\otimes Dens(V)$
defined in Section \ref{construction} one has
\begin{eqnarray}\label{3-2}
\mu_A\ast\mu_B=\mu_{A+B}
\end{eqnarray}
for any $A,B\in\ck^{sm}(V^*)$.
\end{proposition}
This proposition implies that our convolution $\ast$ coincides
with the Bernig-Fu convolution, and hence there is no abuse of
notation. Before we prove this proposition we will prove another
proposition.
\begin{proposition}\label{3-3}
Let us fix a Euclidean metric on $V$ and an orientation. Consider
the isomorphisms $V\tilde\to V^*,\,Dens(V)\tilde\to \CC,\, or
V\tilde \to\CC$ induced by these choices. With these identifications
consider $\bar\FF_V\colon Val^{sm}(V)\tilde\to Val^{sm}(V)$. Then
for any $A\in \ck^{sm}(V)$
\begin{eqnarray}\label{3-4}
\bar\FF_V(V(\bullet,A))=V(\bullet,J^{-1}A)
\end{eqnarray}
where $J\colon V\to V$ is the rotation by $\frac{\pi}{2}$
counterclockwise.
\end{proposition}
{\bf Proof.} First let us remind that for any smooth functions
$h_1,h_2$ on the unit circle one can define the mixed volume
$V(h_1,h_2)$ (see e.g. \cite{goodey-weil-84}) which is bilinear
with respect to $h_1,h_2$. The idea is as follows. First if
$h_1,h_2$ are supporting functionals of convex sets $A_1,A_2$
respectively, let us define $V(h_1,h_2)$ to be equal to
$V(A_1,A_2)$. Since every smooth function on the circle is a
difference of supporting functions of convex bodies, let us extend
this expression by bilinearity. We get a well defined notion. Thus
we may and will identify  convex set $A$ with its supporting
function. We may assume that it is smooth on $S^1=\PP_+(V)$.

\underline{Case 1.} Let us assume first that $A$ is an even
function. Then the value at $l\in \PP(V)$ of the Klain imbedding
of $V(\bullet,A)$ to $C^{+,\infty}(\PP(V))$ is equal to $\kappa
A(l^\perp)$ where $\kappa$ is a normalizing constant. By the
definition of $\bar\FF_V$ on even valuations, the value at $l\in
\PP(V)$ of the Klain imbedding of $\bar \FF_V(V(\bullet,A))$ is
equal to $\kappa A(l^\perp)=\kappa (J^{-1} A)(l)$. Thus the
proposition is proved in the even case.

\underline{Case 2.} Let us assume now that $A$ is an odd function.
Recall that for $K\in \ck(V)$
$$V(K,A)=\frac{1}{2} \int_{l\in \PP_+(V)}A(l)dS_1(K,l)$$
where $S_1(K,\cdot)$ is the first area measure of $K$ (see e.g.
\cite{schneider-book}, formula (5.1.18)). By the construction of
$\bar\FF_V$ on odd valuations
$$\bar\FF_V(V(\bullet,A))(K)=\frac{1}{2}
\int_{l\in\PP_+(V)}A(Jl)dS_1(K,l)=\frac{1}{2}
\int_{l\in\PP_+(V)}(J^{-1}A)(l)dS_1(K,l).$$ Thus the proposition
is proved. \qed

{\bf Proof of Proposition \ref{3-2}.} Let us fix a Euclidean
metric and an orientation on $V$. It is easy to check (see e.g.
\cite{bernig-fu}, Corollary 1.3) that the equality (\ref{3-2}) for
all $A,B\in \ck^{sm}(V^*)$ is equivalent to the following two
conditions:

(a) $vol$ is the unit element with respect to $\ast$, i.e.
$vol\ast x=x$ for any $x$;

(b) $V(\bullet,A)\ast V(\bullet,
B)=\frac{1}{2}V(A,B)\chi(\bullet)$ for any $A,B\in \ck^{sm}(V^*)$.

The property (a) holds due to the corresponding property of the
product $\star$ on $\ca$. Let us prove (b). We will prove it in a
more general form when $A$ and $B$ are smooth functions on the
unit circle. Since by (\ref{2-11}) $\ca_1^+\star\ca_1^-=0$, the
proof splits into two cases:

(a) both $A$ and $B$ are even;

(b) both $A$ and $B$ are odd.

Though the even case was considered in \cite{bernig-fu} (for all
dimensions), we will prove it here in two dimensions for the sake
of completeness. Thus let us first consider the case (a). Since
$A,B$ are even, we have
\begin{eqnarray*}
V(\bullet,A)\ast V(\bullet,B)=\bar\FF_V\left(\bar\FF_V^{-1}(V(\bullet,A))\cdot\bar\FF_V^{-1}(V(\bullet,B))\right)=\\
\bar\FF_V\left(V(\bullet,J A)\cdot V(\bullet,J
B)\right)\overset{\mbox{Example }
\ref{E:product-2dim}}{=}\bar\FF_V\left(\frac{1}{2}V(J A,-J B)\cdot
vol(\bullet)\right)=\\
\frac{1}{2}V(A,-B)\cdot\chi(\bullet)=\frac{1}{2}V(A,B)\cdot
\chi(\bullet).
\end{eqnarray*}
Now let us consider the case (b), i.e. $A,B$ are odd. Then
$V(\bullet,A),\, V(\bullet,B)\in \ca_1^-\otimes or V$. Hence we
have by the definition of $\ast$
$$V(\bullet,A)\ast
V(\bullet,B)=-\bar\FF_V\left(\bar\FF_V^{-1}(V(\bullet,A))\cdot\bar\FF_V^{-1}(V(\bullet,B))\right).$$
Similarly to the previous case the last expression is equal to
$$-\frac{1}{2}V(A,-B)\chi(\bullet)=\frac{1}{2}V(A,B)\chi(\bullet).$$
Proposition is proved. \qed


\def\valsm{Val^{sm}}
\subsection{Isomorphisms of algebras $\valsm(V)$ and $\valsm(V^*)\otimes
Dens(V)$.}\label{two-dim-isomorphism} We are going to prove the
following result.
\begin{theorem}\label{T:ISO}
Let $V$ be a two dimensional real vector space. There exists an
isomorphism of topological vector spaces
$$\FF_V\colon \valsm(V)\tilde \to\valsm(V^*)\otimes Dens(V)$$
which satisfies the following properties:

(1) $\FF_V$ commutes with the natural $GL(V)$-action;

(2) $\FF_V$ is an isomorphism of algebras, i.e. $\FF_V(\phi\cdot
\psi)=\FF_V(\phi)\ast \FF_V(\psi)$ for any $\phi,\psi\in
\valsm(V)$;

(3) This isomorphism takes real valued valuations to real valued.


\end{theorem}
{\bf Proof.}
The isomorphism on even valuations $\FF_V\colon
Val^{+,sm}(V)\tilde\to Val^{+,sm}(V^*)\otimes Dens(V)$ was defined
by the author in \cite{alesker-jdg-03}. Its homomorphism property
was proved by Bernig and Fu in \cite{bernig-fu}. We will extend
this map to odd valuations. The construction will depend on a
choice of orientation of $V$. Thus let us fix an orientation of
$V$.

Consider the action of the subgroup $SL(V)\subset GL(V)$ on
$Val^{-,sm}_1(V^*)\otimes Dens(V)$. Then
\begin{eqnarray*}
Val^{-,sm}_1(V^*)\otimes Dens(V)=W_1\oplus W_2
\end{eqnarray*}
where $W_1,W_2\subset Val^{-,sm}_1(V^*)\otimes Dens(V)$ are closed
$SL(V)$-invariant $SL(V)$-irreducible infinite dimensional
subspaces. $W_1$ and $W_2$ are not isomorphic to each other as
$SL(V)$-modules. Moreover one can choose a Cartan subalgebra and a
root system of the Lie algebra of $sl(V)$ (compatible in an
appropriate way with the orientation of $V$) so that $W_1$ is a
highest weight module, and $W_2$ is a lowest weight module. The
property of being either highest or lowest module depends only on
the orientation of $V$, and not on a Cartan subalgebra and a
positive root system.
\begin{lemma}\label{L:gl(2)}
The map $$\LL:=Id_{W_1}\oplus (-Id_{W_2})\colon
Val^{-,sm}_1(V^*)\otimes Dens(V)\to Val^{-,sm}_1(V^*)\otimes
Dens(V)\otimes or(V)$$ is an isomorphism of $GL(V)$-modules.
\end{lemma}
Here is a warning: the notation of the lemma is a bit misleading.
The target space seems to be different form the source space since
it is twisted by $or(V)$.  The meaning is that we just consider the
target space to be equal to the source space as a vector space, but
the action of $GL(V)$ on the target is different: is it twisted by
the sign of the determinant of a matrix.

Let us postpone the proof of this lemma and finish the proof of the
theorem. Let us define $\FF_V\colon Val^{-,sm}_1(V)\to
Val^{-,sm}_1(V^*)\otimes Dens(V)$ by
\begin{eqnarray*}
\FF_V:=\LL^{-1}\circ \bar \FF_V \mbox{ on } Val^{-,sm}_1(V)
\end{eqnarray*}
where $\bar \FF_V$ was defined in (\ref{2-10}).

In order to prove that $\FF_V$ is an isomorphism of algebras it
remains to show that for any $\phi,\psi\in Val^{-,sm}_1(V)$ one
has
$$\FF_V(\phi\cdot \psi)=\FF_V(\phi)\ast\FF_V(\psi).$$
In the notation of Section \ref{construction} one has
\begin{eqnarray*}
\FF_V(\phi\cdot \psi)=\bar\FF_V\phi\star\bar\FF_V\psi;\\
\FF_V(\phi)\ast\FF_V(\psi)=-\FF_V(\phi)\star\FF_V(\psi)=
-\LL^{-1}(\bar\FF_V\phi)\star\LL^{-1}(\bar\FF_V\psi).
\end{eqnarray*}
Thus we have to show that for any $u,v\in Val^{-,sm}_1(V^*)\otimes
Dens(V)$ one has
\begin{eqnarray}\label{starrr}
u\star v=-\LL^{-1}(u)\star \LL^{-1}(v).
\end{eqnarray}
It is clear that if $u\in W_1,v\in W_2$ then (\ref{starrr}) holds.
For $u,v\in W_i$, $i=1,2$, the equality (\ref{starrr}) is
equivalent to $u\star v=0$. Let us prove it for $i=1$; the case
$i=2$ is considered similarly. Let us observe that $\star$ induces
the $SL(V)$-equivariant map $W_1\to W_1^*\otimes Dens(V)$. Since
$W_1$ is a highest weight irreducible infinite dimensional
$SL(V)$-module, $W_1^*\otimes Dens(V)$ is a lowest weight
irreducible infinite dimensional $SL(V)$-module. Hence they cannot
be isomorphic. This proves that $\FF_V$ is a homomorphism of
algebras.

It it clear from the construction of $\FF_V$ that $\FF_V$ maps
real valued valuations to real valued. Hence theorem is proved
modulo Lemma \ref{L:gl(2)}. \qed

{\bf Proof of Lemma \ref{L:gl(2)}.} It follows from general
representation theory of the group $GL_2(\RR)$ (see e.g.
\cite{jacquet-langlands}, Ch. I, Theorem 5.11(VI)) that the
irreducible representations $Val^{-,sm}_1(V^*)\otimes Dens(V)$ and
$Val^{-,sm}_1(V^*)\otimes Dens(V)\otimes or(V)$ of the group
$GL(V)\simeq GL_2(\RR)$ are isomorphic. Let us fix such an
isomorphism $\LL'$. Since $W_1,W_2$ are irreducible non-isomorphic
$SL(V)$-modules, $\LL'$ must have the form $$\LL'=\alp
Id_{W_1}\oplus \beta Id_{W_2}$$ where $\alp,\beta\in \CC^*$.
Dividing by $\alp$ we may assume that $\alp=1$. The composition
$$(\LL'\otimes Id_{or(V)})\circ \LL'\colon Val^{-,sm}_1(V)\otimes
Dens(V)\to Val^{-,sm}_1(V)\otimes Dens(V)$$ is an automorphism of
$GL(V)$-module, hence it must have the form $\gamma
(Id_{W_1}\oplus Id_{W_2})$. On the other hand this composition is
equal to $Id_{W_1}\oplus \beta^2 Id_{W_2}$. Hence $\beta^2=1$,
i.e. $\beta=\pm 1$. Since the identity map $Id_{W_1\oplus W_2}$ is
{\itshape not} a morphism of $GL(V)$-modules
$Val^{-,sm}_1(V^*)\otimes Dens(V)$ and $Val^{-,sm}_1(V^*)\otimes
Dens(V)\otimes or(V)$, it implies that $\beta=-1$. Lemma is
proved. \qed
\subsection{Plancherel type formula in two dimensions.}\label{Ss:composition-F-2d}
Let us consider the map
\begin{eqnarray*}
\FF_{V^*}\otimes Id_{Dens(V)}\colon Val^{sm}(V^*)\otimes
Dens(V)\tilde\to Val^{sm}(V).
\end{eqnarray*}
Let us denote by $\ce_V\colon Val^{sm}(V)\to Val^{sm}(V)$ the
operator given by $(\ce_V\phi)(K)=\phi(-K)$ for any $\phi\in
Val^{sm}(V)$, $K\in\ck(V)$.
\begin{proposition}\label{P:composition-2-dim} Let $V$ be a two
dimensional vector space. Then
$$(\FF_{V^*}\otimes Id_{Dens(V)})\circ \FF_V=\ce_V.$$
\end{proposition}
{\bf Proof.} Obviously, on 0- and 2-homogeneous valuations both
operators are obviously equal to identity. For 1-homogeneous
valuations the result follows from Proposition \ref{3-3}. \qed

\section{Fourier transform on valuations in higher
dimensions.}\label{S:fourier-high-dim} The goal of this section is
to construct the Fourier transform in higher dimensions and to
prove its main properties.

\subsection{Construction of the Fourier transform.}\label{Ss:construction-fourier}
Let $0\leq k\leq n$. Let us consider the (infinite dimensional)
vector bundle
$$\ct^0_{k,V}\to Gr_{n-k}(V)$$ whose fiber over $F\in Gr_{n-k}(V)$ is
equal to $Val^{sm}(V/F)$. Similarly let $\ct_{k,V;i}^{0}\to
Gr_{n-k}(V),\,\ct_{k,V;i}^{0, \pm}\to Gr_{n-k}(V)$ be the vector
bundles whose fiber over $F\in Gr_{n-k}(V)$ is equal to
$Val^{sm}_i(V/F)$ or $Val^{\pm,sm}_i(V/F)$ respectively. Let
\begin{eqnarray}\label{w4}
\ct_{k,V}:=\ct^0_{k,V}\otimes |\ome_{Gr_{n-k}(V)}|,\\
\ct_{k,V;i}:=\ct_{k,V;i}^0\otimes
|\ome_{Gr_{n-k}(V)}|,\\\label{w5}
\ct^\pm_{k,V;i}:=\ct^{0,\pm}_{k,V;i}\otimes |\ome_{Gr_{n-k}(V)}|.
\end{eqnarray}
Also for a subspace $F\subset V$ let us denote by $p_F$ the
canonical map
$$p_F\colon V\to V/F.$$

Consider the natural map
$$\Xi_{k,V}\colon C^\infty (Gr_{n-k}(V),\ct_{k,V;i})\to
Val_i^{sm}(V)$$ defined by $\xi\mapsto \int_{F\in
Gr_{n-k}(V)}p^*_F(\xi(F)).$ Note that the map $\Xi_{k,V}$ is
$GL(V)$-equivariant, hence its image is indeed contained in
$Val_i^{sm}(V)$.
\begin{lemma}\label{L:pull-cosm}
Let $f_i\colon V_i\to W_i,\, i=1,2$ be two linear maps such that
$f_2$ is injective. Then for any $\phi_1\in Val(W_1),\, \phi_2\in
Val^{sm}(W_2)$ the pushforward $f_2^*(\phi_2)\in Val^{sm}(V_2)$
and
\begin{eqnarray}\label{E:pull-cosm}
(f_1\boxtimes f_2)^*(\phi_1\boxtimes\phi_2)=f_1^*\phi_1\boxtimes
f_2^*\phi_2.
\end{eqnarray}
\end{lemma}
\begin{remark}
In the statement of the proposition we use the construction of the
product of a continuous valuation and a smooth one from the
appendix.
\end{remark}
{\bf Proof} of Lemma \ref{L:pull-cosm}. It is clear that
$f_2^*(\phi_2)\in Val^{sm}(V_2)$. Next both sides of
(\ref{E:pull-cosm}) are continuous with respect to $\phi_1\in
Val(W_1),\,\phi_2\in Val^{sm}(W_2)$. Hence by the McMullen's
conjecture we may assume that
$$\phi_i(\bullet)=vol_i(\bullet +A_i),\, i=1,2,$$
where $vol_i$ is a Lebesgue measure on $W_i$, $A_i\in
\ck^{sm}(W_i)$. Then one has
\begin{eqnarray}\label{E:pc1}
(\phi_1\boxtimes \phi_2)(K)=(vol_1\boxtimes vol_2)(K+(A_1\times
A_2)),\\\label{E:pc2} (f_1\boxtimes f_2)^*(\phi_1\boxtimes
\phi_2)(K)=(vol_1\boxtimes vol_2)\left((f_1\boxtimes
f_2)(K)+(A_1\times A_2)\right).
\end{eqnarray}
On the other hand
\begin{eqnarray*}
(f_i^*\phi_i)(K)=vol_i(f_i(K)+A_i), i=1,2.
\end{eqnarray*}

Note that $(f_1\boxtimes f_2)^*=(f_1\boxtimes Id)^*\circ
(Id\boxtimes f_2)^*$. Let us compute first $(Id\boxtimes
f_2)^*(\phi_1\boxtimes \phi_2)$. First we will identify $V_2$ with
its image in $W_2$. Let us fix Lebesgue measures $vol_{W_2}$ on
$W_2$ and $vol_{W_2/V_2}$ on $V_{W_2/V_2}$ such that
$vol_2=vol_{V_2}\otimes vol_{W_2/V_2}$. We have
\begin{eqnarray}\label{E:pc3}
(f_2^*\phi_2)(\bullet)=\int_{z\in W_2/V_2} vol_{V_2}(\bullet
+(A_2\cap z))dvol_{W_2/V_2}(z).
\end{eqnarray}

Next for any $K\in \ck(V_1\times V_2)$ we have
\begin{eqnarray}
(Id\times f_2)^*(\phi_1\boxtimes \phi_2)(K)=\\\int_{z\in
W_2/V_2}(vol_1\boxtimes vol_2)\left((K+(A_1\times A_2))\cap
(V_1\times z)\right)dvol_{W_2/V_2}(z)=\\\label{E:pc4} \int_{z\in
W_2/V_2}(vol_1\boxtimes vol_2)\left(K+A_1\times (A_2\cap
z)\right)dvol_{W_2/V_2}(z).
\end{eqnarray}
Now observe that the map $W_2/V_2\to Val(W_2)$ defined by
$$z\mapsto vol_{V_2}(\bullet +(A_2\cap z))=:\phi(z)(\bullet)$$
is a bounded map which is continuous almost everywhere (more
precisely, this map is continuous in the interior of the image of
$A_2$ in $W_2/V_2$). The expression (\ref{E:pc4}) is equal to
\begin{eqnarray}\label{E:pc5}
\int_{z\in W_2/V_2}(\phi_1\boxtimes \phi(z))(K)dvol_{W_2/V_2}(z),
\end{eqnarray}
and by the mentioned above continuity we may rewrite the
expression (\ref{E:pc5}) as $$\left(\phi_1\boxtimes \int_{z\in
W_2/V_2}\phi(z)dvol_{W_2/v_2}(z)\right)(K).$$

By (\ref{E:pc3}) $\int_{z\in W_2/V_2}\phi(z)
dvol_{W_2/V_2}(z)=f_2^*\phi_2$. Thus we have proven that
$(Id\times f_2)^*(\phi_1\boxtimes \phi_2)=\phi_1\boxtimes
f_2^*\phi_2$.

Now we may assume that $f_2=Id$, and it remains to show that
$$(f_1\times Id)^*(\phi_1\boxtimes\phi_2)=f_1^*\phi_1\boxtimes
\phi_2.$$ Let us decompose $f_1=g\circ h$ where $h$ is surjection
and $g$ is injection. Then $(f_1\boxtimes Id)^*=(h\boxtimes
Id)^*\circ (g\boxtimes Id)^*$. Since we have assumed that both
$\phi_1$ and $\phi_2$ can be chosen smooth, by the above proven
case of injections one has $(g\boxtimes
Id)^*(\phi_1\boxtimes\phi_2)=g^*\phi_1\boxtimes \phi_2$.

Thus it remains to show that for any surjection $h\colon V_1\surj
W_1$, and any $A_i\in \ck^{sm}(W_i)$,
$\phi_i(\bullet)=vol_i(\bullet +A_i),\, i=1,2,$ one has
$$(h\boxtimes Id)^*(\phi_1\boxtimes\phi_2) =h^*\phi_1\boxtimes
\phi_2.$$

Let $M:=Ker (h)$, $m:=\dim M$. We identify $W_1$ with $V_1/M$. Let
us fix a Lebesgue measure $vol_M$ on $M$. Let
$vol_{V_1}:=vol_M\otimes vol_1$. Let us fix $\tilde A_1\in
\ck^{sm}(V_1)$ such that $h(\tilde A_1)=A_1$. Let us fix
$S\in\ck(M)$ with $vol_M(S)=1$. Then by Lemma \ref{L:aLemma}
$$(h^*\phi_1)(\bullet)=vol_1(h(\bullet)+A_1)=
\frac{1}{k!}\frac{d^k}{d\eps^k}\big|_0 vol_{V_1}(\bullet+\tilde
A_1+\eps S).$$ Hence for $K\in \ck(V_1\times W_2)$ we have
\begin{eqnarray}\label{E:pc10}(h^*\phi_1\boxtimes
\phi_2)(K)=\frac{1}{k!}\frac{d^k}{d\eps^k}\big|_0
(vol_{V_1}\boxtimes vol_2)\left(K+((\tilde A_1+\eps S)\times
A_2)\right).\end{eqnarray} Again by Lemma \ref{L:aLemma} the right
hand side of (\ref{E:pc10}) is equal to $$(vol_1\boxtimes
vol_2)\left((h\boxtimes Id)(K)+(A_1\times A_2)\right)=(h\boxtimes
Id)^*(\phi_1\boxtimes \phi_2)(K).$$ Hence $(h\boxtimes
Id)^*(\phi_1\boxtimes \phi_2)=h^*\phi_1\boxtimes \phi_2$. Hence
lemma is proved. \qed


\begin{proposition}\label{exter-pull}
Let $0\leq k,l\leq n$. Let
\begin{eqnarray*}
\phi=\Xi_{k,V}(\xi),\, \psi=\Xi_{l,W}(\eta)
\end{eqnarray*}
with $\xi\in C^\infty(Gr_{n-k}(V),\ct_{k,V}),\, \eta\in
C^\infty(Gr_{n-l}(V),\ct_{l,V})$.

Then
\begin{eqnarray}\label{w1}
\phi\boxtimes\psi=\int_{F\in Gr_{n-k}(V)}\int_{E\in
Gr_{n-l}(W)}(p_F\boxtimes p_E)^*(\xi(F)\boxtimes \eta(E))
\end{eqnarray}
where $p_F\boxtimes p_E\colon V\times W\to V/F\times W/E$.
\end{proposition}
\begin{remark}
In (\ref{w1}) $\xi(F)\boxtimes \eta(E)$ is considered as an
element of $Val(V/F\times V/E)\otimes
|\ome_{Gr_{n-k}(V)}|\big|_F\otimes |\ome_{Gr_{n-l}(V)}|\big|_E$.
\end{remark}
{\bf Proof} of Proposition \ref{exter-pull}. We have
\begin{eqnarray*}
\phi\boxtimes \psi=\int_{F\in
Gr_{n-k}(V)}p_F^*(\xi(F))\boxtimes\psi\overset{\mbox{Lemma
}\ref{L:pull-cosm}}{=}\\
\int_{F\in Gr_{n-k}(V)}(p_F\boxtimes Id_W)^*\left(\xi(F)\boxtimes
\psi\right)\overset{\mbox{Lemma } \ref{L:pull-cosm}}{=}\\
\int_{F\in Gr_{n-k}(V)}(p_F\boxtimes Id_W)^*\int_{E\in
Gr_{n-l}(W)}(Id_F\boxtimes p_E)^*\left(\xi(F)\boxtimes
\eta(E)\right)=\\
\int_{F\in Gr_{n-k}(V)}\int_{E\in Gr_{n-l}(W)}(p_F\boxtimes
p_E)^*\left(\xi(F)\boxtimes \eta(E)\right).
\end{eqnarray*}
Proposition \ref{exter-pull} is proved. \qed

\begin{proposition}\label{P:w5}
For any $0\leq k\leq n$
\begin{eqnarray*}
\Xi_{k,V}\left(C^\infty(Gr_{n-k}(V),\ct^+_{k,V;k})\right)=Val^{+,sm}_k(V),\\
\Xi_{k+1,V}\left(C^\infty(Gr_{n-k-1}(V),\ct^-_{k+1,V;k})\right)=Val^{-,sm}_k(V).
\end{eqnarray*}
\end{proposition}
{\bf Proof.} The map $\Xi_{k,V}$ is $GL(V)$-equivariant non-zero
map. Hence the result follows from the Irreducibility theorem and
the Casselman-Wallach theorem. \qed

Let us remind the construction of an isomorphism
\begin{eqnarray*}
\FF_V\colon Val_k^{+,sm}(V)\tilde\to Val_{n-k}^{+,sm}(V^*)\otimes
Dens(V)
\end{eqnarray*}
on even valuations from \cite{alesker-jdg-03} (where it was
denoted by $\DD$). First let us do it for $k=n$. Then
$$Val^+_n(V)=Val_n(V)=Dens(V).$$
On the other hand
$$Val_0^{+,sm}(V^*)\otimes Dens(V)=Val_0(V^*)\otimes
Dens(V)=\CC\otimes Dens(V)=Dens(V).$$ Take
$$\FF_V\colon Val_n(V)\tilde\to Val_0(V)\otimes Dens(V)$$
to be the identity isomorphism.

Let us consider the case $k<n$. Let $F\in Gr_{n-k}(V)$. We have
\begin{eqnarray}\label{w6}
\FF_{V/F}\colon Val_k^{sm}(V/F)\tilde\to
Val_0^{sm}(F^\perp)\otimes Dens(F^{\perp *}).
\end{eqnarray}
Let $p_F^\vee\colon F^\perp\inj V^*$ be the map dual to the
projection $p_F\colon V\to V/F$. Let us consider the map
\begin{eqnarray}\label{w7}
S^+_k\colon C^\infty(Gr_{n-k}(V),\ct^+_{k,V;k})\to
Val_{n-k}^{+,sm}(V^*)\otimes Dens(V)
\end{eqnarray}
given by
\begin{eqnarray*}
S^+_k(\xi)=\int_{F\in Gr_{n-k}(V)}p^\vee_{F*} (\FF_{V/F}(\xi(F))).
\end{eqnarray*}
\begin{theorem}[\cite{alesker-jdg-03}]\label{T:w8}
There exists a unique map
$$\FF_V\colon Val_k^{+,sm}(V)\to Val_{n-k}^{+,sm}(V^*)\otimes
Dens(V)$$ which makes the following diagram commutative
\def\aa1{C^\infty(Gr_{n-k}(V),\ct^+_{k,V;k})}
\def\bb1{Val^{+,sm}_k(V)}
\def\cc1{Val_{n-k}^{+,sm}(V^*)\otimes Dens(V)}
\def\ff1{\Xi_{k,V}}
\def\gg1{S^+_k}
\def\hh1{\FF_V}
$$\Vtriangle[\aa1`\bb1`\cc1;\ff1`\gg1`\hh1].$$
This map $\FF_V$ is a $GL(V)$-equivariant isomorphism of
topological vector spaces.
\end{theorem}
Let us construct $\FF_V$ on {\itshape odd} $k$-homogeneous
valuations. Thus let $1\leq k\leq n-1$.

First let us consider the case $k=1$. We have
$$\Xi_{2,V}\colon C^\infty(Gr_{n-2}(V),\ct^-_{2,V;1})\surj
Val_1^{-,sm}(V).$$ Consider the map
\begin{eqnarray*}
S^-_2\colon C^\infty(Gr_{n-2}(V),\ct^-_{2,V;1})\to
Val^{-,sm}_{n-1}(V)
\end{eqnarray*}
given by
$$S^-_2(\xi)=\int_{F\in
Gr_{n-2}(V)}p^\vee_{F*}(\FF_{V/F}(\xi(F)))$$ where
$\FF_{V/F}\colon Val_1^{-,sm}(V/F)\tilde\to
Val_1^{-,sm}(F^\perp)\otimes Dens(F^{\perp *})$ is the Fourier
transform defined in the two dimensional case by Theorem
\ref{T:ISO}.

\begin{proposition}\label{P:w9}
There exists a unique map
$$\FF_V\colon Val_1^{-,sm}(V)\to Val_{n-1}^{-,sm}(V^*)\otimes
Dens(V)$$ making the following diagram commutative
$$\Vtriangle[C^\infty (Gr_{n-2}(V),\ct^{-}_{2,V;1})`Val^{-,sm}_1(V)`
Val^{-,sm}_{n-1}(V^*)\otimes Dens(V);\Xi_{2,V}`S_2^-`\FF_V].$$
This map $\FF_V$ is a $GL(V)$-equivariant isomorphism of
topological vector spaces.
\end{proposition}
Before we prove this proposition let us prove the following
result.
\begin{proposition}\label{P:w10}
Let $1\leq k\leq n-1$. The $GL(V)$-module
$C^\infty(Gr_{n-k},\ct^-_{k+1,V;k})$ is admissible of finite
length, and the Jordan-H\"older series of it contains
$Val_k^{-,sm}(V)$ with multiplicity one.
\end{proposition}
{\bf Proof.} Let us fix an isomorphism $V\simeq \RR^n$. Let $Q_0$
be the subgroup of $GL_n(\RR)$ consisting of matrices
\begin{eqnarray*}
Q_0=\left\{\left[\begin{array}{c|c|c}
                 A&*&*\\\hline
                 0&B&*\\\hline
                 0&0&c
                \end{array}\right]\big| A\in GL_{n-k-1}(\RR),B\in
                GL_k(\RR),c\in \RR^*.\right\}.
\end{eqnarray*}
$Q_0$ is a parabolic subgroup of $GL_n(\RR)$. Let us denote by
$\cv$ the complex line bundle over $\cf_{n-k-1,n-1}(\RR^n)\simeq
GL_n(\RR)/Q_0$ whose fiber over $(F,E)\in\cf_{n-k-1,n-1}(V)$ is
equal to $Dens(E/F)\otimes or(V/E)$. $\cv$ is
$GL_n(\RR)$-equivariant in a natural way.

Let $p\colon \cf_{n-k-1,n-1}(\RR^n)\to Gr_{n-k-1}(\RR^n)$ be the
natural map, i.e. $p(F,E)=F$. For any $F\in Gr_{n-k-1}(\RR^n)$ the
space $Val_k^{-,sm}(\RR^n/F)$ is canonically a quotient of
$C^\infty(p^{-1}(F),\cv)$. Hence the $GL_n(\RR)$-module
$C^\infty(Gr_{n-k-1}(V),\ct^-_{k+1,V;k})$ is a quotient of the
$GL_n(\RR)$-module $C^\infty(GL_n(\RR)/Q_0,\cv\otimes
p^*(|\ome_{Gr_{n-k}(V)}|))$.

Let now $P_0\subset GL_n(\RR)$ be the parabolic subgroup of
matrices
$$P_0=\left\{\left[\begin{array}{cc|c}
                 B&*&*\\
                 0&c&*\\\hline
                 0&0&A
                \end{array}\right]\big| A\in GL_{n-k-1}(\RR),B\in
                GL_k(\RR),c\in \RR^*.\right\}.$$

Recall that $Val^{-,sm}_k(\RR^n)$ is a subquotient of
$Ind^{GL_n(\RR)}_{P_0}\xi$ where $\xi\colon P_0\to \CC$ is the
character given by
$$\xi\left(\left[\begin{array}{cc|c}
                 B&*&*\\
                 0&c&*\\\hline
                 0&0&A
                \end{array}\right]\right)=sgn (c)\cdot |\det B|^{-1}.$$
Moreover the Jordan-H\"older series of $Ind^{GL_n(\RR)}_{P_0}\xi$
contains $Val^{-,sm}_k(\RR^n)$ with multiplicity one by Remark
\ref{R:mult1thm}. The $GL(V)$-modules
$C^\infty(GL_n(\RR)/Q_0,\cv\otimes p^*(|\ome_{Gr_{n-k}(V)}|))$ and
$Ind^{GL_n(\RR)}_{P_0}\xi$ have the same Jordan-H\"older series by
Corollary \ref{C:arepr}. Hence Proposition \ref{P:w10} is proved.
\qed

{\bf Proof} of Proposition \ref{P:w9}. It is clear that the map
$$S^-_2\colon C^\infty(Gr_{n-2},\ct^-_{2,V;1})\to
Val^{-,sm}_{n-1}(V^*)\otimes Dens(V)$$ is non-vanishing and
$GL(V)$-equivariant. By Proposition \ref{isomorphism} the
$GL(V)$-modules $Val^{-,sm}_1(V)$ and
$Val^{-,sm}_{n-1}(V^*)\otimes Dens(V)$ are isomorphic. Hence by
Proposition \ref{P:w10} all the irreducible subquotients of
$C^\infty(Gr_{n-2},\ct^-_{2,V;1})$ which are non-isomorphic to
$Val_1^{-,sm}(V)$, are contained in $Ker(S^-_2)$. Hence $\FF_V$
does exist indeed. Moreover $\FF_V$ is continuous by the open
mapping theorem for Fr\'echet spaces (see e.g. \cite{schaefer},
Ch. III, \S 2). The rest of the statements of the theorem follow
from the Irreducibility theorem and the Casselman-Wallach theorem.
\qed

Now let us consider the case $k=n-1$. By Proposition \ref{P:w9} we
have the isomorphism
$$\FF_{V^*}\colon Val_1^{-,sm}(V^*)\tilde\to
Val_{n-1}^{-,sm}(V)\otimes Dens(V^*).$$ Define
$$\FF_V\colon Val_{n-1}^{-,sm}(V)\to Val_1^{-,sm}(V^*)\otimes
Dens(V)$$ by
\begin{eqnarray}\label{D:w11}
\FF_V:=\ce_{V^*}\circ (\FF_{V^*}^{-1}\otimes Id_{Dens(V)})
\end{eqnarray}
where $(\ce_{V^*}\phi)(K)=\phi(-K)$ (thus $\ce_{V^*}$ is just
multiplication by -1 on $Val^{-,sm}(V^*)\otimes Dens(V)$).
\begin{remark}
For $n=2,k=1$ both definitions of $\FF_V\colon
Val_1^{-,sm}(V)\tilde\to Val_1^{-,sm}(V^*)\otimes Dens(V)$
coincide with the construction discussed in Section
\ref{two-dim-isomorphism} due to Proposition
\ref{P:composition-2-dim}.
\end{remark}
Now let us consider the general case $1\leq k\leq n-1$. We have
the epimorphism
$$\Xi_{k+1,V}\colon C^\infty(Gr_{n-k-1}(V),\ct^-_{k+1,V;k})\to
Val_k^{-,sm}(V).$$ Consider the map
$$S^-_{k+1}\colon C^\infty(Gr_{n-k-1}(V),\ct^-_{k+1,V;k})\to
Val_{n-k}^{-,sm}(V^*)\otimes Dens(V)$$ defined by
$S^-_{k+1}(\xi)=\int_{F\in
Gr_{n-k-1}(V)}p^\vee_{F*}(\FF_{V/F}(\xi))$ where $\FF_{V/F}$ in the
right hand side is defined by (\ref{D:w11}).
\begin{theorem}\label{T:w12}
Let $1\leq k\leq n-1$. There exists a unique map $\FF_V\colon
Val_k^{-,sm}(V)\to Val_{n-k}^{-,sm}(V^*)\otimes Dens(V)$ making
the following diagram commutative
$$\Vtriangle[C^\infty(Gr_{n-k-1}(V),\ct^-_{k+1,V;k})`Val_k^{-,sm}(V)`
Val_{n-k}^{-,sm}(V^*)\otimes
Dens(V);\Xi_{k+1,V}`S^-_{k+1}`\FF_V].$$ Thus map $\FF_V$ is a
$GL(V)$-equivariant isomorphism of linear topological spaces.
\end{theorem}
\begin{remark}
It is clear that the construction of $\FF_V$ given in Theorem
\ref{T:w12} coincides with the previous constructions for
$k=1,n-1$.
\end{remark}
{\bf Proof} of Theorem \ref{T:w12}.  It is clear that the map
$$S^-_{k+1}\colon C^\infty(Gr_{n-k-1},\ct^-_{k+1,V,k})\to
Val^{-,sm}_k(V^*)\otimes Dens(V)$$ is non-vanishing and
$GL(V)$-equivariant. By Proposition \ref{isomorphism} the
$GL(V)$-modules $Val^{-,sm}_k(V)$ and
$Val^{-,sm}_{n-k}(V^*)\otimes Dens(V)$ are isomorphic. Hence by
Proposition \ref{P:w10} all the irreducible subquotients of
$C^\infty(Gr_{n-k-1},\ct^-_{k+1,V;k})$ which are non-isomorphic to
$Val_k^{-,sm}(V)$, are contained in $Ker(S^-_{k+1})$. Hence
$\FF_V$ does exist indeed. Moreover $\FF_V$ is continuous by the
open mapping theorem for Fr\'echet spaces (see e.g.
\cite{schaefer}, Ch. III, \S 2). The rest of the theorem follow
from the Irreducibility theorem and the Casselman-Wallach theorem.
\qed

\subsection{Relations of the Fourier transform to the pullback
and pushforward.}

\begin{theorem}\label{T:z3}
Let $i\colon L\inj V$ be an injection of linear spaces. Let
$\phi\in Val^{sm}(V)$. Then $i^*\phi\in Val^{sm}(L)$, and
\begin{eqnarray}\label{z4}
\FF_L(i^*\phi)=i^\vee_*(\FF_V\phi).
\end{eqnarray}
\end{theorem}
{\bf Proof.} It is clear that $i^*\phi\in Val^{sm}(L)$. Let us
prove the equality (\ref{z4}). Obviously this equality is true if
$i$ is an isomorphism.

\underline{Case 1.} Assume that $\phi\in Val^{+,sm}_k(V)$. If
$\dim L<k$ then both sides of (\ref{z4}) vanish.

\underline{Case 1a.} Let us assume in addition that $\dim L=k$. By
Proposition \ref{P:w5} there exists $\xi\in
C^\infty(Gr_{n-k}(V),\ct^+_{k,V;k})$ such that
$$\phi=\Xi_{k,V}(\xi).$$
More explicitly
$$\phi=\int_{F\in Gr_{n-k}(V)}p^*_{F}(\xi(F)).$$
Then
\begin{eqnarray*}
i^*\phi=\int_{F\in Gr_{n-k}(V)}(p_F\circ i)^*(\xi(F))=\int_{F\in
Gr_{n-k}(V),F\cap L=\{0\}}(p_F\circ i)^*(\xi(F)).
\end{eqnarray*}
Now observe that for any $F\in Gr_{n-k}(V)$ such that $F\cap
L=\{0\}$ the map $p_F\circ i\colon L\to V/F$ is an isomorphism.
Hence
\begin{eqnarray*}
\FF_L(i^*\phi)=\int_{F\in Gr_{n-k}(V),F\cap
L=\{0\}}\FF_L((p_F\circ i)^*(\xi(F)))=\\\int_{F\in
Gr_{n-k}(V),F\cap
L=\{0\}}(i^\vee_*\circ p^\vee_{F*})(\FF_{V/F}(\xi(F)))=\\
i^\vee_*\left(\int_{F\in Gr_{n-k}(V),F\cap
L=\{0\}}p_{F*}^\vee(\FF_{V/F}(\xi(F)))\right)=\\i^\vee_*\left(\int_{F\in
Gr_{n-k}(V)}
p_{F*}^\vee(\FF_{V/F}(\xi(F)))\right)=i^\vee_*(\FF_V\phi)
\end{eqnarray*}
where the last equality holds by the definition of $\FF_V$. Case
1a is proved.

\underline{Case 1b.} Let us assume that $\dim L>k$. Let $E\subset
L$ be an arbitrary $k$-dimensional subspace of $L$, and let
$j\colon E\inj L$ denote the imbedding map. Using Case 1a we have
\begin{eqnarray*}
\FF_E(j^*(i^*\phi))=\FF_E((ij)^*\phi)=
(ij)^\vee_*(\FF_V\phi)=j^\vee_*(i^\vee_*(\FF_V\phi)).
\end{eqnarray*}
On the other hand using again Case 1a, we have
\begin{eqnarray*}
\FF_E(j^*(i^*\phi))=j^\vee_*(\FF_L(i^*\phi)).
\end{eqnarray*}
Thus we get
$$j^\vee_*(i^\vee_*(\FF_V\phi))=j^\vee_*(\FF_L(i^*\phi)).$$
Thus Theorem \ref{T:z3} in the even case follows from the
following lemma.
\begin{lemma}\label{L:z5}
Let $L$ be an $l$-dimensional vector space. Let $0\leq k<l$. Let
$\phi\in Val^+_{l-k}(L)\otimes Dens(L^*)$. Assume that for any
surjection $p\colon L\surj K$ of rank $k$
$$p_*\psi=0.$$
Then $\psi=0$.
\end{lemma}
{\bf Proof.} Let $\cv$ be the subspace of $Val^+_{l-k}(L)$
consisting of $\psi$ such that $p_*\psi=0$ for any surjection $p$
of rank $k$. It is easy to see that $\cv$ is $GL(V)$-invariant
closed subspace, $\cv\ne Val^+_{l-k}(L)$. Hence by the
Irreducibility theorem $\cv=0$. \qed

Theorem \ref{T:z3} is proved in the even case.

\underline{Case 2.} Assume that $\phi\in Val_k^{-,sm}(V)$. The
proof of this case will be similar to the proof of Case 1. First
observe that if $\dim L\leq k$ then both sides of (\ref{z4})
vanish.

\underline{Case 2a.} Let us assume in addition that $\dim L=k+1$.
By Proposition \ref{P:w5} there exists $\eta\in
C^\infty(Gr_{n-k-1}(V),\ct^-_{k+1,V;k})$ such that
$$\phi=\Xi_{k+1,V}(\eta).$$
More explicitly
$$\phi=\int_{F\in Gr_{n-k-1(V)}}p_F^*(\eta(F)).$$
Analogously to Case 1a we have
\begin{eqnarray*}
i^*\phi=\int_{F\in Gr_{n-k-1}(V),F\cap L=\{0\}}(p_F\circ
i)^*(\eta(F)).
\end{eqnarray*}
For any $F\in Gr_{n-k-1}(V)$ such that $F\cap L=\{0\}$ the map
$p_F\circ i\colon L\to V/F$ is an isomorphism. Hence
\begin{eqnarray*}
\FF_L(i^*\phi)=\int_{F\in Gr_{n-k-1}(V),F\cap
L=\{0\}}\FF_L\left((p_F\circ i)^*(\eta(F))\right)=\\
\int_{F\in Gr_{n-k-1}(V),F\cap L=\{0\}}(p_F\circ
i)^\vee_*\left(\FF_{V/F}(\eta(F))\right)=\\i^\vee_*\left(\int_{F\in
Gr_{n-k-1}(V)}p^\vee_{F*}(\FF_{V/F}(\eta(F))\right)=i^\vee_*(\FF_V\phi)
\end{eqnarray*}
where the last equality holds by the definition of $\FF_V$.

\underline{Case 2b.} Let us assume now that $\dim L>k+1$. Let
$E\subset L$ be an arbitrary subspace of dimension $k+1$, and let
$j\colon E\inj L$ be the imbedding map. By Case 2a we have
\begin{eqnarray*}
\FF_E(j^*(i^*\phi))=\FF_E((ij)^*\phi)=\\
(ij)^\vee_*(\FF_V\phi)=j^\vee_*(i^\vee_*(\FF_V\phi)).
\end{eqnarray*}
Using again Case 2a we have
\begin{eqnarray*}
\FF_E(j^*(i^*\phi))=j^\vee_*(\FF_L(i^*\phi)).
\end{eqnarray*}
Thus we obtain
$$j^\vee_*(i^\vee_*(\FF_V\phi))=j^\vee_*(\FF_L(i^*\phi)).$$
Then the proof of Case 2 follows immediately from the following
lemma.
\begin{lemma}\label{L:z6}
Let $L$ be an $l$-dimensional vector space. Let $1\leq k\leq l-1$.
Let $\psi\in Val^-_{l-k}(L)\otimes Dens(L^*)$. Assume that for any
surjection $p\colon L\surj N$ of rank $k+1$
$$p_*\psi=0.$$ Then $\psi=0$.
\end{lemma}
{\bf Proof.} The proof is completely analogous to the proof of
Lemma \ref{L:z5} and will be omitted. \qed

Thus Theorem \ref{T:z3} is proved. \qed

\begin{theorem}\label{T:z1}
Let $\xi\in C^\infty(Gr_{n-k}(V),\ct_{k,V})$. Then
\begin{eqnarray}\label{z2}
\FF_V(\Xi_{k,V}(\xi))=\int_{F\in
Gr_{n-k}(V)}p^\vee_{F*}(\FF_{V/F}(\xi(F))).
\end{eqnarray}
\end{theorem}
{\bf Proof.} We may and will assume that $\xi\in
C^\infty(Gr_{n-k}(V),\ct^\pm_{k,V;i})$.

\underline{Case 1.} Let us assume that $\xi\in
C^\infty(Gr_{n-k}(V),\ct^+_{k,V;i})$. If $k<i$ then both sides of
(\ref{z2}) vanish. If $k=i$ then (\ref{z2}) is just the definition
of $\FF_V$. Thus let us assume that $k>i$.
\def\cfki{\cf_{n-k,n-i}(V)}
Let us denote by $\cf_{n-k,n-i}(V)$ the manifold of partial flags
$$\cf_{n-k,n-i}(V):=\{(F,E)|\, F\subset E,\dim F=n-k,\dim
E=n-i\}.$$ Let
\begin{eqnarray*}
p_{n-k}\colon \cfki\to Gr_{n-k}(V),\\
p_{n-i}\colon \cfki\to Gr_{n-i}(V)
\end{eqnarray*}
be the natural projections defined by $p_{n-k}((F,E))=F,\,
p_{n-i}((F,E))=E$. Recall that we denote by $\ct^{0,+}_{k,V;i}$
the vector bundle over $Gr_{n-k}(V)$ whose fiber over $F\in
Gr_{n-k}(V)$ is equal to $Val_i^{+,sm}(V/F)$, and by
$$\ct^+_{k,V;i}:=\ct^{0,+}_{k,V,i}\otimes |\ome_{Gr_{n-k}(V)}|.$$
Note that the fiber of $\ct^{0,+}_{i,V;i}(V)$ over $L\in
Gr_{n-i}(V)$ is equal to $Dens(V/L)$. Let us denote
\def\tlciz{\tilde\ct^{0,+}_{i,V;i}}
\def\tlci{\tilde\ct^0_{i,V;i}}
\begin{eqnarray*}
\tilde\ct^{0,+}_{i,V;i}:=p^*_{n-i}(\ct^{0,+}_{i,V;i}),\\
\tilde\ct^0_{i,V;i}:=\tlciz\otimes |\ome_{\cf_{n-k,n-i}(V)}|.
\end{eqnarray*}
We have the canonical map
\begin{eqnarray}\label{z7}
\Psi^+_{k,i}\colon
C^\infty(\cf_{n-k,n-i}(V),\tilde\ct^+_{i,V;i})\to
C^\infty(Gr_{n-k}(V),\ct^+_{k,V;i})
\end{eqnarray}
given by
$$(\Psi^+_{k,i}(\zeta))(F)=\int_{L\in Gr_{n-i}(V),L\supset F}
q^*_{L,F}(\zeta(L,F))$$ where $q_{L,F}\colon V/F\to V/L$ is the
canonical map. It is clear that $\Psi^+_{k,i}$ is
$GL(V)$-equivariant continuous map.
\begin{claim}\label{C:z8}
$\Psi^+_{k,i}$ is onto.
\end{claim}
{\bf Proof.} This claim easily follows from Propositions
\ref{P:w5} and \ref{part1-epi}. \qed

We have
\begin{eqnarray*}
\Xi_{k,V}(\Psi^+_{k,i}(\zeta))=\int_{F\in
Gr_{n-k}(V)}p_F^*\left(\int_{L\in Gr_{n-i}(V),L\supset
F}q^*_{L,F}(\zeta(L,F))\right)=\\
\int_{F\in Gr_{n-k}(V)}\left(\int_{L\in Gr_{n-i}(V),L\supset
F}p_L^*(\zeta(L,F))\right)=\int_{L\in
Gr_{n-i}(V)}p_L^*\left(\int_{F\in Gr_{n-k}(L)}\zeta(L,F)\right)=\\
\int_{L\in Gr_{n-i}(V)}p_L^*(\eta(L))
\end{eqnarray*}
where $\eta(L):=\int_{F\in Gr_{n-k}(L)}\zeta(L,F)$, thus $\eta\in
C^\infty(Gr_{n-i}(V),\ct^+_{i,V;i})$. Hence by the definition of
$\FF_V$ one has
\begin{eqnarray*}
\FF_V\left(\Xi_{k,V}(\Psi^+_{k,i}\zeta)\right)=\int_{L\in
Gr_{n-i}(V)}p^\vee_{L*}(\FF_{V/L}(\eta(L))).
\end{eqnarray*}
But by the continuity of the Fourier transform
$\FF_{V/L}(\eta(L))=\int_{F\in Gr_{n-i}(L)}\FF_{V/L}(\zeta(L,F)).$
Also
\begin{eqnarray*}
p_L=q_{L,F}\circ p_F,\\
p_L^\vee=p_F^\vee\circ q^\vee_{L,F}.
\end{eqnarray*}
Hence
\begin{eqnarray}\label{z9}
\FF_V(\Xi_{k,V}(\Psi_{k,i}(\zeta)))=
\int_{\cf_{n-k,n-i}(V)}p^\vee_{F*}\left(q^\vee_{L,F*}(\FF_{V/L}(\zeta(L)))\right).
\end{eqnarray}
By the definition of $\FF_{V/F}$ we have
\begin{eqnarray*}
\int_{L\in Gr_{n-i}(V),\, L\supset F}
q^\vee_{L,F*}\left(\FF_{V/L}(\zeta(L,F))\right)=
\FF_{V/F}\left(\int_{L\in Gr_{n-i}(V),\, L\supset F}
q^*_{L,F}(\zeta(L,F))\right)=\\\FF_{V/F}(\Psi^+_{k,i}(\zeta)(F)).
\end{eqnarray*}
Substituting this into (\ref{z9}) we get
\begin{eqnarray*}
\FF_V\left(\Xi_{k,V}(\Psi^+_{k,V}(\zeta))\right)= \int_{F\in
Gr_{n-k}(V)}p^\vee_{F*}\left(\FF_{V/F}(\Psi^+_{k,V}(\zeta)(F))\right).
\end{eqnarray*}
Since $\Psi^+_{k,i}$ is onto (by Claim \ref{C:z8}), Case 1 is
proved.

\underline{Case 2.} Now we assume that $\xi\in
C^\infty(Gr_{n-k}(V),\ct^-_{k,V;i})$.

If $k<i+1$ then both sides of (\ref{z2}) vanish. If $k=i+1$ then
the result is just the definition of $\FF_V$. Let us assume that
$k>i+1$. The proof will be analogous to Case 1. Similarly to Case
1 we denote by $\cf_{n-k,n-i-1}(V)$ the manifold of partial flags
$$\cf_{n-k,n-i-1}(V):=\{(F,E)|\,F\subset E\subset V, \dim F=n-k,\,
\dim E=n-i-1\}.$$ We denote
\def\fll{\cf_{n-k,n-i-1}(V)}
\begin{eqnarray*}
p_{n-k}\colon\fll\to Gr_{n-k}(V),\\
p_{n-i-1}\colon \fll\to Gr_{n-i-1}(V)
\end{eqnarray*}
the natural projections. Define
\begin{eqnarray*}
\tilde\ct^{0,-}_{i+1,V;i}:=p^*_{n-i-1}(\ct^{0,-}_{i+1,V;i}),\\
\tilde\ct^-_{i+1,V;i}:=\tilde\ct^{0,-}_{i+1,V;i}\otimes
|\ome_{\fll}|.
\end{eqnarray*}
We have the canonical map
\begin{eqnarray}\label{z10}
\Psi^-_{k,i}\colon C^\infty(\fll,\tilde\ct^-_{i+1,V;i})\to
C^\infty(Gr_{n-k},\ct^-_{k,V;i})
\end{eqnarray}
given, for any $F\in Gr_{n-k}(V)$, by
$$(\Psi^-_{k,i}(\zeta))(F)=\int_{L\in
Gr_{n-i-1}(V),L\supset F}q^*_{L,F}(\zeta(L,F))$$ where as
previously $q_{L,F}\colon V/F\to V/L$ is the canonical map. It is
clear that $\Psi^-_{k,i}$ is a $GL(V)$-equivariant continuous map.
\begin{claim}\label{C:z11}
$\Psi^-_{k,i}$ is onto.
\end{claim}
{\bf Proof} easily follows from Proposition \ref{P:w5}. \qed

The rest of the proof is parallel to Case 1. We have
\begin{eqnarray*}
\Xi_{k,V}(\Psi^-_{k,i}(\zeta))=\int_{F\in
Gr_{n-k}(V)}p^*_F\left(\int_{L\in Gr_{n-i-1}(V),L\supset
F}q^*_{L,F}(\zeta(L,F))\right)=\\
\int_{F\in Gr_{n-k}(V)}\int_{L\in Gr_{n-i-1}(V),L\supset
F}p_L^*(\zeta(L,F))=\int_{L\in Gr_{n-i-1}(V)}p_L^*\left(\int_{F\in
Gr_{n-k}(L)}\zeta(L,F)\right)=\\
\int_{L\in Gr_{n-i-1}(V)}p_L^*(\eta(L)).
\end{eqnarray*}
where $\eta(L):=\int_{F\in Gr_{n-k}(L)}\zeta(L,F)$, thus $\eta\in
C^\infty(Gr_{n-i-1}(V),\ct^-_{i+1,V;i})$. Hence by the definition
of $\FF_V$ we have
\begin{eqnarray*}
\FF_V\left(\Xi_{k,V}(\Psi^-_{k,i}(\zeta))\right)=\int_{L\in
Gr_{n-i-1}(V)}p^\vee_{L*}\left(\FF_{V/L}(\eta(L))\right).
\end{eqnarray*}
But by the continuity of the Fourier transform one has
$\FF_{V/L}(\eta(L))=\int_{F\in Gr_{n-k}(L)}\FF_{V/L}(\zeta(L,F))$.
Hence
\begin{eqnarray}\label{z12}
\FF_V\left(\Xi_{k,V}(\Psi^-_{k,i}(\zeta))\right)=
\int_{\fll}p^\vee_{F*}\left(q^\vee_{L,F*}(\FF_{V/L}(\zeta(L,F)))\right).
\end{eqnarray}
By the definition of $\FF_{V/F}$ we have
\begin{eqnarray*}
\int_{L\in Gr_{n-i-1}(V),L\supset F}
q^\vee_{L,F*}\left(\FF_{V/L}(\zeta(L))\right)=\FF_{V/F}\left(\int_{L\in
Gr_{n-i-1}(V),L\supset F}q^*_{L,F}(\zeta(L,F))\right)=\\
\FF_{V/F}\left(\Psi^-_{k,i}(\zeta)(F)\right).
\end{eqnarray*}
Substituting this into (\ref{z12}) we get
$$\FF_V\left(\Xi_{k,V}(\Psi^-_{k,i}(\zeta))\right)=\int_{F\in
Gr_{n-k(V)}}p^\vee_{F*}\left(\FF_{V/F}(\Psi^-_{k,i}(\zeta)(F))\right).$$
Since $\Psi^-_{k,i}$ is onto (by Claim \ref{C:z11}), Case 2 is
proved. Thus Theorem \ref{T:z1} is proved. \qed

\subsection{A Plancherel type formula in higher
dimensions}\label{Ss:composition-F-high} The main result of this
section is Theorem \ref{T:com02} below.

Let us introduce more notation. Let
\begin{eqnarray}\label{D:com1}
\cs_{k,V;i}^{0,\pm}\to Gr_k(V)
\end{eqnarray}
denote the vector bundle whose fiber over $E\in Gr_k(V)$ is equal
to $Val^{\pm,sm}_i(E)\otimes Dens(E^*)$. Let
\begin{eqnarray}\label{D:com2}
\cs_{k,V;i}^\pm:=\cs_{k,V;i}^{0,\pm}\otimes |\ome_{Gr_k(V)}|.
\end{eqnarray}
Similarly let
\begin{eqnarray}\label{D:com1.1}
\cs_{k,V}^{0}\to Gr_k(V)
\end{eqnarray}
denote the vector bundle whose fiber over $E\in Gr_k(V)$ is equal
to $Val^{sm}(E)\otimes Dens(E^*)$. Let
\begin{eqnarray}\label{D:com1.2}
\cs_{k,V}:=\cs_{k,V}^{0}\otimes |\ome_{Gr_k(V)}|.
\end{eqnarray}

For a vector space $W$ we can consider the map
\begin{eqnarray*}
\FF_W\otimes Id_{Dens(W^*)}\colon Val^{sm}(W)\otimes
Dens(W^*)\to\\
(Val^{sm}(W^*)\otimes Dens(W))\otimes Dens(W^*)=Val^{sm}(W^*).
\end{eqnarray*}
\begin{lemma}\label{L:com3}
 Let $V$ be an $n$-dimensional vector space. Let $\eta\in
 C^\infty(Gr_2(V),\cs^-_{2,V;1})$. Then
\begin{eqnarray}\label{com4}
(\FF_V\otimes Id_{Dens(V^*)})\left(\int_{E\in
Gr_2(V)}i_{E*}(\eta(E))\right)=\int_{Gr_2(V)}(i^\vee_E)^*(\FF_E\otimes
Id_{Dens(E^*)})(\eta(E))
\end{eqnarray}
where $i_E\colon E\inj V$ denotes the imbedding map.
\end{lemma}
{\bf Proof.} First note that $i_{E*}(\eta(E))\in
Val_{n-1}^-(V)\otimes Dens(V^*)$. Recall that by the definition of
the Fourier transform $\FF_V\colon Val^{-,sm}_{n-1}(V)\to
Val_1^{-,sm}(V^*)\otimes Dens(V)$
$$\FF_V=\ce_V\circ (\FF^{-1}_{V^*}\otimes Id_{Dens(V)}).$$
Hence the equality (\ref{com4}) is equivalent to
\begin{eqnarray}\label{com5}
\int_{E\in Gr_2(V)}i_{E*}(\eta(E))=(\ce_V\circ
\FF_{V^*})\left(\int_{E\in Gr_2(V)}(i^\vee_E)^*(\FF_E\otimes
Id_{Dens(E^*)})(\eta(E))\right).
\end{eqnarray}
The right hand side of (\ref{com5}) is equal to
$$\FF_{V^*}\left(\int_{E\in
Gr_2(V)}(i^\vee_E)^*\left(\ce_E\circ(\FF_E\otimes
Id_{Dens(E^*)})\right)(\eta(E))\right).$$ But by the Plancherel
formula in dimension two (Proposition \ref{P:composition-2-dim})
we know that
$$\ce_E\circ(\FF_E\otimes Id_{Dens(E^*)})=\FF_{E^*}^{-1}.$$
Hence (\ref{com5}) is equivalent to
\begin{eqnarray*}
\int_{E\in Gr_2(V)}i_{E*}(\eta(E))=\FF_{V^*}\left(\int_{E\in
Gr_2(V)}(i^\vee_E)^*\FF_E^{-1}(\eta(E))\right).
\end{eqnarray*}
The last equality is a special case of Theorem \ref{T:z1}. \qed.

\begin{lemma}\label{L:com01}
The composition map
$$Val_1^{-,sm}(V)\overset{\FF_V}{\to} Val_{n-1}^{-,sm}(V^*)\otimes
Dens(V)\overset{\FF_{V^*}\otimes Id_{Dens(V)}}{\to}
Val_1^{-,sm}(V)$$ is multiplication by -1.
\end{lemma}
{\bf Proof.} Let $\phi\in Val_1^{-,sm}(V)$. By Proposition
\ref{P:w5} there exists $\xi\in C^\infty(Gr_2(V),\ct^-_{2,V;1})$
such that $\phi=\Xi_{2,V}(\xi)$. Then by the definition of $\FF_V$
\begin{eqnarray*}
\FF_V(\phi)=\int_{F\in
Gr_{n-2}(V)}p^\vee_{F*}(\FF_{V/F}(\xi(F)))=\int_{E\in
Gr_2(V^*)}i_{E*}(\FF_{V/E^\perp}(\xi(E^\perp)).
\end{eqnarray*}
Then applying Lemma \ref{L:com3} we get
\begin{eqnarray}\label{com7}
(\FF_{V^*}\otimes Id_{Dens(V)})(\FF_V(\phi))=\int_{E\in
Gr_2(V^*)}(i^\vee_E)^*\left((\FF_E\otimes Id_{Dens(E^*)})\circ
\FF_{V/E^\perp}\right)(\xi(E^\perp))
\end{eqnarray}
But $(\FF_E\otimes Id_{Dens(E^*)})\circ
\FF_{V/E^\perp}=(\FF_E\otimes Id_{Dens(E^*)})\circ
\FF_{E^*}=\ce_{E^*}$ by the two-dimensional case (Proposition
\ref{P:composition-2-dim}). Substituting this into (\ref{com7}) we
get
\begin{eqnarray*}
(\FF_{V^*}\otimes Id_{Dens(V)})(\FF_V(\phi))=-\int_{E\in
Gr_2(V^*)}(i^\vee_E)^*(\xi(E^\perp))=-\int_{F\in
Gr_{n-2}(V)}p_F^*(\xi(F))=-\phi.
\end{eqnarray*}
Lemma is proved. \qed
\begin{lemma}\label{L:com01.1}
The composition map
$$Val^{+,sm}(V)\overset{\FF_V}{\to}Val^{+,sm}(V^*)\otimes
Dens(V)\overset{\FF_{V^*}\otimes Id_{Dens(V)}}{\to}
Val^{+,sm}(V)$$ is equal to the identity.
\end{lemma}
{\bf Proof.} Let $0\leq k\leq n$. Let $\cl_k(V)\to Gr_k(V)$ denote
the line bundle whose fiber over $F\in Gr_k(V)$ is equal to
$Dens(F)$. We have the Klain imbeddings
\begin{eqnarray*}
\tau\colon Val^{+,sm}_k(V)\inj C^\infty(Gr_k(V),\cl_k(V)),\\
\tau'\colon Val_{n-k}^{+,sm}\otimes Dens(V)\inj
C^\infty(Gr_{n-k}(V^*),\cl_{n-k}(V^*))\otimes Dens(V).
\end{eqnarray*}
Taking the orthogonal complement we have the identification
$$Gr_k(V)\simeq Gr_{n-k}(V^*).$$ Also for any $F\in Gr_k(V)$ we
have
$$Dens(F)=Dens(F^\perp)\otimes Dens(V).$$
Thus we have a $GL(V)$-equivariant isomorphism
$$\gamma\colon C^\infty(Gr_k(V),\cl_k(V))\tilde\to
C^\infty(Gr_{n-k}(V^*),\cl_{n-k}(V))\otimes Dens(V).$$ Let
$\phi\in Val^{+,sm}_k(V)$. By \cite{alesker-jdg-03}
$$\tau'(\FF_V\phi)=\gamma(\tau\phi).$$
It is easy to see that for any $\psi\in
Val_{n-k}^{+,sm}(V^*)\otimes Dens(V)$ one has
$$\tau\left((\FF_{V^*}\otimes
Id_{Dens(V)})(\psi)\right)=\gamma^{-1}(\tau'\psi).$$ Hence
\begin{eqnarray*}
\tau\left(((\FF_{V^*}\otimes Id_{dens(V)})\circ
\FF_V)(\phi)\right)=\gamma^{-1}(\tau'(\FF_V\phi))=
\gamma^{-1}(\gamma(\tau\phi))=\tau\phi.
\end{eqnarray*}
Lemma is proved. \qed


\begin{proposition}\label{T:fourier-id}
Let $p\colon V\surj W$ be a surjective linear map of vector
spaces. Then the following diagram is commutative
\begin{eqnarray}\label{fourier-id}
\square<1`1`1`1;1700`600>[Val^{sm}(V)\otimes Dens(V^*)`
Val^{sm}(V^*)`Val^{sm}(W)\otimes Dens(W^*)`
Val^{sm}(W^*);\FF_V\otimes
Id_{Dens(V^*)}`p_*`p^{\vee*}`\FF_W\otimes Id_{Dens(W^*)}]
\end{eqnarray}
\end{proposition}
{\bf Proof.} \underline{Case 1.} Let us consider the case of even
valuations.

By Lemma \ref{L:com01.1}, on even valuations $\FF_V\otimes
Id_{Dens(V^*)}=\FF_{V^*}^{-1},\,\FF_W\otimes
Id_{Dens(W^*)}=\FF_{W^*}^{-1}$. Hence the commutativity of the
diagram (\ref{fourier-id}) is equivalent to commutativity of the
following diagram
\begin{eqnarray*}
\square<1`1`1`1;1700`600>[Val^{+,sm}(V^*)`Val^{+,sm}(V)\otimes
Dens(V^*)` Val^{+,sm}(W^*)`Val^{+,sm}(W)\otimes
Dens(W^*);\FF_{V^*}`p^{\vee*}`p_*`\FF_{W^*}].
\end{eqnarray*}
But the last diagram is a special case of Theorem \ref{T:z3}. Thus
the theorem is proved in the even case.

\underline{Case 2.} Let us consider the odd case.

Let us denote $m:=\dim W$. Let $1\leq j\leq n-1$. We have to show
that the following diagram is commutative:
\begin{eqnarray}\label{fourier-id2}
\square<1`1`1`1;1800`600>[Val^{-,sm}_j(V)\otimes Dens(V^*)`
Val^{-,sm}_{n-j}(V^*)`Val^{-,sm}_{j-n+m}(W)\otimes Dens(W^*)`
Val^{-,sm}_{n-j}(W^*);\FF_V\otimes
Id_{Dens(V^*)}`p_*`p^{\vee*}`\FF_W\otimes Id_{Dens(W^*)}]
\end{eqnarray}

\underline{Case 2a.} Let us assume in addition that $\dim W=2$.

Then necessarily $j=n-1$. In this case, by Lemma \ref{L:com01}, we
have $\FF_V\otimes Id_{Dens(V^*)}=-\FF_{V^*}^{-1},\, \FF_W\otimes
Id_{Dens(W^*)}=-\FF_{W^*}^{-1}$. Hence the commutativity of the
diagram (\ref{fourier-id2}) is equivalent to the commutativity of
the diagram
\begin{eqnarray}\label{fourier-id3}
\square<1`1`1`1;1600`600>[Val^{-,sm}_{1}(V^*)`Val^{-,sm}_{n-1}(V)\otimes
Dens(V^*)`Val^{-,sm}_{1}(W^*)`Val^{-,sm}_{1}(W)\otimes
Dens(W^*);\FF_{V^*}`p^{\vee*}`p_*`\FF_{W^*}].
\end{eqnarray}
The last diagram is again a special case of Theorem \ref{T:z3}.

\underline{Case 2b.} Let us assume that $j-n+m=1$.

In other words we have to show that the following diagram is
commutative
\begin{eqnarray}\label{diag-case2b}
\square<1`1`1`1;1800`600>[Val^{-,sm}_j(V)\otimes Dens(V^*)`
Val^{-,sm}_{n-j}(V^*)`Val^{-,sm}_{1}(W)\otimes Dens(W^*)`
Val^{-,sm}_{m-1}(W^*);\FF_V\otimes
Id_{Dens(V^*)}`p_*`p^{\vee*}`\FF_W\otimes Id_{Dens(W^*)}].
\end{eqnarray}

Let us fix an arbitrary 2-dimensional subspace $i\colon E\inj W$.
Let us form a Cartesian square
$$\square<1`3`3`1;500`400>[Z`V`E`W;\tilde i`\tilde p`p`i].$$ Note that $p,\tilde p$
are surjections, $i,\tilde i$ are injections. Set $z:=\dim
Z=2+n-m$.

Let us consider the following cube diagram where we denote
$Dens(\bullet)$ by $D(\bullet)$ for brevity:


$$
\bfig \putsquare<1`1`1`1;2250`1500>(0,1800)[Val^{-,sm}_j(V)\otimes
D(V^*)`Val^{-,sm}_{n-j}(V^*)`Val_j^{-,sm}(Z)\otimes
D(V^*)`Val^{-,sm}_{z-j}(Z^*)\otimes D(Z)\otimes
D(V^*);\FF_V\otimes Id_{D(V^*)}`\tilde i^*\otimes
Id_{D(V^*)}`\tilde i^{\vee}_*\otimes Id_{D(V^*)}`\FF_Z\otimes
Id_{D(V^*)}]
\putsquare<1`1`1`1;2250`1500>(825,900)[Val^{-,sm}_1(W)\otimes
D(W^*) `Val^{-,sm}_{n-j}(W^*)`Val^{-,sm}_1(E)\otimes
D(W^*)`Val^{-,sm}_1(E^*)\otimes D(E)\otimes D(W^*);\FF_W\otimes
Id_{D(W^*)}`i^*\otimes Id_{D(W^*)}`i^\vee_*\otimes
Id_{D(W^*)}`\FF_E\otimes Id_{D(W^*)}]
\putmorphism(0,3300)(1,-1)[``p_*]{900}1r
\putmorphism(2250,3300)(1,-1)[`` p^{\vee *}]{900}1r
\putmorphism(2200,1800)(1,-1)[``"\tilde p_*"]{900}1l
\putmorphism(0,1800)(1,-1)[``"\tilde p^{\vee }_*"]{900}1l\efig$$

The following facets of this cube commute:

$\bullet$ the left facet, by base change Theorem
\ref{T:base-change2};

$\bullet$ the right facet, by base change Theorem
\ref{T:base-change};

$\bullet$ the bottom facet, by Case 2a (since $\dim E=2$);

$\bullet$ the back facet, by Theorem \ref{T:z3};

$\bullet$ the front facet, by Theorem \ref{T:z3}.

These properties and a straightforward diagram chasing imply that
\begin{eqnarray*}
(i^\vee_*\otimes Id_{D(W^*)})\circ\left(p^{\vee *}\circ
(\FF_V\otimes Id_{D(V^*)})-(\FF_W\circ Id_{D(W^*)})\circ
p_*)\right)=0.
\end{eqnarray*}
Since $i^\vee$ is an arbitrary surjection of rank two, Lemma
\ref{L:z6} implies that
$$p^{\vee *}\circ (\FF_V\otimes
Id_{D(V^*)})-(\FF_W\circ Id_{D(W^*)})\circ p_*)=0.$$ This means
that the diagram (\ref{diag-case2b}) is commutative. Thus Case 2b
is proved.


\underline{Case 2c.} Let us consider the general odd case.

Let us fix an {\itshape arbitrary} surjection $q\colon W\surj X$
with $\dim X=n-j+1$. Let us consider the following diagram:
$$\bfig
\putsquare<1`1`1`1;1700`600>(0,750)[Val^{-,sm}_j(V)\otimes
Dens(V^*)` Val^{-,sm}_{n-j}(V^*)`Val^{-,sm}_{j-n+m}(W)\otimes
Dens(W^*)` Val^{-,sm}_{n-j}(W^*);\FF_V\otimes
Id_{Dens(V^*)}`p_*`p^{\vee*}`\FF_W\otimes Id_{Dens(W^*)}]
\putsquare<0`1`1`1;1700`600>(0,150)[\phantom{Val^{-,sm}_{j-n+m}(W)\otimes
Dens(W^*)}`\phantom{Val^{-,sm}_{n-j}(W^*)}`Val^{-,sm}_1(X)\otimes
Dens(X^*)`Val^{-,sm}_{n-j}(X^*);`q_*`q^{\vee*}`\FF_X\otimes
Id_{Dens(X^*)}] \efig$$ By Case 2b, the lower square of this
diagram commutes. Also by Case 2b, the big exterior contour of
this diagram commutes. It follows that
$$q^{\vee *}\circ
\left(p^{\vee ^*}\circ (\FF_V\otimes Id_{Dens(V^*)})-(\FF_W\otimes
Id_{Dens(W^*)})\circ p_*\right)=0.$$ Since the last equality holds
for an arbitrary surjection $q$ of rank $n-j+1$, $q^\vee$ is an
arbitrary injection of an $(n-j+1)$-dimensional subspace. But by
Schneider's theorem \cite{schneider-simple} odd $k$-homogeneous
valuations are uniquely determined by their restrictions to all
$(k+1)$-dimensional subspaces. Hence
$$p^{\vee*}\circ(\FF_V\otimes Id_{Dens(V^*)})-(\FF_W\otimes
Id_{Dens(W^*)})\circ p_*=0.$$ This means that the diagram
(\ref{fourier-id2}) commutes. Proposition is proved. \qed


The following theorem generalizes Lemma \ref{L:com3}.
\begin{theorem}\label{T:com8}
Let $0\leq k\leq n$. Let $\eta\in C^\infty(Gr_k(V),\cs_{k,V})$
($\cs_{k,V}$ was defined in (\ref{D:com1.1}), (\ref{D:com1.2})).
Then
\begin{eqnarray}\label{tttt}
(\FF_V\otimes Id_{Dens(V^*)})\left(\int_{E\in
Gr_k(V)}i_{E*}(\eta(E))\right)=\int_{E\in
Gr_k(V)}(i_E^\vee)^*(\FF_E\otimes Id_{Dens(E^*)})(\eta(E)).
\end{eqnarray}
\end{theorem}
{\bf Proof.} \underline{Case 1.} Let us consider the even case,
i.e. $\eta(E)\in Val_i^{+,sm}(E)\otimes Dens(E^*)$ for any $E\in
Gr_k(V)$, $0\leq i\leq k$.

\underline{Case 1a.} Assume in addition that $i=0$. Thus
$\eta(E)\in Dens(E^*)$ for any $E\in Gr_k(V)$.

We have the canonical identification
\begin{eqnarray}\label{se1}
Dens(E^*)=Dens(V/E)\otimes Dens(V^*).
\end{eqnarray}
It is easy to see from the definitions that under this
identification the map
$$i_{E*}\colon Dens(E^*)\to Val(V)\otimes Dens(V^*)$$
coincides with the map
$$p_E^*\otimes Id_{Dens(V^*)}\colon Dens(V/E)\otimes Dens(V^*)\to
Val(V)\otimes Dens(V^*)$$ where $p_E\colon V\to V/E$ is the
canonical projection.

Let us denote by $\tilde \eta\in
C^\infty(Gr_k(V),\ct_{k,V;k})\otimes Dens(V^*)$ the section
corresponding to $\eta$ under the isomorphism (\ref{se1}). Thus
\begin{eqnarray}\label{se2}
\int_{E\in Gr_k(V)}i_{E*}(\eta(E))=\int_{E\in
Gr_k(V)}(p_E^*\otimes Id_{Dens(V^*)})(\tilde \eta(E)).
\end{eqnarray}

Let us fix a Lebesgue measure $vol_V$ on $V$. Set
$$\hat \eta:=\tilde \eta\cdot vol_V\in
C^\infty(Gr_k(V),\ct_{k,V;k}).$$ Then
$$\int_{E\in Gr_k(V)}i_{E*}(\eta(E))=\left(\int_{E\in
Gr_k(V)}p_E^*(\hat \eta(E))\right)\otimes vol_V^{-1}.$$ Hence by
the definition of the Fourier transform
\begin{eqnarray*}
(\FF_V\otimes Id_{Dens(V^*)})\left(\int_{E\in
Gr_k(V)}i_{E*}(\eta(E))\right)= \left(\int_{E\in
Gr_k(V)}p_{E*}^\vee(\FF_{V/E}(\hat\eta(E)))\right)\cdot
vol_V^{-1}.
\end{eqnarray*}
Hence it suffices to check that for any $E\in Gr_k(V)$
\begin{eqnarray}\label{se3}
\left(p_{E*}^\vee(\FF_{V/E}(\hat\eta(E)))\right)\cdot
vol_V^{-1}=(i_{E*}^\vee)^*(\FF_E\otimes Id_{Dens(E^*)})(\eta(E)).
\end{eqnarray}
To prove the identity (\ref{se3}) let us fix a Lebesgue measure
$vol_E$ on $E$. Let $vol_{V/E}:=\frac{vol_V}{vol_E}$ be the
corresponding Lebesgue measure on $V/E$. It is sufficient to prove
(\ref{se3}) for $\eta(E)=vol_E^{-1}$. Then under the isomorphism
(\ref{se1}) we have
$$\eta(E)=vol_{V/E}\otimes vol_V^{-1}.$$ Hence
$$\hat\eta(E)=vol_{V/E}.$$ We have
$$(\FF_E\otimes Id_{Dens(E^*)}(\eta(E))=vol_E^{-1}\in Dens(E^*).$$
Hence
\begin{eqnarray}\label{se4}
(i_{E*}^\vee)^*(\FF_E\otimes
Id_{Dens(E^*)})(\eta(E))=(i_{E*}^\vee)^*(vol_E^{-1}).
\end{eqnarray}
Next
\begin{eqnarray*}
\FF_{V/E}(\hat \eta(E))=\FF_{V/E}(vol_{V/E})=\chi_{(V/E)^*}\otimes
vol_{V/E}\in Val_0((V/E)^*)\otimes Dens(V/E).
\end{eqnarray*}
Hence
\begin{eqnarray*}
\left(p_{E*}^\vee(\FF_{V/E}(\hat\eta(E)))\right)\cdot
vol_V^{-1}=\left(p_{E*}^\vee(\chi_{(V/E)^*}\otimes
vol_{V/E})\right)\cdot vol_V^{-1}.
\end{eqnarray*}
It is easy to see that the last expression is equal to the right
hand side of (\ref{se4}). This implies (\ref{se3}). Hence Case 1a
is proved.

\underline{Case 1b.} Assume that $i>0$.

Recall that we denote by $\cf_{k-i,k}(V)$ the partial flag
manifold
$$\cf_{k-i,k}(V):=\{(L,E)|\, L\subset E,\, L\in Gr_{k-i}(V),\,
E\in Gr_k(V)\}.$$ For $(L,E)\in \cf_{k-i,k}(V)$ let us denote by
$i_{E,L}$ the imbedding $L\inj E$. Let us denote by $\cg^0$ the
line bundle over $\cf_{k-i,k}(V)$ whose fiber over $(L,E)$ is
equal to $Dens(L^*)$. Let
$$\cg:=\cg^0\otimes |\ome_{\cf_{k-i,i}(V)}|.$$

\begin{claim}\label{Cl:seg3}
For any $\eta\in C^\infty(Gr_k(V),\cs^+_{k,V;i}),k>i,$ there
exists $\xi\in C^\infty(\cf_{k-i,k}(V),\cg)$ such that for any
$E\in Gr_k(V)$
\begin{eqnarray}\label{seg4}
\eta(E)=\int_{L\in Gr_{k-i}(E)}i_{E,L*}(\xi(L,E)).
\end{eqnarray}
\end{claim}
{\bf Proof.} Define the map $C^\infty(\cf_{k-i,k}(V),\cg)\to
C^\infty(Gr_k(V),\cs^+_{k,V;i})$ by the right hand side of
(\ref{seg4}). This map is $GL(V)$-equivariant. Then Irreducibility
Theorem, Casselman-Wallach theorem, and Proposition
\ref{part1-epi} imply the claim. \qed

Let us continue proving Case 1b. We have
\begin{eqnarray*}
\int_{E\in Gr_k(V)}i_{E*}(\eta(E))=\int_{E\in Gr_k(V)}\int_{L\in
Gr_{k-i}(E)}i_{E*}(i_{E,L*}(\xi(L,E)))=\\\int_{E\in
Gr_k(V)}\int_{L\in Gr_{k-i}(E)}i_{L*}(\xi(L,E))= \int_{L\in
Gr_{k-i}(V)}i_{L*}\left(\int_{E\in Gr_k(V),E\supset
L}\xi(L,E)\right).
\end{eqnarray*}
Let us denote by $\hat\xi(L)$ the inner integral in the last
expression. Thus $\hat\xi\in C^\infty(Gr_{k-i}(V),\cs_{k-i,V;0})$.
Hence by Case 1a
\begin{eqnarray}\label{sss5}
(\FF_V\otimes Id_{Dens(V^*)})\left(\int_{E\in
Gr_k(V)}i_{E*}(\eta(E))\right)=\\\label{sss6} (\FF_V\otimes
Id_{Dens(V^*)})\left(\int_{L\in
Gr_{k-i}(V)}i_{L*}(\hat\xi(L))\right)=\\\label{sss7} \int_{L\in
Gr_{k-i}(V)}(i_L^\vee)^*(\FF_L\otimes Id_{Dens(L^*)})(\hat\xi(L)).
\end{eqnarray}
Using again Case 1a we obtain that the right hand side of the
equality (\ref{tttt}) is equal to
\begin{eqnarray*}
\int_{E\in Gr_k(V)}(i_{E}^\vee)^*(\FF_E\otimes
Id_{Dens(E^*)})\left(\int_{L\in
Gr_{k-i}(E)}i_{E,L*}(\xi(L,E))\right)=\\
\int_{E\in Gr_k(V)}(i_{E}^\vee)^*\left(\int_{L\in
Gr_{k-i}(E)}(i_{E,L}^\vee)^*(\FF_L\otimes
Id_{Dens(L^*)})(\xi(L,E))\right)=\\
\int_{E\in Gr_k(V)}\int_{L\in
Gr_{k-i}(E)}(i_L^\vee)^*(\FF_L\otimes Id_{Dens(L^*)})(\xi(L,E))=\\
\int_{L\in Gr_{k-i}(E)}(i_L^\vee)^*(\FF_L\otimes
Id_{Dens(L^*)})(\hat\xi(L))\overset{\mbox{by } (\ref{sss7})}{=}\\
(\FF_V\otimes Id_{Dens(V^*)})\left(\int_{E\in
Gr_k(V)}i_{E*}(\eta(E))\right).
\end{eqnarray*}
Thus Case 1b is proved. Hence Case 1 is proved too.

\underline{Case 2.} Let us consider the odd case, i.e. $\eta(E)\in
Val^{-,sm}_i(E)\otimes Dens(E^*)$ for any $E\in Gr_k(V)$, $1\leq
i\leq k-1$.


Then the equality (\ref{tttt}) becomes equivalent to
\begin{eqnarray*}
\int_{E\in Gr_2(V)}i_{E*}(\eta(E))=\FF_V\left(\int_{E\in
Gr_2(V)}(i_E^\vee)^*(\FF_E^{-1}(\eta(E)))\right).
\end{eqnarray*}
The last equality is a special case of Theorem \ref{T:z1}. Case 2a
is proved.

\underline{Case 2a.} Assume that $i=1$.

Let us fix an arbitrary surjection $p\colon V\surj W$ of rank $k$.
We have
\begin{eqnarray}\label{sem0.5}
p^{\vee*}\left((\FF_V\otimes Id_{Dens(V^*)})\left(\int_{E\in
Gr_k(V)}(\eta(E))\right)\right)\overset{\mbox{Prop.
}(\ref{T:fourier-id})}{=}\\
(\FF_W\otimes Id_{Dens(W^*)})\left(p_*\left(\int_{E\in
Gr_k(V)}i_{E*}i_{E*}(\eta(E))\right)\right)=\\
(\FF_W\otimes Id_{Dens(W^*)})\left(\int_{E\in Gr_k(V)}(p\circ
i_E)_*(\eta(E))\right)=\\
(\FF_W\otimes Id_{Dens(W^*)})\left(\int_{E\in Gr_k(V),\, E\cap Ker
(p)=\{0\}}(p\circ i_E)_*(\eta(E))\right)=\\\label{sem1} \int_{E\in
Gr_k(V),\, E\cap Ker (p)=\{0\}}(p\circ i_E)^{\vee*}(\FF_E\otimes
Id_{Dens(E^*)})(\eta(E))
\end{eqnarray}
where the last equality follows from the fact that $p\circ
i_E\colon E\to W$ is an isomorphism for any $E\in Gr_k(V)$ such
that $E\cap Ker (p)=\{0\}$. Next the expression (\ref{sem1}) is
equal to
\begin{eqnarray}
\int_{E\in Gr_k(V)}(p\circ i_E)^{\vee*}\left((\FF_E\otimes
Id_{Dens(E^*)})(\eta(E))\right)=\\
\int_{E\in Gr_k(V)}p^{\vee*}\left(i_E^{\vee*}\left((\FF_E\otimes
Id_{Dens(E^*)})(\eta(E))\right)\right)=\\\label{sem2}
p^{\vee*}\left(\int_{E\in Gr_k(V)}i_E^{\vee*}\left((\FF_E\otimes
Id_{Dens(E^*)})(\eta(E))\right)\right).
\end{eqnarray}
Thus we have shown that (\ref{sem0.5})$=$(\ref{sem2}). In other
words for any surjection $p\colon V\surj W$ of rank $k$
$$p^{\vee*}(\mbox{l.h.s. of }(\ref{tttt}))=
p^{\vee*}(\mbox{r.h.s. of }(\ref{tttt})).$$ Since $i=1$, this
implies (\ref{tttt}). Thus Case 2a is proved.

\underline{Case 2b.} Assume now that $i>1$.

Let us consider the (infinite dimensional) vector bundle $\ch^0$
over the partial flag manifold $\cf_{k-i+1,k}(V)$ whose fiber over
$(L,E)\in \cf_{k-i+1,k}(V)$ is equal to $Val_1^{-,sm}(L)\otimes
Dens(L^*)$. Let
$$\ch:=\ch^0\otimes |\ome|_{\cf_{k-i+1,k}(V)}|.$$
\begin{claim}\label{Cl:seg1}
For any $\eta\in C^\infty(Gr_k(V),\cs^-_{k,V;i}),\, k>l>1,$ there
exists $\xi\in C^\infty(\cf_{k-i+1,k}(V),\ch)$ such that for any
$E\in Gr_k(V)$
\begin{eqnarray}\label{seg2}
\eta(E)=\int_{L\in Gr_{k-i+1}(E)}i_{E,L*}(\xi(L,E)).
\end{eqnarray}
\end{claim}
{\bf Proof.} The map $C^\infty(\cf_{k-i+1,k}(V),\ch)$ given by the
right hand side of (\ref{seg2}) is $GL(V)$-equivariant. Then the
Irreducibility theorem, the Casselman-Wallach theorem, and
Proposition \ref{part1-epi} imply the claim. \qed

Next we have
\begin{eqnarray*}
\int_{E\in Gr_k(V)}i_{E*}(\eta(E))=\int_{E\in Gr_k(V)}\int_{L\in
Gr_{k-i+1}(E)}i_{E*}(i_{E,L*}(\xi(L,E)))=\\\int_{E\in
Gr_k(V)}\int_{L\in Gr_{k-i+1}(E)}i_{L*}(\xi(L,E))= \int_{L\in
Gr_{k-i+1}(V)}i_{L*}\left(\int_{E\in Gr_k(V),E\supset
L}\xi(L,E)\right).
\end{eqnarray*}
Let us denote by $\hat\xi(L)$ the inner integral in the last
expression. Thus $\hat\xi\in
C^\infty(Gr_{k-i}(V),\cs_{k-i+1,V;1}^-)$. Hence by Case 2a
\begin{eqnarray}\label{sst5}
(\FF_V\otimes Id_{Dens(V^*)})\left(\int_{E\in
Gr_k(V)}i_{E*}(\eta(E))\right)=\\\label{sst6} (\FF_V\otimes
Id_{Dens(V^*)})\left(\int_{L\in
Gr_{k-i+1}(V)}i_{L*}(\hat\xi(L))\right)=\\\label{sst7} \int_{L\in
Gr_{k-i+1}(V)}(i_L^\vee)^*(\FF_L\otimes
Id_{Dens(L^*)})(\hat\xi(L)).
\end{eqnarray}
Using again Case 2a we obtain that the right hand side of the
equality (\ref{tttt}) is equal to
\begin{eqnarray*}
\int_{E\in Gr_k(V)}(i_{E}^\vee)^*(\FF_E\otimes
Id_{Dens(E^*)})\left(\int_{L\in
Gr_{k-i+1}(E)}i_{E,L*}(\xi(L,E))\right)=\\
\int_{E\in Gr_k(V)}(i_{E}^\vee)^*\left(\int_{L\in
Gr_{k-i+1}(E)}(i_{E,L}^\vee)^*(\FF_L\otimes
Id_{Dens(L^*)})(\xi(L,E))\right)=\\
\int_{E\in Gr_k(V)}\int_{L\in
Gr_{k-i+1}(E)}(i_L^\vee)^*(\FF_L\otimes Id_{Dens(L^*)})(\xi(L,E))=\\
\int_{L\in Gr_{k-i+1}(E)}(i_L^\vee)^*(\FF_L\otimes
Id_{Dens(L^*)}(\hat\xi(L))\overset{\mbox{by } (\ref{sst7})}{=}\\
(\FF_V\otimes Id_{Dens(V^*)})\left(\int_{E\in
Gr_k(V)}i_{E*}(\eta(E))\right).
\end{eqnarray*}
Thus Case 2a is proved. Hence Case 2 is proved too, and Theorem
\ref{T:com8} is proved. \qed

The following immediate reformulation of Theorem \ref{T:com8} will
be useful in the proof of Theorem \ref{T:hfh1} below.
\begin{corollary}\label{Cor:com8.7}
Let $0\leq k\leq n$. Let $\tilde\eta\in
C^\infty(Gr_k(V),\cs_{k,V})\otimes Dens(V)$. Then
\begin{eqnarray}
\FF_V\left(\int_{E\in Gr_k(V)}(i_{E*}\otimes
Id_{Dens(V)})(\tilde\eta(E))\right)=\\\label{wow}\int_{E\in
Gr_k(V)}(i_E^{\vee*}\otimes Id_{Dens(V)})\left((\FF_E\otimes
Id_{Dens(V/E)})(\tilde\eta(E))\right)
\end{eqnarray}
where $i_{E*}\otimes Id_{Dens(V)}$ in the left hand side is the
map $\left(Val^{sm}(E)\otimes Dens(E^*)\right)\otimes Dens(V)\to
Val^{sm}(V)$, and $\FF_E\otimes Dens(V/E)$ in the right hand side
is the map
\begin{eqnarray*}Val^{sm}(E)\otimes Dens(E^*)\otimes
Dens(V)=Val^{sm}(E)\otimes Dens(V/E)\overset{\FF_E\otimes
Id_{Dens(V/E)}}{\to}
\\\left(Val^{sm}(E^*)\otimes Dens(E)\right)\otimes
Dens(V/E)=Val^{sm}(E^*)\otimes Dens(V).\end{eqnarray*}
\end{corollary}



\begin{theorem}[Plancherel formula]\label{T:com02}
$$(\FF_{V^*}\otimes Id_{Dens(V)})\circ \FF_V=\ce_V.$$
\end{theorem}
{\bf Proof.} In the even case this theorem is precisely Lemma
\ref{L:com01.1}. Let us consider the odd $i$-homogeneous case
$$(\FF_{V^*}\otimes Id_{Dens(V)})\circ \FF_V\colon
Val_i^{-,sm}(V)\to Val_i^{-,sm}(V).$$ If $i=1$ the result is just
Lemma \ref{L:com01}. If $i=n-1$ then the result follows
immediately from the definition of the Fourier transform
$\FF_V\colon Val_{n-1}^{-,sm}(V)\to Val_1^{-,sm}(V^*)\otimes
Dens(V)$ since $\FF_V:=-(\FF_{V^*}\otimes Id_{Dens(V)})^{-1}$.

Let us assume that $2\leq i\leq n-2$. Fix $\phi\in
Val_i^{-,sm}(V)$. By Proposition \ref{P:w5} there exists $\xi\in
C^\infty(Gr_{n-i-1}(V),\ct_{i+1,V;i}^-)$ such that
$\phi=\Xi_{i+1,V}(\xi)$. By the definition of $\FF_V$ in this
case,
$$\FF_V(\xi)=\int_{F\in Gr_{n-i-1}(V)}p_{F*}^{\vee}(\FF_{V/F}(\xi(F))).$$

By Theorem \ref{T:com8} we have
\begin{eqnarray}
(\FF_V\otimes Id_{Dens(V^*)})(\FF_V(\phi))=\\
\int_{F\in Gr_{n-i-1}(V)}p_F^*\left((\FF_{(V/F)^*}\otimes
Id_{Dens(V/F)})(\FF_{V/F}(\xi(F)))\right)=\\\label{oo} -\int_{F\in
Gr_{n-i-1}(V)}p_F^*(\xi(F))
\end{eqnarray}
where the last equality holds by the proved above case
corresponding to $n=i+1$. But the expression (\ref{oo}) is equal
to $-\Xi_{i+1,V}(\xi)=-\phi$. Hence theorem is proved. \qed


Let us fix now a Euclidean metric on $V$. This induces
isomorphisms
\begin{eqnarray}\label{pppp}
V^*\tilde\to V, Dens(V)\tilde\to\CC.\end{eqnarray} Under these
identifications $\FF_V\colon Val^{sm}(V)\tilde \to Val^{sm}(V)$.
From Theorem \ref{T:com02} we easily deduce
\begin{corollary}\label{C:plancher}
Under the identifications (\ref{pppp}) one has $\FF_V^4=Id$.
Moreover $\FF_V^2\ne Id$ provided $\dim V>1$.
\end{corollary}

\subsection{Homomorphism property of the Fourier transform in
higher dimensions.}\label{Ss:homomorphism-F-high}
\begin{theorem}\label{T:hfh1}
The map
$$\FF_V\colon Val^{sm}(V)\to Val^{sm}(V^*)\otimes Dens(V)$$
is an isomorphism of algebras.
\end{theorem}
{\bf Proof.} It remains to show that $\FF_V$ is a {\itshape
homomorphism} of algebras, i.e.
\begin{eqnarray}\label{pr-conv}
\FF_V(\phi\cdot \psi)=\FF_V(\phi)\ast \FF_V(\psi)
\end{eqnarray}
for any $\phi,\psi\in Val^{sm}(V)$.  We will consider a number of
cases.

\underline{Case 1.} Assume that $\phi$ and $\psi$ are even.

This case was proved by Bernig and Fu \cite{bernig-fu}. We will
not prove this case here though it can be proved similarly to Case
2 below, but simpler.

\underline{Case 2.} Assume that $\phi\in Val^{+,sm}_i(V),\,
\psi\in Val^{-,sm}_j(V)$.

We will need a general lemma.
\begin{lemma}\label{L:prod-push}
Let $F\subset V$ be a linear subspace. Let $i_F\colon F\inj V$,
$p_F\colon V\surj V/F$ be the canonical imbedding and projection
maps. Let $E\subset V$ be another linear subspace, and let
$p_E\colon V\surj V/E$ be the canonical projection.

Let $\mu\in Dens(V/F)\subset Val^{sm}(V/F)$, $\zeta\in
Val^{sm}(V/E)$. Then
\begin{eqnarray}\label{forl}
(p_F\times p_E)^*(\mu\boxtimes \zeta)=(i_{F*}\otimes
Id_{Dens(V)})((p_E\circ i_F)^*(\zeta)\otimes \mu)
\end{eqnarray}
where
\begin{eqnarray*}
p_F\times p_E\colon V\to V/F\times V/E,\\
(p_E\circ i_F)^*(\zeta)\otimes\mu\in Val(F)\otimes
Dens(V/F)=\left(Val(F)\otimes Dens(F^*)\right)\otimes Dens(V),\\
i_{F*}\otimes Id_{Dens(V)}\colon\left(Val(F)\otimes
Dens(F^*)\right)\otimes Dens(V)\to Val(V).
\end{eqnarray*}
\end{lemma}
{\bf Proof.} By the McMullen's conjecture we may assume that
$$\zeta(\bullet)=vol_{V/E}(\bullet+A)$$
where $vol_{V/E}$ is a Lebesgue measure on $V/E$, and $A\in
\ck^{sm}(V/E)$. Let us compute first the left hand side of
(\ref{forl}). One has
\begin{eqnarray}\label{op1}
(p_F\times p_E)^*(\mu\boxtimes \zeta)(K)=\\
(\mu\boxtimes\zeta)((p_F\times p_E)(K))=\\\label{op2} \int_{z\in
V/F}vol_{V/E}\left(((p_F\times p_E)(K)\cap(\{z\}\times
V/E))+A\right)d\mu(z).
\end{eqnarray}
Let us compute now the right hand side of (\ref{forl}). One has
\begin{eqnarray*}
(i_{F*}\otimes Id_{Dens(V)})\left((p_E\circ
i_F)^*(\zeta)\otimes \mu\right)(K)=\\
\int_{z\in V/F}\left((p_E\circ i_F)^*\zeta\right)(K\cap z)d
\mu(z)=\\
\int_{z\in V/F}\zeta\left((p_E\circ i_F)(K\cap z)\right)d \mu(z).
\end{eqnarray*}
Hence
\begin{eqnarray}\label{op3}
(i_{F*}\otimes Id_{Dens(V)})\left((p_E\circ i_F)^*(\zeta)\otimes
\mu\right)(K)=\\\label{op4} \int_{z\in V/F}\zeta(p_E(K\cap
z))d\mu(z)=\\\label{op5} \int_{z\in V/F}vol_{V/E}(p_E(K\cap
z)+A)d\mu(z).
\end{eqnarray}
But obviously
$$\{z\}\times p_E(K\cap
z)=(p_F\times p_E)(K)\cap (\{z\}\times V/E).$$ Substituting this
into (\ref{op5}) and comparing with (\ref{op2}) we get that
$$(\ref{op1})=(\ref{op3}).$$
Thus lemma is proved. \qed

Let us return now to the proof of Case 2 of the theorem. First
observe that if $i+j>n-1$ then both sides of (\ref{pr-conv})
vanish. Thus we will assume that $i+j\leq n-1$.

\underline{Case 2a.} Let us assume in addition that $i+j=n-1$.

By Proposition \ref{P:w5} there exist $\mu\in
C^\infty(Gr_{n-i}(V),\ct^+_{i,V;i})$ and $\zeta\in
C^\infty(Gr_{i}(V),\ct^-_{j+1,V;j})$ such that
$\phi=\Xi_{i,V}(\mu),\psi=\Xi_{j+1,V}(\zeta).$ In other words
\begin{eqnarray*}
\phi=\int_{E\in Gr_{n-i}(V)}p_E^*(\mu(E)),\\
\psi=\int_{F\in Gr_{i}(V)}p_F^*(\zeta(F)).
\end{eqnarray*}
Then we have
\begin{eqnarray}\label{2a1}
\FF_V(\phi\cdot \psi)\overset{\mbox{Prop. }\ref{exter-pull}}
{=}\\\label{2a2}\FF_V\left(\int_{E\in Gr_{n-i}(V)}\int_{F\in
Gr_{i}(V)}(p_F\times p_E)^*(\mu(F)\boxtimes
\zeta(E))\right)\overset{\mbox{Lemma }
\ref{L:prod-push}}{=}\\\label{2a3} \FF_V\left(\int_{E\in
Gr_{n-i}(V)}\int_{F\in Gr_{i}(V)}(i_{F*}\otimes
Id_{Dens(V)})\left((p_E\circ i_F)^*(\zeta(E))\otimes
\mu(F)\right)\right)
\end{eqnarray}
Let us denote
\begin{eqnarray}\label{locus}
\cz:=\{(E,F)\in Gr_{n-i}(V)\times Gr_i(V)|\, E\cap F=0\}.
\end{eqnarray}
Thus $\cz\subset Gr_{n-i}(V)\times Gr_i(V)$ is an open dense
subset whose complement is a (singular) subvariety of positive
codimension. Then clearly
\begin{eqnarray*}
(\ref{2a3})=\FF_V\left(\int_\cz(i_{F*}\otimes
Id_{Dens(V)})\left((p_E\circ i_F)^*(\zeta(E))\otimes
\mu(F)\right)\right).
\end{eqnarray*}
Let us denote for brevity
$$\ome(E,F):=(p_E\circ i_F)^*(\zeta(E))\otimes
\mu(F).$$ It is easy to see that

$\bullet$ $\ome$ is a {\itshape continuous} section of the bundle
$\ct^+_{i,V;i}\boxtimes \ct^-_{j+1,V;j}$ over $Gr_{n-i}(V)\times
Gr_i(V)$;

$\bullet$ $\ome$ is smooth over $\cz$.

Let us fix a sequence $\gamma_N\colon Gr_{n-i}(V)\times Gr_i(V)\to
[0,1]$, $N\in \NN$, of $C^\infty$-smooth functions vanishing in a
neighborhood (depending on $N$) of the complement of $\cz$, and
converging uniformly on compact subsets of $\cz$ to the function 1
with all partial derivatives.

From the mentioned properties of $\ome$ and the choice of a
sequence $\{\gamma_N\}$ it easily follows that
$$\lim_{N\to \infty} \int_\cz \gamma_N\cdot (i_{F*}\otimes Id_{Dens(V)})\ome(E,F)=
\int_\cz (i_{F*}\otimes Id_{Dens(V)})\ome(E,F)$$ where the
convergence is understood in the space $Val^{sm}(V)$. Hence
\begin{eqnarray*}
(\ref{2a3})=\lim_{N\to \infty}\FF_V\left(\int_\cz
\gamma_N\cdot(i_{F*}\otimes Id_{Dens(V)})\left((p_E\circ
i_F)^*(\zeta(E))\otimes \mu(F)\right)\right)=\\
\lim_{N\to \infty}\int_{E\in Gr_{n-i}(V)}\FF_V\left(\int_{F\in
Gr_i(V)}\gamma_N\cdot(i_{F*}\otimes Id_{Dens(V)})\left((p_E\circ
i_F)^*(\zeta(E))\otimes \mu(F)\right)\right).
\end{eqnarray*}

We may apply Corollary \ref{Cor:com8.7} to the last expression; it
is equal to
\begin{eqnarray*}
\lim_{N\to \infty}\int_{E\in Gr_{n-i}(V)}\int_{F\in
Gr_i(V)}\gamma_N\cdot (i_F^{\vee*}\otimes
Id_{Dens(V)})\left((\FF_F\otimes Id_{Dens(V/F)})((p_E\circ
i_F)^*(\zeta(E))\otimes \mu(F))\right).
\end{eqnarray*}
The last limit is clearly equal to
\begin{eqnarray}
\int_{(E,F)\in \cz} (i_F^{\vee*}\otimes
Id_{Dens(V)})\left((\FF_F\otimes Id_{Dens(V/F)})((p_E\circ
i_F)^*(\zeta(E))\otimes \mu(F))\right)=\\\label{2a4}
\int_{(E,F)\in \cz} (i_F^{\vee*}\otimes
Id_{Dens(V)})\left(\FF_F\left((p_E\circ
i_F)^*(\zeta(E))\right)\otimes \mu(F)\right).
\end{eqnarray}
Since for $(E,F)\in \cz$ the map $p_E\circ i_F\colon F\to V/E$ is
an isomorphism, we have
$$\FF_F\circ (p_E\circ i_F)^*=(p_E\circ
i_F)^\vee_*\circ \FF_{V/E}.$$ Hence
\begin{eqnarray*}
(\ref{2a4})=\int_\cz(i_F^{\vee*}\otimes
Id_{Dens(V)})\left(i^\vee_{F*}(p^\vee_{E*}(\FF_{V/E}\zeta(E)))\otimes
\mu(F)\right)=\\
\int_{F\in Gr_i(V)}(i_F^{\vee*}\otimes
Id_{Dens(V)})\left(i^\vee_{F*}\left(\int_{E\in
Gr_{n-i}(V)}(p^\vee_{E*}(\FF_{V/E}\zeta(E)))\right)\otimes
\mu(F)\right)=\\
\int_{F\in Gr_i(V)}(i_F^{\vee*}\otimes
Id_{Dens(V)})\left(i^\vee_{F*}(\FF_V\psi)\otimes\mu(F)\right).
\end{eqnarray*}
To summarize the above computation, we have obtained the equality
\begin{eqnarray}\label{2a5}
\FF_V(\phi\cdot\psi)=\int_{F\in Gr_i(V)}(i_F^{\vee*}\otimes
Id_{Dens(V)})\left(i^\vee_{F*}(\FF_V\psi)\otimes\mu(F)\right).
\end{eqnarray}

Hence to finish the proof of Case 2a it remains to prove the
following lemma.
\begin{lemma}\label{L:2a6}
Let $\mu\in C^\infty(Gr_i(V),\ct_{i,V;i})$. Let $\phi=\int_{F\in
Gr_i(V)}p_F^*(\mu(F))$. Let $\xi\in Val^{sm}(V^*)\otimes Dens(V)$.
Then
\begin{eqnarray}\label{2a7}
\int_{F\in Gr_i(V)}(i_F^{\vee*}\otimes
Id_{Dens(V)})\left(i^\vee_{F*}(\xi)\otimes\mu(F)\right)=\xi\ast
\FF_V(\phi).
\end{eqnarray}
\end{lemma}
{\bf Proof.} First recall that $Dens(V/F)=Dens(F^*)\otimes
Dens(V)$. Under this identification, for any $F\in Gr_i(V)$ we
have $\mu(F)\in Dens(F^*)\otimes Dens(V)\otimes
|\ome_{Gr_i(V)}|\big|_F$. Next we have
$$\FF_V(\phi)=\int_{F\in Gr_i(V)}(i_F^{\vee*}\otimes
Id_{Dens(V)})(\mu(F)).$$ Let us fix a Lebesgue measure $vol_V$ on
$V$. Then by Proposition \ref{exter-pull}
\begin{eqnarray*}\xi\boxtimes \FF_V(\phi)=\int_{F\in
Gr_i(V)}\xi\boxtimes (i_F^{\vee*}\otimes Id_{Dens(V)})(\mu(F))=\\
\int_{F\in Gr_i(F)}\xi\boxtimes \left(i_F^{\vee*}(\mu(F)\otimes
vol_V^{-1})\otimes vol_V\right).
\end{eqnarray*}
Hence, using Proposition \ref{P:r2}, we get
\begin{eqnarray*}
\xi\ast\FF_V(\phi)=a_*(\xi\boxtimes \FF_V(\phi))=\\
\int_{F\in Gr_i(F)}a_*\left(\xi\boxtimes
\left(i_F^{\vee*}(\mu(F)\otimes vol_V^{-1})\otimes
vol_V\right)\right)
\end{eqnarray*}
where $a\colon V^*\times V^*\to V^*$ is the addition map. Thus in
order to prove the lemma, it suffices to show that for any $F\in
Gr_i(V)$, any $\xi\in Val^{sm}(V^*)\otimes Dens(V)$, and any
$\mu\in Dens(V/F)=Dens(F^*)\otimes Dens(V)$ one has
\begin{eqnarray}\label{14-1}
a_*\left(\xi\boxtimes \left(i_F^{\vee*}(\mu\otimes
vol_V^{-1})\otimes vol_V\right)\right)=(i_F^{\vee*}\otimes
Id_{Dens(V)})\left(i^\vee_{F*}(\xi)\otimes\mu\right).
\end{eqnarray}

Let us fix also a Lebesgue measure $vol_F$ on $F$. Let us denote
$vol_{V/F}:=\frac{vol_V}{vol_F}$ the corresponding Lebesgue
measure on $V/F$. It is enough to prove the equality (\ref{14-1})
for $\mu=vol_{V/F}$. Then $\mu\otimes vol_V^{-1}=vol_F^{-1}\in
Dens(F^*)$. Furthermore, by the McMullen's conjecture we may
assume that $\xi(\bullet)=vol_V^{-1}(\bullet +A)\otimes vol_V$
where $A\in \ck^{sm}(V^*)$.

Let us fix $K\in \ck(V^*)$. Then we have
\begin{eqnarray}
(\mbox{r.h.s. of } (\ref{14-1}))(K)=(i^\vee_{F*}(\xi)\otimes
\mu)(i^\vee_F(K))=\\\left(vol_F^{-1}(i^\vee_F(K)+i_F^\vee(A))\otimes
vol_F\right)\otimes vol_{V/F}=\\\label{14-2}
vol_F^{-1}(i_F^\vee(K+A))\otimes vol_V.
\end{eqnarray}

On the other hand we have
\begin{eqnarray}
\xi\boxtimes \left(i_F^{\vee*}(\mu\otimes vol_V^{-1})\otimes
vol_V\right)=\\\label{14-3} \left(vol_V^{-1}(\bullet +A)\otimes
vol_V\right)\boxtimes \left(i_F^{\vee*}(vol_F^{-1})\otimes
vol_V\right).
\end{eqnarray}
Let us fix $S\in \ck(Ker(i_F^\vee))=\ck((V/F)^*)$ with
$vol_{V/F}^{-1}(S)=1$. Then by Lemma \ref{L:aLemma} we have
\begin{eqnarray}\label{14-4}
(i_F^{\vee
*}(vol_F^{-1}))(K)=vol_F^{-1}(i_F^\vee(K))=\frac{1}{(n-i)!}
\frac{d^{n-i}}{d\eps^{n-i}}\big|_{\eps=0}vol_V^{-1}(K+\eps S).
\end{eqnarray}
Substituting (\ref{14-4}) into (\ref{14-3}) and using the
definition of pushforward we get
\begin{eqnarray}\label{14-5}
(\mbox{l.h.s. of }(\ref{14-1}))= vol_V^{-1}(i_F^\vee(K+A))\otimes
vol_V.
\end{eqnarray}
Comparing (\ref{14-5}) and (\ref{14-2}) we conclude the equality
(\ref{14-1}). Hence Lemma \ref{L:2a6} is proved. \qed

Case 2a is proved as well.

\underline{Case 2b.} Let us assume now that $i+j+1<n$.

By Schneider's theorem \cite{schneider-simple} every odd
$(i+j+1)$-homogeneous valuation is uniquely determined by its
restrictions to all $(i+j+1)$-dimensional subspaces. Hence it is
enough to check that for any $(i+j+1)$-dimensional subspace
$i\colon L\inj V$ one has
$$i^*(\phi\cdot \psi)=i^*(\FF_V^{-1}(\FF_V\phi\ast \FF_V\psi)).$$
Let us compute the right hand side of the above equality:
\begin{eqnarray*}
i^*\left(\FF_V^{-1}(\FF_V\phi\ast
\FF_V\psi)\right)\overset{\mbox{Thm. } \ref{T:z3}}{=}
\FF_L^{-1}\left(i_*^\vee(\FF_V\phi\ast\FF_V\psi)\right)\overset{\mbox{Prop.
} \ref{P:homom-of-push}}{=}\\
\FF_L^{-1}(i_*^\vee(\FF_V\phi)\ast
i_*^\vee(\FF_V\phi))\overset{\mbox{Thm. }\ref{T:z3}}{=}
\FF_L^{-1}(\FF_L(i^*\phi)\ast \FF_L(i^*\psi))\overset{\mbox{Case
2a}}{=}\\
i^*\phi\cdot i^*\psi=i^*(\phi\cdot \psi).
\end{eqnarray*}
Thus Case 2b is proved. Hence Case 2 is proved too.

\underline{Case 3.} Assume that $\phi,\psi\in Val_1^{-,sm}(V)$.

We will use the homomorphism property of the Fourier transform in
two dimensions proved in Section \ref{two-dim-isomorphism}. By
Klain's theorem \cite{klain} it is enough to show that for any
2-dimensional space $E$ and an imbedding $i\colon E\inj V$ one has
\begin{eqnarray}\label{general-isomor}
i^*(\phi\cdot
\psi)=i^*\left(\FF_V^{-1}(\FF_V\phi\ast\FF_V\psi)\right).
\end{eqnarray}
Let us compute the right hand side of (\ref{general-isomor}):
\begin{eqnarray*}
i^*\left(\FF_V^{-1}(\FF_V\phi\ast\FF_V\psi)\right)\overset{\mbox{Thm.
}\ref{T:z3}}{=}
\FF_L^{-1}\left(i^\vee_*(\FF_V\phi\ast\FF_V\psi)\right)\overset{\mbox{Prop.
}\ref{P:homom-of-push}}{=}\\
\FF_L^{-1}\left(i^\vee_*(\FF_V\phi)\ast
i^\vee_*(\FF_V\psi)\right)\overset{\mbox{Thm. }\ref{T:z3}}{=}
\FF_L^{-1}\left(\FF_L(i^*\phi)\ast
\FF_L(i^*\psi)\right)\overset{\mbox{Thm. }\ref{T:ISO}(2)}{=}\\
i^*\phi\cdot i^*\psi=i^*(\phi\cdot \psi).
\end{eqnarray*}
Case 3 is proved.

\underline{Case 4.} Let us prove the equality (\ref{pr-conv}) in
general.

The only remaining case is $\phi\in Val^{-,sm}_i(V),\,\psi\in
Val^{-,sm}_j(V)$, and either $i>1$ or $j>1$. For any $k\geq 1$ the
subspace $Val_1^{-,sm}(V)\cdot Val_{k-1}^{+,sm}(V)$ is dense in
$Val_k^{-,sm}(V)$ by the Irreducibility Theorem. Hence we may
assume that $\phi=\phi^-\cdot \phi^+$ where $\phi^-\in
Val_1^{-,sm}(V)$, $\phi^+\in Val_{i-1}^{+,sm}(V)$; and similarly
$\psi=\psi^-\cdot \psi^+$ where $\psi^-\in Val_1^{-,sm}(V)$,
$\psi^+\in Val_{j-1}^{+,sm}(V)$. Then we have
\begin{eqnarray*}
\FF_V(\phi\cdot
\psi)=\FF_V((\phi^+\cdot\psi^+)\cdot(\phi^-\cdot\psi^-))\overset{\mbox{Case
1}}{=}\\
\FF_V(\phi^+\cdot\psi^+)\ast\FF_V(\phi^-\cdot\psi^-)\overset{\mbox{Cases
1,3}}{=}\FF_V(\phi^+)\ast\FF_V(\psi^+)\ast\FF_V(\phi^-)\ast\FF_V(\psi^-)\overset{\mbox{Case
2}}{=}\\ \FF_V(\phi^+\cdot\phi^-)\ast \FF_V(\psi^+\cdot
\psi^-))=\FF_V(\phi)\ast\FF_V(\psi).
\end{eqnarray*}
Theorem \ref{T:hfh1} is proved. \qed

\begin{remark}\label{R:non-unique}
The Fourier transform $\FF_V$ we have constructed is not quite
canonical. More precisely, let us fix $n>1$. Let $\cc_n$ denote
the category whose objects $Ob(\cc_n)$ are $n$-dimensional real
vector spaces, and morphisms between them are linear isomorphisms.
Assume that for any object $V$ of $\cc_n$ we are given an
isomorphism $\FF_V\colon Val^{sm}(V)\tilde\to Val^{sm}(V^*)\otimes
Dens(V)$ of linear topological spaces such that
\newline
$\bullet$ for any morphism $f\colon V\to W$ in $\cc_n$ (i.e. $f$
is just a linear isomorphism) the following diagram is commutative
\begin{eqnarray*}
\square<1`1`1`1;1300`400>[Val^{sm}(V)`Val^{sm}(W)`Val^{sm}(V^*)\otimes
Dens(V)`Val^{sm}(W^*)\otimes Dens(W); `\FF_V`\FF_W`]
\end{eqnarray*}
where the horizontal arrows are obvious isomorphisms induced by
the isomorphisms $V\overset{f}{\to}W$ and
$W^*\overset{f^\vee}{\to}V^*$ where $f^\vee$ is the dual of $f$;
\newline
$\bullet$ for any $V\in Ob(\cc_n)$ the map $\FF_V$ is an
isomorphism of algebras when the source is equipped with the
product and the target with the convolution;
\newline
$\bullet$ for any $V\in Ob(\cc_n)$ one has the Plancherel type
formula as in Theorem \ref{T:main1}(3).

Then one can show that there exist exactly four families of maps
$\{\FF_V\}_{V\in Ob(\cc_n)}$ satisfying the above conditions. The
difficult part (which is the main subject of this article) is to
prove existence of at least one such a family.
\end{remark}

\section{A hard Lefschetz type theorem for
valuations.}\label{lefschetz}\setcounter{subsection}{1}\setcounter{theorem}{0}
\setcounter{equation}{0} Let $V$ be a {\itshape Euclidean}
$n$-dimensional space. Let $V_1\in Val_1(V)$ denote the first
intrinsic volume (see e.g. \cite{schneider-book}, p. 210). This
valuation is invariant under the orthogonal group and it is
smooth. The main result of this section is the following theorem.
\begin{theorem}[hard Lefschetz type theorem]\label{T:lefschetz}
Let $0\leq i<n/2$. Then the map
$$Val_i^{sm}(V)\to Val_{n-i}^{sm}(V)$$
given by $\phi\mapsto (V_1)^{n-2i}\cdot \phi$ is an isomorphism.
\end{theorem}
\begin{remark}\label{R:lefschetz-even}
For even valuations this result was proved first by the author in
\cite{alesker-gafa-sem-04}.
\end{remark}

Before we prove Theorem \ref{T:lefschetz} we need some preparations.
In \cite{alesker-jdg-03} the author has introduced an operator
$$\Lambda\colon Val(V)\to Val(V)$$ defined by
$(\Lambda\phi)(K)=\frac{d}{d\eps}\big|_{\eps=0}\phi(K+\eps \cdot
D)$ for any $\phi\in Val(V),\, K\in \ck(V)$. Note that
$\phi(K+\eps \cdot D)$ is a polynomial in $\eps\geq 0$ of degree
at most $n$ by a result of McMullen \cite{mcmullen-euler}. The
operator $\Lambda$ decreases the degree of homogeneity by 1. We
are going to use the following theorem which was proved by the
author \cite{alesker-jdg-03} for even valuations and by Bernig and
Br\"ocker \cite{bernig-brocker} in general.
\begin{theorem}\label{T:hlold}
Let $n\geq i>n/2$. The operator $$\Lambda^{2i-n}\colon
Val_i^{sm}(V)\to Val_{n-i}^{sm}(V)$$ is an isomorphism.
\end{theorem}

Next, the Euclidean metric on $V$ induces the identifications
$V\tilde\to V^*$ and $Dens(V)\tilde\to \CC$. Under these
identifications the Fourier transform acts $\FF_V\colon
Val^{sm}(V)\tilde \to Val^{sm}(V)$. We will need the following
lemma which was observed by Bernig and Fu in \cite{bernig-fu},
Corollary 1.9, in the case of even valuations.
\begin{lemma}\label{L:hlopers}
For any $\phi\in Val^{sm}(V)$ one has
$$V_1\cdot \phi=\kappa(\FF_V^{-1}\circ \Lambda\circ \FF_V)(\phi)$$
where $\kappa$ is a non-zero constant depending on $n$ only.
\end{lemma}
{\bf Proof.} The proof is essentially the same as in the even case
\cite{bernig-fu}, once one has the Fourier transform. By the
homomorphism property of the Fourier transform we have
\begin{eqnarray}\label{opa4}
\FF_V(V_1\cdot\phi)=\FF_V(V_1)\ast \FF_V(\phi). \end{eqnarray}
Observe that $\FF_V(V_1)$ is an $O(n)$-invariant valuation
homogeneous of degree $n-1$. Hence by the Hadwiger
characterization theorem \cite{hadwiger-book} it must be
proportional to the $(n-1)$-th intrinsic volume $V_{n-1}$, which
is proportional to the valuation $K\mapsto
\frac{d}{d\eps}\big|_{0}vol(K+\eps D)$. Next observe that for any
$A\in \ck^{sm}(V)$ and any $\phi\in Val^{sm}(V)$
\begin{eqnarray}\label{opa1}
vol(\bullet+A)\ast \phi=\phi(\bullet+A).
\end{eqnarray}
Indeed the equality (\ref{opa1}) is easily checked for $\phi$ of
the form $\phi(\bullet)=vol(\bullet+B)$, and the general case
follows from the McMullen's conjecture. Hence $\FF_V(V_1)\ast
\FF_V(\phi)$ is proportional to
\begin{eqnarray}\label{opa3}
\frac{d}{d\eps}\big|_{0}(\FF_V\phi)(\bullet+\eps
D)=\Lambda(\FF_V\phi). \end{eqnarray} Then Lemma \ref{L:hlopers}
follows from (\ref{opa4}), (\ref{opa3}). \qed

{\bf Proof of Theorem \ref{T:lefschetz}.} It follows immediately
from Theorem \ref{T:hlold} and Lemma \ref{L:hlopers}. \qed


\section{Appendix: a remark on exterior product on
valuations.}\label{appendix-ext-product}\setcounter{subsection}{1}\setcounter{theorem}{0}
\setcounter{equation}{0} In this appendix we will prove a slightly
more refined statement on the exterior product of valuations than
it appears in \cite{alesker-gafa-04}. Though for the purposes of
this article we need only the case of translation invariant
valuations, we will prove the result in a greater generality of
polynomial valuations following \cite{alesker-gafa-04}.

Let us remind the definition of a {\itshape polynomial} valuation
introduced by Khovanskii and Pukhlikov in
\cite{khovanskii-pukhlikov1}. Let $V$ be an $n$-dimensional real
vector space.
\begin{definition}
A valuation $\phi$ is called polynomial of degree at most $d$ if
for every $K\in {\cal K}(V)$ the function $x\mapsto \phi (K+x)$ is
a polynomial on $V$ of degree at most $d$.
\end{definition}
Note that valuations polynomial of degree 0 are just translation
invariant valuations. Polynomial valuations have many nice
combinatorial-algebraic properties (\cite{khovanskii-pukhlikov1},
\cite{khovanskii-pukhlikov2}).

 Let $\pvd$ denote the
space of continuous valuations on $V$ which are polynomial of
degree at most $d$. It is a Fr\'echet space (in fact a Banach
space) with the topology of uniform convergence on compact subsets
of $\ck(V)$. Let $\ond$ denote the (finite dimensional) space of
$n$-densities on $V$ with polynomial coefficients of degree at
most $d$ (clearly $\ond$ is canonically isomorphic to $(\oplus
_{i=0}^{d}Sym^{i}V^*) \otimes |\wedge ^{n}V^*|$ where $|\wedge
^nV^*|$ denotes the space of Lebesgue measures on $V$).

The group $GL(V)$ acts naturally on $\pvd$ as usual:
$(g\phi)(K)=\phi(g^{-1}K)$. This action is continuous.  The
subspace of $GL(V)$-smooth vectors is denoted by $PVal_d^{sm}(V)$.

For a vector space $U$, a smooth measure $\mu$ on $U$,and $A\in
\ck(U)$ let us denote by $\mu_A$ the valuation $[K\mapsto
\mu(K+A)]$. Now let us state the main result of this appendix
which refines Proposition 1.10 from \cite{alesker-gafa-04}.
\begin{proposition}\label{P:ext-prod}
Let $V,W$ be finite dimensional real vector spaces. There exists a
continuous bilinear map
$$PVal_d^{sm}(V)\times PVal_{d'}(W)\to PVal_{d+d'}(V\times W)$$
which is uniquely characterized by the property that for any
polynomial measures $\mu,\nu$ on $V,W$ respectively, and any $A\in
\ck(V),\, B\in \ck^{sm}(W)$ one has
$$(\mu_A,\nu_B)\mapsto (\mu\boxtimes\nu)_{A\times B}$$
where $\mu\boxtimes\nu$ denotes the usual product measure. This
map is called the exterior product and is denoted by $\boxtimes$.
\end{proposition}
\begin{remark}
An important difference of this proposition in comparison to
\cite{alesker-gafa-04} is that now we can consider the exterior
product of a smooth valuation by a continuous one (and not just a
product of two smooth valuations).
\end{remark}

Before we prove this proposition let us introduce more notation
and remind some constructions from \cite{alesker-gafa-04}. Let us
denote by $\pl$ the manifold of oriented lines passing through the
origin in $V^*$. Let $L$ denote the line bundle over $\pl$ whose
fiber over an oriented line $l$ consists of linear functionals on
$l$.

We are going to remind the construction from
\cite{alesker-gafa-04} of a natural linear map
$$\Theta_{k,d}\colon \ond \otimes C^{\infty} ((\pl)^k, L^{\boxtimes
k})\to \pvd$$ which commutes with the natural action of the group
$GL(V)$ on both spaces and induces an epimorphism on the subspaces
of smooth vectors.

The construction is as follows. Let $\mu\in \ond,\, A_1,\dots
,A_k\in \ck(V)$. Then $\int_{\sum_{j=1}^k \lam_jA_j} \mu$ is a
polynomial in $\lam_j\geq 0$ of degree at most $n+d$. This can be
easily seen directly, but it was also proved in general for
polynomial valuations by Khovanskii and Pukhlikov
\cite{khovanskii-pukhlikov1}. Also it easily follows that the
coefficients of this polynomial depend continuously on $(A_1,\dots
,A_k)\in \ck(V)^k$ with respect to the Hausdorff metric. Hence we
can define a continuous map $\Theta_{k,d}'\colon\ond
\times\ck(V)^k\to PVal_d(V)$ given by
$$(\Theta_{k,d}'(\mu;A_1,\dots,A_k))(K):=\frac{\pt^k}{\pt
\lam_1\dots\pt\lam_k}\bigg |_{\lam_j =0}\int_{K+\sum_{j=1}^k
\lam_jA_j}\mu.$$

It is clear that $\Theta_{k,d}'$ is Minkowski additive with
respect to each $A_j$. Namely, say for $j=1, \, a,b\geq 0$, one
has
$$\Theta_{k,d}'(\mu;a A_1'+b
A_1'',A_2,\dots,A_k)=a\Theta_{k,d}'(\mu;A_1',A_2,\dots,A_k)+
b\Theta_{k,d}'(\mu;A_1'',A_2,\dots,A_k).$$

Remind that for any $A\in \ck(V)$ one defines the supporting
functional $h_A(y):=\sup_{x\in A}(y,x)$ for any $y\in V^*$. Thus
$h_A\in C(\PP_+(V^*),L)$. Moreover it is well known (and easy to
see) that $A_N\to A$ in the Hausdorff metric if and only if
$h_{A_N}\to h_A$ in $C(\PP_+(V^*),L)$. Also any section $F\in
C^2(\PP_+(V^*),L)$ can be presented as a difference $F=G-H$ where
$G,\,H\in C^2(\PP_+(V^*),L)$ are supporting functionals of some
convex compact sets and $\max\{||G||_{C^2},||H||_{C^2}\}\leq c
||F||_{C^2}$ where $c$ is a constant. (Indeed one can choose
$G=F+R\cdot h_D,\, H=R\cdot h_D$ where $D$ is the unit Euclidean
ball, and $R$ is a large enough constant depending on
$||F||_{C^2}$.) Hence we can uniquely extend $\Theta_{s,d}'$ to a
multilinear continuous map (which we will denote by the same
letter):
$$\Theta_{k,d}'\colon\ond \times (C^2(\PP_+(V^*),L))^k\to PVal_d(V).$$
By Theorem \ref{T:bilinear-forms} it follows that this map gives
rise to a continuous linear map
$$\Theta_{k,d}\colon\ond\otimes C^{\infty}(\PP_+(V^*)^k, L^{\boxtimes
k})\to PVal_d(V).$$ Since $\Theta_{k,d}$ commutes with the action
of $GL(V)$, its image is contained in $PVal_d^{sm}(V)$. Thus we
got a continuous map
$$\Theta_{k,d}\colon\ond\otimes C^{\infty}(\PP_+(V^*)^k, L^{\boxtimes
k})\to PVal_d^{sm}(V).$$  which we wanted to construct.

We will study this map $\Theta_{k,d}$. Note that it depends on $k$
and $d$ which will be fixed from now on. Let us denote by
$\Theta_d$ the sum of the maps $\bigoplus _{k=0}^{n}\Theta_{k,d}$.
Thus $$\Theta_d\colon\ond \otimes\left(\bigoplus
_{k=0}^{n}C^{\infty}(\pl^k, L^{\boxtimes k})\right) \to
PVal^{sm}_d(V).$$ The following result was proved by the author in
\cite{alesker-gafa-04}, Corollary 1.9.
\begin{lemma}[\cite{alesker-gafa-04}]\label{theta-onto}
The map $\Theta_d$ is onto $PVal^{sm}_d(V)$.
\end{lemma}
Since the source and the target spaces of $\Theta_d$ are Fr\'echet
spaces, by the open mapping theorem (see e.g. \cite{schaefer}, Ch.
III, \S 2) the topology on $PVal^{sm}_d(V)$ is the quotient
topology on $\ond \otimes\left(\bigoplus
_{k=0}^{n}C^{\infty}(\pl^k, L^{\boxtimes k})\right)$.

{\bf Proof} of Proposition \ref{P:ext-prod}. Denote $n:=\dim V,\,
m:=\dim W$. We have the following claim whose proof is easy and is
omitted.
\begin{claim}\label{C:ap1}
Let $\phi\in PVal_d(V)$. Let $\mu\in \Ome^m_{d'}(W)$. Then the map
$\Psi_{\phi,\mu}\colon \ck(V\times W)\to \CC$ given by

$$\Psi_{\phi,\mu}(K):=\int_{w\in W}\phi(K\cap(V\times \{w\})d\mu(w)$$
is a continuous valuation polynomial of degree at most $d+d'$.
\end{claim}
Hence by a result of Khovanskii and Pukhlikov
\cite{khovanskii-pukhlikov1},
$\Psi_{\phi,\mu}(\sum_{i=1}^s\lam_iK_i)$ is a polynomial in
$\lam_1,\dots,\lam_s\geq 0$ of degree at most $d+d'+n+m$  for any
$K_1,\dots,K_s\in \ck(V\times W)$ (for translation invariant
valuations this fact was proved earlier by McMullen
\cite{mcmullen-euler}). Hence for any $A_1,\dots,A_k\in \ck(W)$,
$K\in \ck(V\times W)$ the expression
\begin{eqnarray*}
\Psi_{\phi,\mu}\left(K+\left(\{0\}\times
\sum_{i=1}^k\lam_iA_i\right)\right)=\int_{w\in
W}\phi\left(\left(K+\left(\{0\}\times\sum_{i=1}^k\lam_iA_i\right)\right)\cap
\left(V\times\{w\}\right)\right)d\mu(w)
\end{eqnarray*}
is a polynomial in $\lam_1,\dots,\lam_k\geq 0$ of degree at most
$d+d'+n+m$ (in particular, there is a uniform bound on the
degree).

It easily follows that
$\frac{\pt^k}{\pt\lam_1\dots\pt\lam_k}\big|_0
\Psi_{\phi,\mu}\left(K+\left(\{0\}\times\sum_{i=1}^k\lam_iA_i\right)\right)$
is a continuous valuation with respect to $K\in \ck(V\times W)$.
Moreover the map
\begin{eqnarray}\label{ap2}
PVal_d(V)\times\Ome_{d'}^m(W)\times \ck(W)^k\to
PVal_{d+d'}(V\times W)
\end{eqnarray}
given by
$$(\phi,\mu;A_1,\dots,A_k)\mapsto \left[K\mapsto \frac{\pt^k}{\pt\lam_1\dots\pt\lam_k}\big|_0
\Psi_{\phi,\mu}\left(K+\left(\{0\}\times\sum_{i=1}^k\lam_iA_i\right)\right)\right]$$
is continuous. Also this map is Minkowski additive with respect to
each $A_j\in \ck(W)$. By the argument used in the construction of
the map $\Theta_d$, the map (\ref{ap2}) extends (uniquely) to a
multilinear {\itshape continuous} map
\begin{eqnarray}\label{ap3}
PVal_d(V)\times\Ome_{d'}^m(W)\times C^\infty(\PP_+(W^*),L)^k\to
PVal_{d+d'}(V\times W).
\end{eqnarray}
By Theorem \ref{T:bilinear-forms} the map (\ref{ap3}) gives rise
to a bilinear {\itshape continuous} map
\begin{eqnarray}\label{ap4}
PVal_d(V)\times\Ome_{d'}^m(W)\otimes
C^\infty(\PP_+(W^*)^k,L^{\boxtimes k})\to PVal_{d+d'}(V\times W).
\end{eqnarray}
Summing up over $k=0,\dots, m$ we obtain a bilinear {\itshape
continuous} map
\begin{eqnarray}\label{ap5}
PVal_d(V)\times\left(\Ome_{d'}^m(W)\otimes \left(\bigoplus_{k=0}^m
C^\infty\left(\PP_+(W^*)^k,L^{\boxtimes k}\right)\right)\right)\to
PVal_{d+d'}(V\times W).
\end{eqnarray}
\begin{lemma}\label{L:ap6}
The map (\ref{ap5}) factorizes (uniquely) as
\begin{eqnarray*}
\Vtriangle<1`1`-1;900>[PVal_d(V)\times\left(\Ome_{d'}^m(W)\otimes
\left(\bigoplus_{k=0}^m C^\infty\left(\PP_+(W^*)^k,L^{\boxtimes
k}\right)\right)\right)`PVal_{d+d'}(V\times W)`PVal_d(V)\times
PVal_{d'}^{sm}(W);`Id\boxtimes \Theta_d`].
\end{eqnarray*}
\end{lemma}
{\bf Proof.} If $PVal_d(V)$ is replaced by $PVal_d^{sm}(V)$, the
corresponding result was proved in \cite{alesker-gafa-04}, and the
obtained map $$PVal_d^{sm}(V)\times PVal_{d'}^{sm}(W)\to
PVal_{d+d'}(V\times W)$$ was exactly the exterior product. Our
lemma follows from this fact and the continuity of the map
(\ref{ap5}) because $PVal^{sm}_d(V)\subset PVal_d(V)$ is a dense
subspace. Lemma is proved. \qed

The map $$PVal_d(V)\times PVal^{sm}_{d'}(W)\to PVal_{d+d'}(V\times
W)$$ from Lemma \ref{L:ap6} is the map we need. Proposition
\ref{P:ext-prod} is proved. \qed

\end{document}

%% file: diagram.tex
\makeatletter

\def\diagram{\m@th\leftwidth=\z@ \rightwidth=\z@ \topheight=\z@
\botheight=\z@ \setbox\@picbox\hbox\bgroup}

\def\enddiagram{\egroup\wd\@picbox\rightwidth\unitlength
\ht\@picbox\topheight\unitlength \dp\@picbox\botheight\unitlength
\hskip\leftwidth\unitlength\box\@picbox}

\def\bfig{\begin{diagram}}
\def\efig{\end{diagram}}
\newcount\wideness \newcount\leftwidth \newcount\rightwidth
\newcount\highness \newcount\topheight \newcount\botheight

\def\ratchet#1#2{\ifnum#1<#2 \global #1=#2 \fi}

\def\putbox(#1,#2)#3{%
\horsize{\wideness}{#3} \divide\wideness by 2
{\advance\wideness by #1 \ratchet{\rightwidth}{\wideness}}
{\advance\wideness by -#1 \ratchet{\leftwidth}{\wideness}}
\vertsize{\highness}{#3} \divide\highness by 2
{\advance\highness by #2 \ratchet{\topheight}{\highness}}
{\advance\highness by -#2 \ratchet{\botheight}{\highness}}
\put(#1,#2){\makebox(0,0){$#3$}}}

\def\putlbox(#1,#2)#3{%
\horsize{\wideness}{#3}
{\advance\wideness by #1 \ratchet{\rightwidth}{\wideness}}
{\ratchet{\leftwidth}{-#1}}
\vertsize{\highness}{#3} \divide\highness by 2
{\advance\highness by #2 \ratchet{\topheight}{\highness}}
{\advance\highness by -#2 \ratchet{\botheight}{\highness}}
\put(#1,#2){\makebox(0,0)[l]{$#3$}}}

\def\putrbox(#1,#2)#3{%
\horsize{\wideness}{#3}
{\ratchet{\rightwidth}{#1}}
{\advance\wideness by -#1 \ratchet{\leftwidth}{\wideness}}
\vertsize{\highness}{#3} \divide\highness by 2
{\advance\highness by #2 \ratchet{\topheight}{\highness}}
{\advance\highness by -#2 \ratchet{\botheight}{\highness}}
\put(#1,#2){\makebox(0,0)[r]{$#3$}}}

\def\adjust[#1]{} 

\newcount \coefa
\newcount \coefb
\newcount \coefc
\newcount\tempcounta
\newcount\tempcountb
\newcount\tempcountc
\newcount\tempcountd
\newcount\xext
\newcount\yext
\newcount\xoff
\newcount\yoff
\newcount\gap%
\newcount\arrowtypea
\newcount\arrowtypeb
\newcount\arrowtypec
\newcount\arrowtyped
\newcount\arrowtypee
\newcount\height
\newcount\width
\newcount\xpos
\newcount\ypos
\newcount\run
\newcount\rise
\newcount\arrowlength
\newcount\halflength
\newcount\arrowtype
\newdimen\tempdimen
\newdimen\xlen
\newdimen\ylen
\newsavebox{\tempboxa}%
\newsavebox{\tempboxb}%
\newsavebox{\tempboxc}%

\newdimen\w@dth

\def\setw@dth#1#2{\setbox\z@\hbox{\m@th$#1$}\w@dth=\wd\z@
\setbox\@ne\hbox{\m@th$#2$}\ifnum\w@dth<\wd\@ne \w@dth=\wd\@ne \fi
\advance\w@dth by 1.2em}

\def\t@^#1_#2{\allowbreak\def\n@one{#1}\def\n@two{#2}\mathrel
{\setw@dth{#1}{#2}
\mathop{\hbox to \w@dth{\rightarrowfill}}\limits
\ifx\n@one\empty\else ^{\box\z@}\fi
\ifx\n@two\empty\else _{\box\@ne}\fi}}
\def\t@@^#1{\@ifnextchar_{\t@^{#1}}{\t@^{#1}_{}}}
\def\to{\@ifnextchar^{\t@@}{\t@@^{}}}

\def\t@left^#1_#2{\def\n@one{#1}\def\n@two{#2}\mathrel{\setw@dth{#1}{#2}
\mathop{\hbox to \w@dth{\leftarrowfill}}\limits
\ifx\n@one\empty\else ^{\box\z@}\fi
\ifx\n@two\empty\else _{\box\@ne}\fi}}
\def\t@@left^#1{\@ifnextchar_{\t@left^{#1}}{\t@left^{#1}_{}}}
\def\toleft{\@ifnextchar^{\t@@left}{\t@@left^{}}}

\def\two@^#1_#2{\allowbreak
\def\n@one{#1}\def\n@two{#2}\mathrel{\setw@dth{#1}{#2}
\mathop{\vcenter{\lineskip\z@\baselineskip\z@
                 \hbox to \w@dth{\rightarrowfill}%
                 \hbox to \w@dth{\rightarrowfill}}%
       }\limits
\ifx\n@one\empty\else ^{\box\z@}\fi
\ifx\n@two\empty\else _{\box\@ne}\fi}}
\def\tw@@^#1{\@ifnextchar _{\two@^{#1}}{\two@^{#1}_{}}}
\def\two{\@ifnextchar ^{\tw@@}{\tw@@^{}}}

\def\tofr@^#1_#2{\def\n@one{#1}\def\n@two{#2}\mathrel{\setw@dth{#1}{#2}
\mathop{\vcenter{\hbox to \w@dth{\rightarrowfill}\kern-1.7ex
                 \hbox to \w@dth{\leftarrowfill}}%
       }\limits
\ifx\n@one\empty\else ^{\box\z@}\fi
\ifx\n@two\empty\else _{\box\@ne}\fi}}
\def\t@fr@^#1{\@ifnextchar_ {\tofr@^{#1}}{\tofr@^{#1}_{}}}
\def\tofro{\@ifnextchar^ {\t@fr@}{\t@fr@^{}}}

\def\mon{\mathop{\m@th\hbox to
      14.6\P@{\lasyb\char'51\hskip-2.1\P@$\arrext$\hss
$\mathord\rightarrow$}}\limits} 
\def\leftmono{\mathrel{\m@th\hbox to
14.6\P@{$\mathord\leftarrow$\hss$\arrext$\hskip-2.1\P@\lasyb\char'50%
}}\limits} 
\mathchardef\arrext="0200       

\def\settypes(#1,#2,#3){\arrowtypea#1 \arrowtypeb#2 \arrowtypec#3}
\def\settoheight#1#2{\setbox\@tempboxa\hbox{#2}#1\ht\@tempboxa\relax}%
\def\settodepth#1#2{\setbox\@tempboxa\hbox{#2}#1\dp\@tempboxa\relax}%
\def\settokens`#1`#2`#3`#4`{%
     \def\tokena{#1}\def\tokenb{#2}\def\tokenc{#3}\def\tokend{#4}}
\def\setsqparms[#1`#2`#3`#4;#5`#6]{%
\arrowtypea #1
\arrowtypeb #2
\arrowtypec #3
\arrowtyped #4
\width #5
\height #6
}
\def\setpos(#1,#2){\xpos=#1 \ypos#2}

\def\settriparms[#1`#2`#3;#4]{\settripairparms[#1`#2`#3`1`1;#4]}%

\def\settripairparms[#1`#2`#3`#4`#5;#6]{%
\arrowtypea #1
\arrowtypeb #2
\arrowtypec #3
\arrowtyped #4
\arrowtypee #5
\width #6
\height #6
}

\def\resetparms{\settripairparms[1`1`1`1`1;500]\width 500}

\resetparms

\def\mvector(#1,#2)#3{
\put(0,0){\vector(#1,#2){#3}}%
\put(0,0){\vector(#1,#2){26}}%
}
\def\evector(#1,#2)#3{{
\arrowlength #3
\put(0,0){\vector(#1,#2){\arrowlength}}%
\advance \arrowlength by-30
\put(0,0){\vector(#1,#2){\arrowlength}}%
}}

\def\horsize#1#2{%
\settowidth{\tempdimen}{$#2$}%
#1=\tempdimen
\divide #1 by\unitlength
}

\def\vertsize#1#2{%
\settoheight{\tempdimen}{$#2$}%
#1=\tempdimen
\settodepth{\tempdimen}{$#2$}%
\advance #1 by\tempdimen
\divide #1 by\unitlength
}

\def\putvector(#1,#2)(#3,#4)#5#6{{%
\ifnum3<\arrowtype
\putdashvector(#1,#2)(#3,#4)#5\arrowtype
\else
\ifnum\arrowtype<-3
\putdashvector(#1,#2)(#3,#4)#5\arrowtype
\else
\xpos=#1
\ypos=#2
\run=#3
\rise=#4
\arrowlength=#5
\ifnum \arrowtype<0
    \ifnum \run=0
        \advance \ypos by-\arrowlength
    \else
        \tempcounta \arrowlength
        \multiply \tempcounta by\rise
        \divide \tempcounta by\run
        \ifnum\run>0
            \advance \xpos by\arrowlength
            \advance \ypos by\tempcounta
        \else
            \advance \xpos by-\arrowlength
            \advance \ypos by-\tempcounta
        \fi
    \fi
    \multiply \arrowtype by-1
    \multiply \rise by-1
    \multiply \run by-1
\fi
\ifcase \arrowtype
\or \put(\xpos,\ypos){\vector(\run,\rise){\arrowlength}}%
\or \put(\xpos,\ypos){\mvector(\run,\rise)\arrowlength}%
\or \put(\xpos,\ypos){\evector(\run,\rise){\arrowlength}}%
\fi\fi\fi
}}

\def\putsplitvector(#1,#2)#3#4{
\xpos #1
\ypos #2
\arrowtype #4
\halflength #3
\arrowlength #3
\gap 140
\advance \halflength by-\gap
\divide \halflength by2
\ifnum\arrowtype>0
   \ifcase \arrowtype
   \or \put(\xpos,\ypos){\line(0,-1){\halflength}}%
       \advance\ypos by-\halflength
       \advance\ypos by-\gap
       \put(\xpos,\ypos){\vector(0,-1){\halflength}}%
   \or \put(\xpos,\ypos){\line(0,-1)\halflength}%
       \put(\xpos,\ypos){\vector(0,-1)3}%
       \advance\ypos by-\halflength
       \advance\ypos by-\gap
       \put(\xpos,\ypos){\vector(0,-1){\halflength}}%
   \or \put(\xpos,\ypos){\line(0,-1)\halflength}%
       \advance\ypos by-\halflength
       \advance\ypos by-\gap
       \put(\xpos,\ypos){\evector(0,-1){\halflength}}%
   \fi
\else \arrowtype=-\arrowtype
   \ifcase\arrowtype
   \or \advance \ypos by-\arrowlength
       \put(\xpos,\ypos){\line(0,1){\halflength}}%
       \advance\ypos by\halflength
       \advance\ypos by\gap
       \put(\xpos,\ypos){\vector(0,1){\halflength}}%
   \or \advance \ypos by-\arrowlength
       \put(\xpos,\ypos){\line(0,1)\halflength}%
       \put(\xpos,\ypos){\vector(0,1)3}%
       \advance\ypos by\halflength
       \advance\ypos by\gap
       \put(\xpos,\ypos){\vector(0,1){\halflength}}%
   \or \advance \ypos by-\arrowlength
       \put(\xpos,\ypos){\line(0,1)\halflength}%
       \advance\ypos by\halflength
       \advance\ypos by\gap
       \put(\xpos,\ypos){\evector(0,1){\halflength}}%
   \fi
\fi
}

\def\putmorphism(#1)(#2,#3)[#4`#5`#6]#7#8#9{{%
\run #2
\rise #3
\ifnum\rise=0
  \puthmorphism(#1)[#4`#5`#6]{#7}{#8}#9%
\else\ifnum\run=0
  \putvmorphism(#1)[#4`#5`#6]{#7}{#8}#9%
\else
\setpos(#1)%
\arrowlength #7
\arrowtype #8
\ifnum\run=0
\else\ifnum\rise=0
\else
\ifnum\run>0
    \coefa=1
\else
   \coefa=-1
\fi
\ifnum\arrowtype>0
   \coefb=0
   \coefc=-1
\else
   \coefb=\coefa
   \coefc=1
   \arrowtype=-\arrowtype
\fi
\width=2
\multiply \width by\run
\divide \width by\rise
\ifnum \width<0  \width=-\width\fi
\advance\width by60
\if l#9 \width=-\width\fi
\putbox(\xpos,\ypos){#4}
{\multiply \coefa by\arrowlength
\advance\xpos by\coefa
\multiply \coefa by\rise
\divide \coefa by\run
\advance \ypos by\coefa
\putbox(\xpos,\ypos){#5} }%
{\multiply \coefa by\arrowlength
\divide \coefa by2
\advance \xpos by\coefa
\advance \xpos by\width
\multiply \coefa by\rise
\divide \coefa by\run
\advance \ypos by\coefa
\if l#9%
   \putrbox(\xpos,\ypos){#6}%
\else\if r#9%
   \putlbox(\xpos,\ypos){#6}%
\fi\fi }%
{\multiply \rise by-\coefc
\multiply \run by-\coefc
\multiply \coefb by\arrowlength
\advance \xpos by\coefb
\multiply \coefb by\rise
\divide \coefb by\run
\advance \ypos by\coefb
\multiply \coefc by70
\advance \ypos by\coefc
\multiply \coefc by\run
\divide \coefc by\rise
\advance \xpos by\coefc
\multiply \coefa by140
\multiply \coefa by\run
\divide \coefa by\rise
\advance \arrowlength by\coefa
\ifcase\arrowtype
\or \put(\xpos,\ypos){\vector(\run,\rise){\arrowlength}}%
\or \put(\xpos,\ypos){\mvector(\run,\rise){\arrowlength}}%
\or \put(\xpos,\ypos){\evector(\run,\rise){\arrowlength}}%
\fi}\fi\fi\fi\fi}}

\newcount\numbdashes \newcount\lengthdash \newcount\increment

\def\howmanydashes{
\numbdashes=\arrowlength \lengthdash=40
\divide\numbdashes by \lengthdash
\lengthdash=\arrowlength
\divide\lengthdash by \numbdashes
\increment=\lengthdash
\multiply\lengthdash by 3
\divide\lengthdash by 5
}

\def\putdashvector(#1)(#2,#3)#4#5{%
\ifnum#3=0 \putdashhvector(#1){#4}#5
\else
\ifnum#2=0
\putdashvvector(#1){#4}#5\fi\fi}

\def\putdashhvector(#1,#2)#3#4{{%
\arrowlength=#3 \howmanydashes
\multiput(#1,#2)(\increment,0){\numbdashes}%
{\vrule height .4pt width \lengthdash\unitlength}
\arrowtype=#4 \xpos=#1
\ifnum\arrowtype<0 \advance\arrowtype by 7 \fi
\ifcase\arrowtype
\or \advance\xpos by 10
    \put(\xpos,#2){\vector(-1,0){\lengthdash}}
    \advance\xpos by 40
    \put(\xpos,#2){\vector(-1,0){\lengthdash}}
\or \advance \xpos by 10
    \put(\xpos,#2){\vector(-1,0){\lengthdash}}
    \advance\xpos by  \arrowlength
    \advance\xpos by  -50
    \put(\xpos,#2){\vector(-1,0){\lengthdash}}
\or \advance\xpos by 10
    \put(\xpos,#2){\vector(-1,0){\lengthdash}}
\or \advance\xpos by \arrowlength
    \advance\xpos by -\lengthdash
    \put(\xpos,#2){\vector(1,0){\lengthdash}}
\or {\advance\xpos by 10
    \put(\xpos,#2){\vector(1,0){\lengthdash}}}
    \advance\xpos by \arrowlength
    \advance\xpos by -\lengthdash
    \put(\xpos,#2){\vector(1,0){\lengthdash}}
\or \advance\xpos by \arrowlength
    \advance\xpos by -\lengthdash
    \put(\xpos,#2){\vector(1,0){\lengthdash}}
    \advance\xpos by -40
    \put(\xpos,#2){\vector(1,0){\lengthdash}}
   \fi
}}

\def\putdashvvector(#1,#2)#3#4{{%
\arrowlength=#3 \howmanydashes
\ypos=#2 \advance\ypos by -\arrowlength
\multiput(#1,#2)(0,\increment){\numbdashes}%
    {\vrule width .4pt height \lengthdash\unitlength}
\arrowtype=#4 \ypos=#2
\ifnum\arrowtype<0 \advance\arrowtype by 7 \fi
\ifcase\arrowtype
\or \advance\ypos by \arrowlength \advance\ypos by -40
    \put(#1,\ypos){\vector(0,1){\lengthdash}}
    \advance\ypos by -40
    \put(#1,\ypos){\vector(0,1){\lengthdash}}
\or \advance\ypos by 10
    \put(#1,\ypos){\vector(0,1){\lengthdash}}
    \advance\ypos by \arrowlength \advance\ypos by -40
    \put(#1,\ypos){\vector(0,1){\lengthdash}}
\or \advance\ypos by \arrowlength \advance\ypos by -40
    \put(#1,\ypos){\vector(0,1){\lengthdash}}
\or \advance\ypos by 10
    \put(#1,\ypos){\vector(0,-1){\lengthdash}}
\or \advance\ypos by 10
    \put(#1,\ypos){\vector(0,-1){\lengthdash}}
    \advance\ypos by \arrowlength \advance\ypos by -40
    \put(#1,\ypos){\vector(0,-1){\lengthdash}}
\or \advance\ypos by 10
    \put(#1,\ypos){\vector(0,-1){\lengthdash}}
    \advance\ypos by 40
    \put(#1,\ypos){\vector(0,-1){\lengthdash}}
\fi
}}

\def\puthmorphism(#1,#2)[#3`#4`#5]#6#7#8{{%
\xpos #1
\ypos #2
\width #6
\arrowlength #6
\arrowtype=#7
\putbox(\xpos,\ypos){#3\vphantom{#4}}%
{\advance \xpos by\arrowlength
\putbox(\xpos,\ypos){\vphantom{#3}#4}}%
\horsize{\tempcounta}{#3}%
\horsize{\tempcountb}{#4}%
\divide \tempcounta by2
\divide \tempcountb by2
\advance \tempcounta by30
\advance \tempcountb by30
\advance \xpos by\tempcounta
\advance \arrowlength by-\tempcounta
\advance \arrowlength by-\tempcountb
\putvector(\xpos,\ypos)(1,0)\arrowlength\arrowtype
\divide \arrowlength by2
\advance \xpos by\arrowlength
\vertsize{\tempcounta}{#5}%
\divide\tempcounta by2
\advance \tempcounta by20
\if a#8 %
   \advance \ypos by\tempcounta
   \putbox(\xpos,\ypos){#5}%
\else
   \advance \ypos by-\tempcounta
   \putbox(\xpos,\ypos){#5}%
\fi}}

\def\putvmorphism(#1,#2)[#3`#4`#5]#6#7#8{{%
\xpos #1
\ypos #2
\arrowlength #6
\arrowtype #7
\settowidth{\xlen}{$#5$}%
\putbox(\xpos,\ypos){#3}%
{\advance \ypos by-\arrowlength
\putbox(\xpos,\ypos){#4}}%
{\advance\arrowlength by-140
\advance \ypos by-70
\ifdim\xlen>0pt
   \if m#8%
      \putsplitvector(\xpos,\ypos)\arrowlength\arrowtype
   \else
   \putvector(\xpos,\ypos)(0,-1)\arrowlength\arrowtype
   \fi
\else
   \putvector(\xpos,\ypos)(0,-1)\arrowlength\arrowtype
\fi}%
\ifdim\xlen>0pt
   \divide \arrowlength by2
   \advance\ypos by-\arrowlength
   \if l#8%
      \advance \xpos by-40
      \putrbox(\xpos,\ypos){#5}%
   \else\if r#8%
      \advance \xpos by40
      \putlbox(\xpos,\ypos){#5}%
   \else
      \putbox(\xpos,\ypos){#5}%
   \fi\fi
\fi
}}

\def\putsquarep<#1>(#2)[#3;#4`#5`#6`#7]{{%
\setsqparms[#1]%
\setpos(#2)%
\settokens`#3`%
\puthmorphism(\xpos,\ypos)[\tokenc`\tokend`{#7}]{\width}{\arrowtyped}b%
\advance\ypos by \height
\puthmorphism(\xpos,\ypos)[\tokena`\tokenb`{#4}]{\width}{\arrowtypea}a%
\putvmorphism(\xpos,\ypos)[``{#5}]{\height}{\arrowtypeb}l%
\advance\xpos by \width
\putvmorphism(\xpos,\ypos)[``{#6}]{\height}{\arrowtypec}r%
}}

\def\putsquare{\@ifnextchar <{\putsquarep}{\putsquarep%
   <\arrowtypea`\arrowtypeb`\arrowtypec`\arrowtyped;\width`\height>}}
\def\square{\@ifnextchar< {\squarep}{\squarep
   <\arrowtypea`\arrowtypeb`\arrowtypec`\arrowtyped;\width`\height>}}
\def\squarep<#1>[#2`#3`#4`#5;#6`#7`#8`#9]{{
\setsqparms[#1]
\diagram
\putsquarep<\arrowtypea`\arrowtypeb`\arrowtypec`
\arrowtyped;\width`\height>
(0,0)[#2`#3`#4`{#5};#6`#7`#8`{#9}]
\enddiagram
}}                                                 
\def\putptrianglep<#1>(#2,#3)[#4`#5`#6;#7`#8`#9]{{%
\settriparms[#1]%
\xpos=#2 \ypos=#3
\advance\ypos by \height
\puthmorphism(\xpos,\ypos)[#4`#5`{#7}]{\height}{\arrowtypea}a%
\putvmorphism(\xpos,\ypos)[`#6`{#8}]{\height}{\arrowtypeb}l%
\advance\xpos by\height
\putmorphism(\xpos,\ypos)(-1,-1)[``{#9}]{\height}{\arrowtypec}r%
}}

\def\putptriangle{\@ifnextchar <{\putptrianglep}{\putptrianglep
   <\arrowtypea`\arrowtypeb`\arrowtypec;\height>}}
\def\ptriangle{\@ifnextchar <{\ptrianglep}{\ptrianglep
   <\arrowtypea`\arrowtypeb`\arrowtypec;\height>}}
\def\ptrianglep<#1>[#2`#3`#4;#5`#6`#7]{{
\settriparms[#1]
\diagram
\putptrianglep<\arrowtypea`\arrowtypeb`
\arrowtypec;\height>
(0,0)[#2`#3`#4;#5`#6`{#7}]
\enddiagram
}}                                            

\def\putqtrianglep<#1>(#2,#3)[#4`#5`#6;#7`#8`#9]{{%
\settriparms[#1]%
\xpos=#2 \ypos=#3
\advance\ypos by\height
\puthmorphism(\xpos,\ypos)[#4`#5`{#7}]{\height}{\arrowtypea}a%
\putmorphism(\xpos,\ypos)(1,-1)[``{#8}]{\height}{\arrowtypeb}l%
\advance\xpos by\height
\putvmorphism(\xpos,\ypos)[`#6`{#9}]{\height}{\arrowtypec}r%
}}

\def\putqtriangle{\@ifnextchar <{\putqtrianglep}{\putqtrianglep
   <\arrowtypea`\arrowtypeb`\arrowtypec;\height>}}
\def\qtriangle{\@ifnextchar <{\qtrianglep}{\qtrianglep
   <\arrowtypea`\arrowtypeb`\arrowtypec;\height>}}
\def\qtrianglep<#1>[#2`#3`#4;#5`#6`#7]{{
\settriparms[#1]
\width=\height                                
\diagram
\putqtrianglep<\arrowtypea`\arrowtypeb`
\arrowtypec;\height>
(0,0)[#2`#3`#4;#5`#6`{#7}]
\enddiagram
}}

\def\putdtrianglep<#1>(#2,#3)[#4`#5`#6;#7`#8`#9]{{%
\settriparms[#1]%
\xpos=#2 \ypos=#3
\puthmorphism(\xpos,\ypos)[#5`#6`{#9}]{\height}{\arrowtypec}b%
\advance\xpos by \height \advance\ypos by\height
\putmorphism(\xpos,\ypos)(-1,-1)[``{#7}]{\height}{\arrowtypea}l%
\putvmorphism(\xpos,\ypos)[#4``{#8}]{\height}{\arrowtypeb}r%
}}

\def\putdtriangle{\@ifnextchar <{\putdtrianglep}{\putdtrianglep
   <\arrowtypea`\arrowtypeb`\arrowtypec;\height>}}
\def\dtriangle{\@ifnextchar <{\dtrianglep}{\dtrianglep
   <\arrowtypea`\arrowtypeb`\arrowtypec;\height>}}
\def\dtrianglep<#1>[#2`#3`#4;#5`#6`#7]{{
\settriparms[#1]
\width=\height                                
\diagram
\putdtrianglep<\arrowtypea`\arrowtypeb`
\arrowtypec;\height>
(0,0)[#2`#3`#4;#5`#6`{#7}]
\enddiagram
}}

\def\putbtrianglep<#1>(#2,#3)[#4`#5`#6;#7`#8`#9]{{%
\settriparms[#1]%
\xpos=#2 \ypos=#3
\puthmorphism(\xpos,\ypos)[#5`#6`{#9}]{\height}{\arrowtypec}b%
\advance\ypos by\height
\putmorphism(\xpos,\ypos)(1,-1)[``{#8}]{\height}{\arrowtypeb}r%
\putvmorphism(\xpos,\ypos)[#4``{#7}]{\height}{\arrowtypea}l%
}}

\def\putbtriangle{\@ifnextchar <{\putbtrianglep}{\putbtrianglep
   <\arrowtypea`\arrowtypeb`\arrowtypec;\height>}}
\def\btriangle{\@ifnextchar <{\btrianglep}{\btrianglep
   <\arrowtypea`\arrowtypeb`\arrowtypec;\height>}}
\def\btrianglep<#1>[#2`#3`#4;#5`#6`#7]{{
\settriparms[#1]
\width=\height                               
\diagram
\putbtrianglep<\arrowtypea`\arrowtypeb`
\arrowtypec;\height>
(0,0)[#2`#3`#4;#5`#6`{#7}]
\enddiagram
}}

\def\putAtrianglep<#1>(#2,#3)[#4`#5`#6;#7`#8`#9]{{%
\settriparms[#1]%
\xpos=#2 \ypos=#3
{\multiply \height by2
\puthmorphism(\xpos,\ypos)[#5`#6`{#9}]{\height}{\arrowtypec}b}%
\advance\xpos by\height \advance\ypos by\height
\putmorphism(\xpos,\ypos)(-1,-1)[#4``{#7}]{\height}{\arrowtypea}l%
\putmorphism(\xpos,\ypos)(1,-1)[``{#8}]{\height}{\arrowtypeb}r%
}}

\def\putAtriangle{\@ifnextchar <{\putAtrianglep}{\putAtrianglep
   <\arrowtypea`\arrowtypeb`\arrowtypec;\height>}}
\def\Atriangle{\@ifnextchar <{\Atrianglep}{\Atrianglep
   <\arrowtypea`\arrowtypeb`\arrowtypec;\height>}}
\def\Atrianglep<#1>[#2`#3`#4;#5`#6`#7]{{
\settriparms[#1]
\width=\height                                     
\diagram
\putAtrianglep<\arrowtypea`\arrowtypeb`
\arrowtypec;\height>
(0,0)[#2`#3`#4;#5`#6`{#7}]
\enddiagram
}}

\def\putAtrianglepairp<#1>(#2)[#3;#4`#5`#6`#7`#8]{{%
\settripairparms[#1]%
\setpos(#2)%
\settokens`#3`%
\puthmorphism(\xpos,\ypos)[\tokenb`\tokenc`{#7}]{\height}{\arrowtyped}b%
\advance\xpos by\height
\puthmorphism(\xpos,\ypos)[\phantom{\tokenc}`\tokend`{#8}]%
{\height}{\arrowtypee}b%
\advance\ypos by\height
\putmorphism(\xpos,\ypos)(-1,-1)[\tokena``{#4}]{\height}{\arrowtypea}l%
\putvmorphism(\xpos,\ypos)[``{#5}]{\height}{\arrowtypeb}m%
\putmorphism(\xpos,\ypos)(1,-1)[``{#6}]{\height}{\arrowtypec}r%
}}

\def\putAtrianglepair{\@ifnextchar <{\putAtrianglepairp}{\putAtrianglepairp%
   <\arrowtypea`\arrowtypeb`\arrowtypec`\arrowtyped`\arrowtypee;\height>}}
\def\Atrianglepair{\@ifnextchar <{\Atrianglepairp}{\Atrianglepairp%
   <\arrowtypea`\arrowtypeb`\arrowtypec`\arrowtyped`\arrowtypee;\height>}}

\def\Atrianglepairp<#1>[#2;#3`#4`#5`#6`#7]{{
\settripairparms[#1]
\settokens`#2`
\width=\height                                
\diagram
\putAtrianglepairp                            
<\arrowtypea`\arrowtypeb`\arrowtypec`
\arrowtyped`\arrowtypee;\height>
(0,0)[{#2};#3`#4`#5`#6`{#7}]
\enddiagram
}}

\def\putVtrianglep<#1>(#2,#3)[#4`#5`#6;#7`#8`#9]{{%
\settriparms[#1]%
\xpos=#2 \ypos=#3
\advance\ypos by\height
{\multiply\height by2
\puthmorphism(\xpos,\ypos)[#4`#5`{#7}]{\height}{\arrowtypea}a}%
\putmorphism(\xpos,\ypos)(1,-1)[`#6`{#8}]{\height}{\arrowtypeb}l%
\advance\xpos by\height
\advance\xpos by\height
\putmorphism(\xpos,\ypos)(-1,-1)[``{#9}]{\height}{\arrowtypec}r%
}}

\def\putVtriangle{\@ifnextchar <{\putVtrianglep}{\putVtrianglep
   <\arrowtypea`\arrowtypeb`\arrowtypec;\height>}}
\def\Vtriangle{\@ifnextchar <{\Vtrianglep}{\Vtrianglep
   <\arrowtypea`\arrowtypeb`\arrowtypec;\height>}}
\def\Vtrianglep<#1>[#2`#3`#4;#5`#6`#7]{{
\settriparms[#1]
\width=\height                                 
\diagram
\putVtrianglep<\arrowtypea`\arrowtypeb`
\arrowtypec;\height>
(0,0)[#2`#3`#4;#5`#6`{#7}]
\enddiagram
}}

\def\putVtrianglepairp<#1>(#2)[#3;#4`#5`#6`#7`#8]{{
\settripairparms[#1]%
\setpos(#2)%
\settokens`#3`%
\advance\ypos by\height
\putmorphism(\xpos,\ypos)(1,-1)[`\tokend`{#6}]{\height}{\arrowtypec}l%
\puthmorphism(\xpos,\ypos)[\tokena`\tokenb`{#4}]{\height}{\arrowtypea}a%
\advance\xpos by\height
\puthmorphism(\xpos,\ypos)[\phantom{\tokenb}`\tokenc`{#5}]%
{\height}{\arrowtypeb}a%
\putvmorphism(\xpos,\ypos)[``{#7}]{\height}{\arrowtyped}m%
\advance\xpos by\height
\putmorphism(\xpos,\ypos)(-1,-1)[``{#8}]{\height}{\arrowtypee}r%
}}

\def\putVtrianglepair{\@ifnextchar <{\putVtrianglepairp}{\putVtrianglepairp%
    <\arrowtypea`\arrowtypeb`\arrowtypec`\arrowtyped`\arrowtypee;\height>}}
\def\Vtrianglepair{\@ifnextchar <{\Vtrianglepairp}{\Vtrianglepairp%
    <\arrowtypea`\arrowtypeb`\arrowtypec`\arrowtyped`\arrowtypee;\height>}}
\def\Vtrianglepairp<#1>[#2;#3`#4`#5`#6`#7]{{
\settripairparms[#1]
\settokens`#2`
\diagram
\putVtrianglepairp                             
<\arrowtypea`\arrowtypeb`\arrowtypec`
\arrowtyped`\arrowtypee;\height>
(0,0)[{#2};#3`#4`#5`#6`{#7}]
\enddiagram
}}

\def\putCtrianglep<#1>(#2,#3)[#4`#5`#6;#7`#8`#9]{{%
\settriparms[#1]%
\xpos=#2 \ypos=#3
\advance\ypos by\height
\putmorphism(\xpos,\ypos)(1,-1)[``{#9}]{\height}{\arrowtypec}l%
\advance\xpos by\height
\advance\ypos by\height
\putmorphism(\xpos,\ypos)(-1,-1)[#4`#5`{#7}]{\height}{\arrowtypea}l%
{\multiply\height by 2
\putvmorphism(\xpos,\ypos)[`#6`{#8}]{\height}{\arrowtypeb}r}%
}}

\def\putCtriangle{\@ifnextchar <{\putCtrianglep}{\putCtrianglep
    <\arrowtypea`\arrowtypeb`\arrowtypec;\height>}}
\def\Ctriangle{\@ifnextchar <{\Ctrianglep}{\Ctrianglep
    <\arrowtypea`\arrowtypeb`\arrowtypec;\height>}}
\def\Ctrianglep<#1>[#2`#3`#4;#5`#6`#7]{{
\settriparms[#1]
\width=\height                               
\diagram
\putCtrianglep<\arrowtypea`\arrowtypeb`
\arrowtypec;\height>
(0,0)[#2`#3`#4;#5`#6`{#7}]
\enddiagram
}}                                           
\def\putDtrianglep<#1>(#2,#3)[#4`#5`#6;#7`#8`#9]{{%
\settriparms[#1]%
\xpos=#2 \ypos=#3
\advance\xpos by\height \advance\ypos by\height
\putmorphism(\xpos,\ypos)(-1,-1)[``{#9}]{\height}{\arrowtypec}r%
\advance\xpos by-\height \advance\ypos by\height
\putmorphism(\xpos,\ypos)(1,-1)[`#5`{#8}]{\height}{\arrowtypeb}r%
{\multiply\height by 2
\putvmorphism(\xpos,\ypos)[#4`#6`{#7}]{\height}{\arrowtypea}l}%
}}

\def\putDtriangle{\@ifnextchar <{\putDtrianglep}{\putDtrianglep
    <\arrowtypea`\arrowtypeb`\arrowtypec;\height>}}
\def\Dtriangle{\@ifnextchar <{\Dtrianglep}{\Dtrianglep
   <\arrowtypea`\arrowtypeb`\arrowtypec;\height>}}
\def\Dtrianglep<#1>[#2`#3`#4;#5`#6`#7]{{
\settriparms[#1]
\width=\height                              
\diagram
\putDtrianglep<\arrowtypea`\arrowtypeb`
\arrowtypec;\height>
(0,0)[#2`#3`#4;#5`#6`{#7}]
\enddiagram
}}                                          
\def\setrecparms[#1`#2]{\width=#1 \height=#2}%

\def\recursep<#1`#2>[#3;#4`#5`#6`#7`#8]{{\m@th
\width=#1 \height=#2
\settokens`#3`
\settowidth{\tempdimen}{$\tokena$}
\ifdim\tempdimen=0pt
  \savebox{\tempboxa}{\hbox{$\tokenb$}}%
  \savebox{\tempboxb}{\hbox{$\tokend$}}%
  \savebox{\tempboxc}{\hbox{$#6$}}%
\else
  \savebox{\tempboxa}{\hbox{$\hbox{$\tokena$}\times\hbox{$\tokenb$}$}}%
  \savebox{\tempboxb}{\hbox{$\hbox{$\tokena$}\times\hbox{$\tokend$}$}}%
  \savebox{\tempboxc}{\hbox{$\hbox{$\tokena$}\times\hbox{$#6$}$}}%
\fi
\ypos=\height
\divide\ypos by 2
\xpos=\ypos
\advance\xpos by \width
\bfig
\putCtrianglep<-1`1`1;\ypos>(0,0)[`\tokenc`;#5`#6`{#7}]%
\puthmorphism(\ypos,0)[\tokend`\usebox{\tempboxb}`{#8}]{\width}{-1}b%
\puthmorphism(\ypos,\height)[\tokenb`\usebox{\tempboxa}`{#4}]{\width}{-1}a%
\advance\ypos by \width
\putvmorphism(\ypos,\height)[``\usebox{\tempboxc}]{\height}1r%
\efig
}}

\def\recurse{\@ifnextchar <{\recursep}{\recursep<\width`\height>}}

\def\puttwohmorphisms(#1,#2)[#3`#4;#5`#6]#7#8#9{{%
\puthmorphism(#1,#2)[#3`#4`]{#7}0a
\ypos=#2
\advance\ypos by 20
\puthmorphism(#1,\ypos)[\phantom{#3}`\phantom{#4}`#5]{#7}{#8}a
\advance\ypos by -40
\puthmorphism(#1,\ypos)[\phantom{#3}`\phantom{#4}`#6]{#7}{#9}b
}}

\def\puttwovmorphisms(#1,#2)[#3`#4;#5`#6]#7#8#9{{%
\putvmorphism(#1,#2)[#3`#4`]{#7}0a
\xpos=#1
\advance\xpos by -20
\putvmorphism(\xpos,#2)[\phantom{#3}`\phantom{#4}`#5]{#7}{#8}l
\advance\xpos by 40
\putvmorphism(\xpos,#2)[\phantom{#3}`\phantom{#4}`#6]{#7}{#9}r
}}

\def\puthcoequalizer(#1)[#2`#3`#4;#5`#6`#7]#8#9{{%
\setpos(#1)%
\puttwohmorphisms(\xpos,\ypos)[#2`#3;#5`#6]{#8}11%
\advance\xpos by #8
\puthmorphism(\xpos,\ypos)[\phantom{#3}`#4`#7]{#8}1{#9}
}}

\def\putvcoequalizer(#1)[#2`#3`#4;#5`#6`#7]#8#9{{%
\setpos(#1)%
\puttwovmorphisms(\xpos,\ypos)[#2`#3;#5`#6]{#8}11%
\advance\ypos by -#8
\putvmorphism(\xpos,\ypos)[\phantom{#3}`#4`#7]{#8}1{#9}
}}

\def\putthreehmorphisms(#1)[#2`#3;#4`#5`#6]#7(#8)#9{{%
\setpos(#1) \settypes(#8)
\if a#9 %
     \vertsize{\tempcounta}{#5}%
     \vertsize{\tempcountb}{#6}%
     \ifnum \tempcounta<\tempcountb \tempcounta=\tempcountb \fi
\else
     \vertsize{\tempcounta}{#4}%
     \vertsize{\tempcountb}{#5}%
     \ifnum \tempcounta<\tempcountb \tempcounta=\tempcountb \fi
\fi
\advance \tempcounta by 60
\puthmorphism(\xpos,\ypos)[#2`#3`#5]{#7}{\arrowtypeb}{#9}
\advance\ypos by \tempcounta
\puthmorphism(\xpos,\ypos)[\phantom{#2}`\phantom{#3}`#4]{#7}{\arrowtypea}{#9}
\advance\ypos by -\tempcounta \advance\ypos by -\tempcounta
\puthmorphism(\xpos,\ypos)[\phantom{#2}`\phantom{#3}`#6]{#7}{\arrowtypec}{#9}
}}

\def\setarrowtoks[#1`#2`#3`#4`#5`#6]{%
\def\toka{#1}
\def\tokb{#2}
\def\tokc{#3}
\def\tokd{#4}
\def\toke{#5}
\def\tokf{#6}
}
\def\hex{\@ifnextchar <{\hexp}{\hexp<1000`400>}}
\def\hexp<#1`#2>[#3`#4`#5`#6`#7`#8;#9]{%
\setarrowtoks[#9]
\yext=#2 \advance \yext by #2
\xext=#1 \advance\xext by \yext
\bfig
\putCtriangle<-1`0`1;#2>(0,0)[`#5`;\tokb``\tokd]
\xext=#1 \yext=#2 \advance \yext by #2
\putsquare<1`0`0`1;\xext`\yext>(#2,0)[#3`#4`#7`#8;\toka```\tokf]
\advance \xext by #2
\putDtriangle<0`1`-1;#2>(\xext,0)[`#6`;`\tokc`\toke]
\efig
}
\makeatother